\documentclass[12pt]{article}%

\usepackage{amsmath}
\usepackage{graphicx}
\usepackage{amsfonts}
\usepackage{amsthm}
\usepackage{amssymb}%
\usepackage{xcolor}
\usepackage[subnum]{cases}
\usepackage{cite}
\usepackage{titling}
\usepackage{authblk}
\usepackage[hidelinks]{hyperref}

\author{M.S. Indulekha\\ \vspace{0.5 cm} Email: \href{mailto:indulekha@cmu.edu}{indulekha@cmu.edu}}
\title{Motion of Elastic Thin Films by Evaporation-Condensation in the Dewetting Regime}
\date{}
\affil{Department of Mathematical Sciences, Carnegie Mellon University, 5000 Forbes Avenue, Pittsburgh, 15217, PA, USA.}

\def\a{\mathcal{A}}

\def\d{\mathcal{D}}
\def\e{\mathscr{E}}
\def\f{\mathcal{F}}
\def\g{\mathcal{G}}
\def\h{\mathcal{H}}

\def\e{\mathcal{E}}

\def\s{\mathcal{S}}

\def\t{\mathcal{T}}
\def\l{\mathcal{L}}
\def\w{\mathcal{W}}

\def\A{\mathbb{A}}

\def\C{\mathbb{C}}

\def\L{\mathcal{L}}

\def\N{\mathbb{N}}

\def\R{\mathbb{R}}

\def\Z{\mathbb{Z}}
\def\htil{\Tilde{h}}
\def\util{\tilde{u}}
\def\pa{\partial}
\def\varep{\varepsilon}
\def\albar{\overline{\alpha}}
\def\bebar{\overline{\beta}}

\providecommand{\U}[1]{\protect\rule{.1in}{.1in}}
\def\theenumi{\arabic{enumi}}

\def\theenumii{\alph{enumii}}
\def\p@enumii{\theenumi.}

\def\theenumiii{\arabic{enumiii}}
\def\p@enumiii{(\theenumi)(\theenumii)}

\def\p@enumiv{\p@enumiii.\theenumiii}
\newtheorem{thm}{Theorem}[section]
\newtheorem{lem}[thm]{Lemma}
\newtheorem{cor}[thm]{Corollary}
\newtheorem{prop}[thm]{Proposition}
\newtheorem{remark}[thm]{Remark}

\numberwithin{equation}{section}

\begin{document}
\maketitle

\section*{Abstract}
In this work, we show the short-time existence of solutions of the evolution equations that represent the solid state dewetting of thin films through evaporation-condensation as a two dimensional sharp interface variational model. The evolution law is established as the $L^2$ gradient flow of surface energies in the presence of epitaxial strain. The main novelty is the presence of moving contact lines when the film formation is governed by the evaporation-condensation method.
\vspace{0.5 cm}\\
\textbf{Keywords:} Minimizing movements, moving contact lines, elastic thin films
\vspace{0.5 cm}\\
\textbf{Mathematics subject classification:} 35K25 $\cdot$ 35B65 $\cdot$ 35J25 $\cdot$ 74K35
\section{Introduction}
Epitaxial growth of thin crystal films plays an important role in the manufacturing of devices and circuits with applications in microelectronics, bio-engineering, optoelectronics and nanotechnology (see \cite{freund2004thin}). When growth happens in the dewetting regime, the film material dewets and destabilizes on the substrate to form an array of islands, due to the properties of the materials involved. Although this type of growth can be deleterious in some cases, it is advantageous in the manufacturing of sensors and catalysts in the formation of carbon nanotubes (see \cite{jiang2016solid}). Hence, mathematical modeling of this type of crystal growth is an important problem in materials science, as a good understanding of the growth phenomena provides us with better control over manufacturing processes, which, in turn, has significant economic and environmental impacts. 
\paragraph{} Here, we model the film formation through evaporation-condensation method, where the film profile is considered to be a grain-vapor interphase. The film growth happens through deposition of vapor on the substrate at temperatures below the melting point and the effect of mass transport that happens through surface diffusion is neglected. Hence, the dewetting happens in solid state and is different from liquid dewetting. The formation of film islands is guided by the minimization of interfacial surface energy and elastic strain energy that occurs due to crystallographic misfit between substrate and film. 

\paragraph{} The problem is modeled in two dimensions, where we restrict our attention to the cross section of a single island forming on a rigid substrate represented by $\R \times (-\infty,0]$. The island profile is represented by the height function $h: [\alpha, \beta] \rightarrow [0, \infty)$ where $h(\alpha)=h(\beta)=0$ and $h(x)>0$ for $x \in (\alpha, \beta)$, where the mass of the island remains constant in time. The region occupied by the island is given by
\begin{equation}
    \Omega_h:= \{(x,y) \in \R^2: \alpha < x< \beta,\, 0< y< h(x)\}.
\end{equation}
We assume that there is a mismatch between the lattices of the film and the substrate that gives rise to elastic displacement, which is represented by the function $u: \Omega_h \rightarrow \R^2$. This mismatch is enforced through the Dirichlet boundary condition $ u(x,0)=(e_0x, 0)$ for $x \in (\alpha, \beta)$ where $e_0>0$ is a constant. The formation of the film island is guided by the minimization of interfacial surface energy and elastic strain energy that occurs due to the crystallographic mismatch between substrate and film, under the constraint that the area of the island profile remains constant in time. Let $\gamma=\gamma_{FV}$, $\gamma_{VS}$ and $\gamma_{FS}$ represent the surface energy densities between film and vapor, vapor and substrate, and film and substrate, respectively. Let $\gamma_0= \gamma_{VS}-\gamma_{FS}$. The film dewetting happens when $\gamma>\gamma_0$. 

\paragraph{}The morphology of thin crystal films has been extensively studied in the static case in the wetting ($\gamma<\gamma_0$) \cite{bella2015study, bonacini2013epitaxially, bonacini2015stability, chambolle2002computing, chambolle2007interaction, crismale2020equilibrium, de2012regularity, fonseca2007equilibrium, fonseca2014shapes, fusco2012equilibrium, goldman2014scaling, kukta1997minimum, braides2007relaxation, capriani2013quantitative, fonseca2011material} as well as the dewetting regime ($\gamma>\gamma_0$) \cite{davoli2019analytical, piovano2022microscopical}. The time evolution of film growth in the wetting regime has been explored in \cite{piovano2014evolution, fonseca2012motion, fonseca2015motion, fusco2018surface, fusco2020surface, siegel2004evolution}. The evolution of film growth in the dewetting regime has been explored in detail by Dal Maso, Fonseca and Leoni in \cite{dal2025motion}, where the film growth is assumed to be governed by surface diffusion. To the best of our knowledge, this is the first work where the film growth due to evaporation condensation has been modeled in the dewetting regime. 

\paragraph{}In our model, the total energy involved can be represented by $\w(u, \Omega_h)+  \gamma \;\text{length}(\Gamma_h)-\gamma_0(\beta-\alpha)$, where $\w(u, \Omega_h)$ is the linearized elastic energy, $\Gamma_h$ is the graph of $h$ and $ \gamma \;\text{length}(\Gamma_h)-\gamma_0(\beta-\alpha)$ is the interfacial surface energy. As in \cite{dal2025motion, fonseca2012motion, fonseca2015motion, piovano2014evolution}, we introduce a curvature regularization term $\frac{\nu_0}{2} \int_{\Gamma_h} \kappa^2 ds$ in the surface energy, where $\nu_0>0$ is a small parameter. Hence, the underlying energy of the system is
\begin{equation}\label{eqn: surface+elastic energy}
    \w(u, \Omega_h)+  \gamma \;\text{length}(\Gamma_h)-\gamma_0(\beta-\alpha)+\frac{\nu_0}{2} \int_{\Gamma_h} \kappa^2 ds.
\end{equation}
The elastic displacement $u$ is the solution of the elastic equilibrium problem in $\Omega_h$ with the aforementioned Dirichlet boundary condition on $(\alpha, \beta) \times \{0\}$ and the natural Neumann boundary condition on $\Gamma_h$. The contact points and island profile function, $(\alpha, \beta, h)$, are obtained as the $L^2$ gradient flow of (\ref{eqn: surface+elastic energy}) using a minimizing movements argument in which approximate solutions at discrete time steps are obtained through minimization of the energy in (\ref{eqn: surface+elastic energy}) along with an incremental energy term that governs the dynamics of $h$. 

\paragraph{} The connection between the movement and curvature of idealized  grain-vapor interphases has been well established in literature. In one of the earliest works in this topic, Mullins (see \cite{mullins1956two, mullins1957theory}) derived the equations governing the phase-boundary formation in the cases of evaporation condensation as well as surface diffusion.  Let $\widetilde{V}$ be the normal velocity of the time dependent phase boundary, $\Gamma_h$. Then, in the case of evaporation condensation, $\tilde{V}$ equals the curvature $\kappa$ of $\Gamma_h$ up to a rescaling, that is, 
\begin{equation}
    \tilde{V}=\kappa.
\end{equation}
In our case, due to the presence of elastic energy and the regularization term, the evolution equation for $\tilde{V}$ is given by
\begin{gather}
    \widetilde{V}= \gamma\kappa-\nu_0\bigg(\partial_{ss}\kappa + \frac{\kappa^3}{2}\bigg)-\widetilde{W}+m, \label{eqn: normal velocity equation}
\end{gather}
 and the evolution of the contact points $\alpha, \beta$ is given by 
 \begin{gather}
     \sigma_0 \dot{\alpha} = \gamma\cos(\theta_\alpha)-\gamma_0+ \nu_0 \partial_s \kappa_\alpha \sin(\theta_\alpha) \label{eqn: contact points ODE, alpha},\\
     \sigma_0 \dot{\beta}=-\gamma\cos(\theta_\beta)+\gamma_0- \nu_0 \partial_s \kappa_\beta \sin(\theta_\beta) \label{eqn: contact points ODE, beta}
\end{gather}
where $\kappa$ is the curvature of $\Gamma_h$, $\sigma_0>0$ is a material constant, $s$ is the arc length parameter, $\widetilde{W}$ is the value of the elastic energy $u$ at the point in $\Gamma_h$ corresponding to $s$, $\theta_\alpha$ and $\theta_\beta$ are the oriented angles between $x$ axis and the tangents to $\Gamma_h$ at $(\alpha, 0)$ and $(\beta, 0)$, $\partial_s \kappa_\alpha, \partial_s \kappa_\beta$ are the derivatives of the curvature $\kappa$ with respect to $s$ at the values of $s$ corresponding to the points $(\alpha, 0)$ and $(\beta, 0)$, and $m$ is the Lagrange multiplier due to the mass constraint 
\begin{equation}\label{eqn: mass constraint}
    \int_{\alpha(t)}^{\beta(t)} h(t,x) dx = A_0,
\end{equation}
for every $t\geq 0$, where $A_0>0$ is a constant. In the absence of curvature regularization, that is, when $\nu_0=0$, (\ref{eqn: contact points ODE, alpha}) and (\ref{eqn: contact points ODE, beta}) reduce to Young's law (\cite{freund2004thin, philippis2015regularity}). We would like to note that Mullins's results has been extended to include the presence of anisotropic surface energies, where the surface energy density $\gamma$ is non-constant and dependent on the normal to the surface at every point, in seminal works like \cite{fried2004unified, gurtin1988multiphase, gurtin1995anisotropic}. We refer to \cite{fonseca2011material, piovano2014evolution} for the study of elastic thin films with anisotropic surface energies. We would also like to refer to \cite{chambolle2007interaction, crismale2020equilibrium, fusco2020surface, fonseca2015motion} for a treatment of the problem in three dimensions, and to \cite{fusco2018surface, fusco2020surface} for an approach without the curvature regularization term in the total energy.
\paragraph{}While following the minimizing movement scheme, at every discrete time step, we need to solve an elliptic equation involving the elastic displacement term in a domain with corners and obtain $W^{2,p}$ regularity estimates with exact time dependence, which is crucial in the convergence of solutions of the generalized Young's law. While we follow the approach of \cite{dal2025motion}, the main difference is that in \cite{dal2025motion} the authors study the $H^{-1}$ gradient flow associated with the energy (\ref{eqn: surface+elastic energy}), and thus the evolution equation for $h$ is a sixth order parabolic equation. In our case, we are considering the $L^2$ gradient flow of (\ref{eqn: surface+elastic energy}), and thus (\ref{eqn: normal velocity equation}) is a fourth order parabolic equation. The loss of two derivatives is a significant challenge while dealing with moving domains with corners.
We have carried out all the $W^{2,p}$ regularity estimates using the $L^1$ regularity of the third spacial derivative of $h$, unlike \cite{dal2025motion} where the estimates use up to $L^1$ regularity of the fourth spacial derivative of $h$. We have also proved the uniform convergence of the trace of the elastic displacement $u$ on the graph of the island profile $h$ away from end points, which in turn is crucially used in the proof of convergence of the elastic energy term leading up to (\ref{eqn: normal velocity equation}), unlike \cite{dal2025motion} where the convergence of trace of $u$ is weak. Further, since we have lesser regularity of $h$ than in \cite{dal2025motion}, we have a constant term in the differential equation satisfied by $h$ at every discrete time step. We have proved the convergence of these constants to a time dependent function in the limit of the time step converging to 0.

\paragraph{}Though morphology of thin films has been studied extensively, what makes this problem challenging is the presence of sharp interface angles between the film and substrate which move in time. While moving contact line problems in fluid mechanics have been studied previously (see \cite{guo2018stability}, \cite{guo2023stability}, \cite{tice2021dynamics}), there has been very few works in literature which deals with such problems in solid state film growth (see \cite{dal2025motion}). 
\paragraph{}The main result in this work is that for every initial configuration $(\alpha_0, \beta_0, h_0)$, we obtain a small time $T>0$ for which weak solutions of (\ref{eqn: normal velocity equation})-(\ref{eqn: contact points ODE, beta}) exists in $[0,T]$. 

\begin{thm}\label{thm: main thm}
    Under the assumptions (\ref{assumption: alpha_0, beta_0, h_0 are in A_s})-(\ref{assumption: Lip h_0< L_0}), there exists $T>0$ such that the following hold:
    \begin{enumerate}
        \item[(i)] There exist two functions $\alpha, \beta \in H^1((0,T))$, with $\alpha(t)<\beta(t)$ for every $t \in [0,T]$, such that $\alpha(0)= \alpha_0, \,\beta(0)=\beta_0$, and 
        \begin{equation}
            \beta(t)-\alpha(t) \geq \sqrt{\frac{2A_0}{L_0}}
        \end{equation}
        for every $t \in [0,T]$.
        \item[(ii)] There exist $1<p_1<6/5$ and a continuous function $h\geq 0$ defined for $t \in [0,T]$ and $x \in [\alpha(t), \beta(t)]$ such that 
        \begin{gather}
            h(0,x)=h_0(x) \text{ for every } x \in [\alpha_0, \beta_0],\label{eqn: h(0)=h_0 mt}\\
            h(t,x)>0  \text{ for every } t \in [0,T] \text{ and } x \in (\alpha(t), \beta(t)),\label{eqn: h>0 in (alpha, beta) mt}\\
            h(t, \alpha(t))=h(t, \beta(t))=0\text{ for every } t \in [0,T] , \label{eqn: h=0 at alpha, beta mt}\\
            h'(t, \alpha(t))>0 \text{ and } h'(t, \beta(t))<0 \text{ for every } t \in [0,T] , \label{eqn: h' sign at alpha, beta mt}
        \end{gather}
        and 
        \begin{gather}
            {\operatorname*{Lip}} \,h(t,\cdot)<L_0, \label{eqn: Lip h strict inequality mt}\\
            \int_{\alpha(t)}^{\beta(t)} h(t,x)dx =A_0, \label{eqn: area under h is constant mt}\\
            h(t,\cdot) \in H^2((\alpha(t), \beta(t)))  \label{eqn: h is in H^2 mt},\\
             h(t,\cdot) \in W^{4,p_1}((\alpha(t), \beta(t))), \label{eqn: h is in W^4, p_1 mt}
        \end{gather}
        for every $t \in [0,T]$.
        Moreover, if we denote by $h_*$ the extension of $h$ to $\R$ obtained by setting $h_*(t,x)=0$ for $t \in [0,T]$ and $x \in \R \setminus (\alpha(t), \beta(t))$, then
        \begin{equation}
            h_* \in C^{0, 1/2}([0,T]; L^2(\R)) \cap H^1((0,T); L^2(\R)). \label{eqn: h_* regularity mt}
        \end{equation}
        \item[(iii)] There exists a function $u$ defined for every $t \in [0,T]$ and $(x,y) \in \overline{\Omega}_{h(t,\cdot)}$ such that
        \begin{equation}\label{eqn: u regularity mt}
           u(t,\cdot,\cdot)\in C^{3,1-1/p_1}(\overline{\Omega}_h^{a,b}) 
        \end{equation}
        for a.e. $t \in [0, T]$ and $\alpha(t)<a<b<\beta(t)$,
        and $u(t,\cdot,\cdot)$ solves the boundary value problem
        \begin{equation}\label{eqn: u(t,.,.) bvp mt}
            \begin{cases}
                -{\operatorname*{div}} \,\mathbb{C} Eu(t,x,y)=0 \text{ in   }  \,\Omega_{h(t,\cdot)},\\
                \mathbb{C} Eu(t,x,h(t,x))\nu^h(t,x)=0 \text{ for } \, x \in (\alpha(t), \beta(t)),\\
            u(t,x,0)=(e_0x,0) \text{ for } \, x \in (\alpha(t), \beta(t)),
            \end{cases}
        \end{equation}
        for a.e. $t \in (0,T)$ where $\nu^h(t,x)= \frac{(-h'(t,x),1)}{(1+(h'(t,x))^2)^{1/2}}$ denotes the outer unit normal to $\pa \Omega_{h(t,\cdot)}$ at $(x, h(t,x))$ for $\alpha(t) <x<\beta(t)$.
        \item[(iv)]For a.e. $t \in [0,T]$ and a.e. $x \in (\alpha(t), \beta(t))$,
        \begin{gather}
             \sigma_0 \dot{\alpha}(t)= \frac{\gamma}{J(t,\alpha(t))}-\gamma_0+ \nu_0 \frac{h'(t,\alpha(t))}{(J(t,\alpha(t))^2}\bigg(\frac{h''(t,\cdot)}{(J(t,\cdot))^3}\bigg)'(\alpha(t)), \label{eqn: alpha ODE final}\\
        \sigma_0 \dot{\beta}(t)= -\frac{\gamma}{J(t,\beta(t))}+\gamma_0- \nu_0 \frac{h'(t,\beta(t))}{(J(t,\beta(t))^2}\bigg(\frac{h''(t,\cdot)}{(J(t,\cdot))^3}\bigg)'(\beta(t)), \label{eqn: beta ODE final}\\
     \dot{h} = J\bigg[  \gamma \bigg(\frac{h'}{J}\bigg)'-\nu_0 \bigg(\frac{h''}{J^5}\bigg)'' - \frac{5 \nu_0}{2}\bigg(\frac{h'(h'')^2}{J^7}\bigg)'-\overline{W}+m\bigg], \label{eqn: V de final}
        \end{gather}
     where 
        \begin{gather}
    J(t, x):= (1+ (h_*'(t,x))^2)^{1/2}, \label{eqn: def J(t,x)}
\end{gather}
for $t\geq 0$ and $x \in \R$, and $m \in L^\infty ((0,T))$ is the Lagrange multiplier due to the mass constraint (\ref{eqn: mass constraint}).
    \end{enumerate}
\end{thm}
We would like to note that an explicit formula for $m$ is given in Theorem \ref{thm: h PDE final}.
\paragraph{}The paper is organized as follows. In Section 2, we introduce the notations, prove the existence of minimizer of the total energy and state a few standard results that would be used later. Section 3 deals with the derivation of the Euler Lagrange equation associated with the energy functional and formation of some bounds on the derivatives of $h$ that are crucial in the regularity estimates. The diffeomorphism between the corners of $\Omega_h$ and a right triangular domain is established in Section 4 and the corner $W^{2,p}$ regularity estimates of $u$ are proven using the regularity estimates in triangular domain in Section 5. In Section 6, we define the discretizations of $\alpha, \beta$ and $h$ in time and prove some estimates, which in turn are used in proving the convergence of the minimization movements scheme in Section 7.

\section{Preliminaries}
 We denote the class of all admissible surface profiles by $\a_s$, which is the set of all $(\alpha, \beta, h)$ such that $ \alpha<\beta$, $h \in H^2((\alpha, \beta)) \cap H^1_0((\alpha, \beta))$, $ h\geq 0$ in $(\alpha, \beta)$, $\text{Lip}(h) \leq L_0$ and
    \begin{equation}\label{eqn: area constraint}
        \int_\alpha^\beta h(x) dx= A_0.  
   \end{equation}
 For every $(\alpha, \beta, h) \in \a_s$, we have the lower bound
 \begin{equation}\label{eqn: beta-alpha lower bound lemma 1}
        \beta- \alpha \geq \sqrt{\frac{2A_0}{L_0}},
    \end{equation}
by Lemma 2.1 in \cite{dal2025motion}. This result is crucially used in several proofs to come. 
\paragraph{}If $(\alpha, \beta, h) \in \a_s$, we define $\htil$ as 
\begin{equation}\label{eqn: def h tilde}
    \htil(x) := \begin{cases}
        h(x), \text{ if } x \in [\alpha, \beta],\\
        0, \hspace{0.54 cm}\text{ otherwise}
    \end{cases}
\end{equation}
for $x \in \R$ as the extension of $h$ by 0 outside $(\alpha, \beta)$. We represent the region between the substrate and the film surface between $(\alpha, \beta)$ by $\Omega_h$, where 
\begin{equation}\label{eqn: define Omega_h}
    \Omega_h := \{(x,y) \in \R^2: \alpha<x< \beta, 0<y< h(x)\}.
\end{equation}
When we want to restrict our attention to the film structure in a sub-interval of $(\alpha, \beta)$, say in $(a,b)$ where $\alpha < a< b <  \beta$, we use the notation $\Omega_h^{a,b}$,
\begin{equation}\label{eqn: Part of Omega_h}
    \Omega_h^{a,b}:=\{(x,y) \in \R^2: a< x< b, 0< y< h(x)\}.
\end{equation}
Further, we define the admissible class of elastic displacements in $\Omega_h$ as $\a_e$, where
 \begin{equation}
     \mathcal{A}_e(\alpha, \beta, h):= \{u \in H^1(\Omega_h;\R^2): u(x,0)=(e_0x,0) \text{ for a.e. } x \in (\alpha, \beta)\},
 \end{equation}
and the admissible class of minimizers of the total energy as $\a$, where
 \begin{equation}
     \mathcal{A}:=\{(\alpha, \beta, h,u): (\alpha, \beta, h) \in \mathcal{A}_s, u \in \mathcal{A}_e(\alpha, \beta, h)\}.
 \end{equation}

For $(\alpha, \beta, h) \in \a_s$, define the surface energy
\begin{equation}\label{eqn: define surface energy}
    \s(\alpha, \beta, h):= \gamma \int_\alpha^\beta \sqrt{1+(h'(x))^2}-\gamma_0(\beta-\alpha)+ \frac{\nu_0}{2} \int_\alpha^\beta \frac{(h''(x))^2}{(1+(h'(x))^2)^{5/2}}dx, 
\end{equation}
where $\gamma, \gamma_0>0$ and $\gamma > \gamma_0$. Note that the surface energy is non-negative since
 \begin{gather}
     \s(\alpha, \beta, h) \geq \gamma \int_\alpha^\beta \sqrt{1+(h'(x))^2}-\gamma_0(\beta-\alpha)+ \frac{\nu_0}{2} \int_\alpha^\beta \frac{(h''(x))^2}{(1+(h'(x))^2)^{5/2}}dx\nonumber \\ \geq \gamma \int_\alpha ^\beta 1 dx -\gamma_0 (\beta- \alpha) 
      \geq 0.  \label{eqn: surface energy is non-negative}
 \end{gather}
 For $(\alpha, \beta, h, u) \in \a$, define the elastic energy as 
 \begin{equation}\label{eqn: def elastic energy}
     \e(\alpha, \beta, h, u):= \int_{\Omega_h} W(Eu(x,y))dxdy,
 \end{equation}
 where $W: \R^{2 \times 2} \rightarrow [0, \infty)$ is given by
 \begin{equation}\label{eqn: def W}
     W(\xi) := \frac{1}{2}\C\xi : \xi, \text{ with } \C \xi:= \mu(\xi+\xi^T)+ \lambda (\text{tr}\xi)I,
 \end{equation}
 where $\lambda, \mu$ are the Lam\'e coefficients, $I$ is the $2 \times 2$ identity matrix and $Eu = \frac{\nabla u+ \nabla^T u}{2}$. Note that $Eu \in  \R^{2 \times 2}_{\text{sym}}$ and $C\xi= C\xi_{\text{sym}}$ where $\xi_{\text{sym}}= \frac{\xi+\xi^T}{2} \in  \R^{2 \times 2}_{\text{sym}}$. Assume $\mu>0$ and $\lambda+\mu>0$. Clearly, there exists a constant $C_W>0$ such that
 \begin{equation}\label{eqn: W upper and lower bounds by|.|^2}
     \frac{1}{C_W} |\xi|^2\leq W(\xi) \leq C_W |\xi|^2 
 \end{equation}
for all $\xi \in \R^{2 \times 2}_{\text{sym}}$.
\paragraph{} In order to apply the minimizing movements scheme, we need to set up an incremental minimization problem, which, in turn, requires a metric term that penalizes the ``distance" between minimizers of adjacent time steps. For $(\alpha, \beta,h), (\alpha_0, \beta_0, h_0)\in \a_s$, we define
 \begin{eqnarray}\label{eqn: define incremental energy functional}
    \t_{\tau} (\alpha, \beta, h; \alpha_0, \beta_0, h_0):= \frac{1}{2 \tau} \int_{\alpha_0}^{\beta_0} \frac{(\tilde{h}-h_0)^2}{\sqrt{1+(h_0'(x))^2}} dx + \frac{\sigma_0}{2 \tau}(\alpha-\alpha_0)^2 + \frac{\sigma_0}{2 \tau} (\beta-\beta_0)^2,
\end{eqnarray}
and investigate the existence and properties of minimizers of the total energy functional 
\begin{equation}\label{eqn: define total energy}
        \f^0(\alpha, \beta, h, u):= \s(\alpha, \beta, h)+ \e(\alpha, \beta, h, u)+ \t_{\tau}(\alpha, \beta, h; \alpha_0, \beta_0, h_0).
    \end{equation}
    Here, $\htil$ is as defined in (\ref{eqn: def h tilde}).
\begin{thm}\label{thm: 2_existence of minimizer}
    For every $\tau>0$ and every $(\alpha_0, \beta_0, h_0) \in \a_s$, there exists a minimizer $(\alpha, \beta, h, u) \in \a$ of the total energy functional in (\ref{eqn: define total energy}). Moreover, there exists a constant $C = C(A_0, e_0, \lambda, \mu,L_0)>0$ independent of $\alpha_0, \beta_0, h_0$ such that
    \begin{equation} \label{eqn: u H^1 upper bound}
        \|u\|_{H^1(\Omega)} \leq C
    \end{equation}
    for every minimizer $(\alpha, \beta, h, u)$ of $\f^0$ in $\a$.
\end{thm}
Before giving the proof of the theorem above, we state a Korn's inequality result as a supporting lemma, without proof. The result is proved as Lemma 2.3 in \cite{dal2025motion}.
\begin{lem} \label{lem: 3}
    Let $(\alpha, \beta, h) \in \a_s$, let $\Omega_h$ be as in (\ref{eqn: define Omega_h}) and let $1 <p<\infty$. Then, there exists a constant $C=C(L_0,p)>0$ such that 
    \begin{equation}\label{eqn: Korn's inequality}
        \int_{\Omega_h} |\nabla u|^p dx dy \leq C \int_{\Omega_h} |Eu|^p dx dy + Ce_0^p A_0
    \end{equation}
    for every $u \in W^{1,p}(\Omega_h;\R^2)$ such that $u(x,0)=(e_0x,0))$ for $x \in (\alpha, \beta)$ (in the sense of traces), where $A_0$ is the constant defined in (\ref{eqn: area constraint}).
\end{lem}
Now, we turn to the proof of the theorem.
\begin{proof}[Proof of Theorem \ref{thm: 2_existence of minimizer}]
    Let $\{(\alpha_n, \beta_n, h_n, u_n )\}_{n \in \N}$ be a minimizing sequence for (\ref{eqn: define total energy}) in $\a$. Extend $h_n$ by zero outside $(\alpha_n, \beta_n)$ to obtain $\tilde{h}_n$ for every $n$. Then using the Korn's inequality in Lemma \ref{lem: 3} and following the proof of Theorem 2 in \cite{dal2025motion}, we can prove that $\alpha_n \rightarrow \alpha$ and $\beta_n \rightarrow \beta$ with $\alpha \leq \beta$, $\tilde{h}_n \rightarrow \tilde{h}$ uniformly in $\R$ where $\tilde{h}=0$ outside $(\alpha, \beta)$ and $(\alpha, \beta, h) \in \a_s$ where $h$ is the restriction of $\tilde{h}$ to  $(\alpha, \beta)$. We also obtain that $\liminf_{n \rightarrow \infty} \s(\alpha_n, \beta_n, h_n) \geq \s(\alpha, \beta, h)$ and $\liminf_{n \rightarrow \infty} \e(\alpha_n, \beta_n, h_n, u_n) \geq \e(\alpha, \beta, h,u)$.
    \paragraph{} Now since $\tilde{h}_n \rightarrow \tilde{h}$ uniformly in $\R$, by Fatou's Lemma,
    \begin{equation*}
        \liminf_{n \rightarrow \infty} \int_{\alpha_0}^{\beta_0} \frac{(\tilde{h}_n-h_0)^2}{\sqrt{1+(h_0')^2}} dx \geq \int_{\alpha_0}^{\beta_0} \frac{(\tilde{h}-h_0)^2}{\sqrt{1+(h_0')^2}} dx.
    \end{equation*}
Further
\begin{equation}
    \lim_{n \rightarrow \infty} \bigg(\frac{\sigma_0}{2 \tau}(\alpha_n- \alpha_0)^2 + \frac{\sigma_0}{2 \tau}(\beta_n- \beta_0)^2\bigg)= \frac{\sigma_0}{2 \tau}(\alpha- \alpha_0)^2 + \frac{\sigma_0}{2 \tau}(\beta- \beta_0)^2
\end{equation}
as $\alpha_n \rightarrow \alpha $ and $\beta_n \rightarrow \beta$ as $n \rightarrow \infty$. Combining everything above, we have
\begin{equation*}
    \f^0(\alpha, \beta, h, u) \leq \liminf_{n \rightarrow \infty} \f^0(\alpha_n, \beta_n, h_n, u_n ).
\end{equation*}
Since $\{(\alpha_n, \beta_n, h_n, u_n )\}_n$ is a minimizing sequence we have that $(\alpha, \beta, h, u) \in \a$ is a minimizer of $\f^0$.
    
\end{proof}

\subsection{Some interpolation results and general Ascoli's theorem}
In this section, we state two interpolation results on Sobolev functions on intervals in $\R$, and a generalized version of Ascoli-Arzel\'a theorem, which would be used in the proofs of results in later sections.
\begin{thm}\label{Thm: 7.40, AFCSS}(Theorem 7.40, \cite{leoni2024first})
    Let $I \subseteq \R$ be an open interval, let $1\leq p,q,r \leq \infty$, $m \in \N $, $k \in \N \cup \{0\}$ with $0 \leq k \leq m$ be such that 
    \begin{equation}\label{eqn: k,m,p,q,r inequality, Thm 7.40 AFCSS}
        \bigg(1- \frac{k}{m} \bigg)\frac{1}{q} + \frac{k}{m} \frac{1}{p} \geq \frac{1}{r},
    \end{equation}
and let $u \in W^{m,1}_{\text{loc}}(I)$. Then there exists a constant $c = c(k,m,p,q,r)>0$ such that 
\begin{equation}
    \|u^{(k)}\|_{L^r(I)} \leq c l^{1/r - k -1/q}\|u\|_{L^q(I)}+ cl^{m-k-1/p+1/r}\|u^{(m)}\|_{L^p(I)}
\end{equation}
for every $0 < l< \l^1(I)$. In particular, for $p=q=r$,
\begin{equation}
    \|u^{(k)}\|_{L^p(I)} \leq c l^{- k} \|u\|_{L^p(I)}+ cl^{m-k}\|u^{(m)}\|_{L^p(I)}.
\end{equation}
\end{thm}

\begin{thm}\label{Thm: 7.41, AFCSS}(Theorem 7.41, \cite{leoni2024first})
    Let $I \subseteq \R$ be an open interval, let $1\leq p,q,r \leq \infty$, $m \in \N $, $k \in \N \cup \{0\}$ with $0 \leq k < m$, satisfy (\ref{eqn: k,m,p,q,r inequality, Thm 7.40 AFCSS}) and let $u \in W^{m,1}_{\text{loc}}(I)$ with $u \in L^q(I)$ and $u^{(m)} \in L^p(I)$. Then there exists a constant $c = c(k,m,p,q,r)>0$ such that
    \begin{equation}
        \|u^{(k)}\|_{L^r(I)} \leq c\|u\|_{L^q(I)}^{\theta}\|u^{(m)}\|_{L^p(I)}^{1-\theta},
    \end{equation}
    if $I$ has infinite length, and
    \begin{equation}
        \|u^{(k)}\|_{L^r(I)} \leq cl^{1/r - k -1/q}\|u\|_{L^q(I)}+c\|u\|_{L^q(I)}^{\theta}\|u^{(m)}\|_{L^p(I)}^{1-\theta},
    \end{equation}
for every $0 < l< \l^1(I)$, if $I$ has finite length, where
\begin{equation}
    \theta := \frac{m-k-1/p+1/r}{m-1/p+1/q}.
\end{equation}
\end{thm}
\begin{thm}\label{thm: general Ascoli Arzela}(Ascoli-Arzel\`a Theorem, Theorem 47.1, \cite{munkrestopology}) Let $(X, \tau)$ be a  topological space and $(Y,d)$ be a metric space. Give $C(X;Y)$ the topology of compact convergence. Let $\g$ be a subset of $C(X,Y)$. 
\begin{itemize}
    \item[(a)] If $\g$ is equicontinuous under d and the set 
    \begin{equation*}
        \g_a = \{f(a): f \in \g\}
    \end{equation*}
    has compact closure for each $a \in X$, then $\g$ is contained in a compact subspace of $C(X,Y)$.
    \item[(b)] The converse holds if $X$ is locally compact Hausdorff.
\end{itemize}
    
\end{thm}

\section{Euler Lagrange Equations}\label{sec: EL}
In this section, we derive the boundary value problem satisfied by the elastic displacement $u$ and the differential equation satisfied by the island profile function $h$ using standard variational techniques. In the following result, we obtain regularity results for $u$ away from the corners of $\Omega_h$ and obtain $C^3$ regularity of $h$ by bootstrapping arguments. 
\begin{thm}\label{thm: 4_EL Theorem}
Let $\tau>0$, $(\alpha_0, \beta_0, h_0) \in \a_s$ and $(\alpha, \beta, h,u) \in \a$ be a minimizer of the total energy functional given in (\ref{eqn: define total energy}). Assume that 
\begin{equation}\label{assumption: Lip h< L_0 and h>0}
    {\operatorname*{Lip}}\, h< L_0 \text{ and } h(x)>0 \text{ for all } x \in (\alpha, \beta).
\end{equation}
Then,
\begin{gather}
    h \in C^{3}((\alpha, \beta)) \cap C^4((\alpha, \beta) \setminus (\alpha_0, \beta_0)), \\
    u \in C^{2,1}(\overline{\Omega}_h^{a,b}; \R^2) \text{ for every } \alpha<a<b<\beta, \label{eqn: u C^2, 1/2 interior regularity}
\end{gather}
   and $u$ satisfies the boundary value problem 
    \begin{equation}\label{eqn: u BVP EL eqn}
         \begin{cases}
            -{\operatorname*{div}} \,\mathbb{C} Eu(x,y)=0 \text{ in   }  \,\Omega_h,\\
            \mathbb{C} Eu(x,h(x))\nu^h(x)=0 \text{ for } \, x \in (\alpha, \beta),\\
            u(x,0)=(e_0x,0) \text{ for } \, x \in (\alpha, \beta),
        \end{cases}
    \end{equation}
  where $\nu^h(x)$ denotes the outer unit normal to $\partial \Omega_h$ at $(x, h(x))$ and $\overline{\Omega}_h^{a,b}$ is defined in (\ref{eqn: Part of Omega_h}). Moreover,
  \begin{equation}\label{eqn: V Diff eq in Thm statement}
      \frac{1}{\tau}(\htil(x)-h_0(x))\chi_{[\alpha_0, \beta_0]}= J_0 \bigg[\gamma \bigg( \frac{h'}{J}\bigg)'-\nu_0\bigg( \frac{h''}{J^5}\bigg)''-\frac{5}{2}\nu_0 \bigg(\frac{h'(h'')^2 }{J^7}\bigg)'-\overline{W}+m\bigg]
  \end{equation}
 for every $x \in (\alpha, \beta)$, where $\htil$ is defined in (\ref{eqn: def h tilde}), $m$ is the Lagrange multiplier due to the mass constraint (\ref{eqn: area constraint}), and $J, J_0$ and $\overline{W}$ are defined by
   \begin{equation}\label{eqn: Def J, J_0}
       J(x):= \sqrt{1+(h'(x))^2}, \quad J_0(x):= \sqrt{1+(h_0'(x))^2}
   \end{equation}
   and 
   \begin{equation}\label{eqn: def W bar}
       \overline{W}(x):= W(Eu(x, h(x))),
   \end{equation}
   for $x \in (\alpha, \beta)$.
\end{thm}
\begin{proof}
Since $u$ minimizes $\mathcal{E}(\alpha, \beta, h,u)$, the standard variations in $\Omega_h$ produces the weak form of the elliptic boundary value problem (\ref{eqn: u BVP EL eqn}). As $h \in H^2((\alpha, \beta))$, $h' \in H^1((\alpha, \beta))$. Then, by the Fundamental theorem of calculus and H\"older's inequality, $h' \in C^{0, 1/2}([\alpha, \beta])$. That is, $h \in C^{1, 1/2}([\alpha, \beta])$. Note that we have $h(x)>0$ for every $x \in (\alpha, \beta)$. Therefore, we can use standard elliptic regularity results as in \cite{agmon1964estimates}, Theorem 9.3 to obtain that
    \begin{equation}\label{eqn: interior regularity of u EL section}
        u \in C^{1,1/2}(\overline{\Omega}_h^{a,b}; \R^2)
    \end{equation}
    for every $\alpha< a < b< \beta$.
    \paragraph{}  Fix $a_0, b_0$ such that 
    \begin{equation}\label{eqn: def a_0, b_0}
        a_0=\alpha + \frac{\beta-\alpha}{4}, \quad b_0=\beta - \frac{\beta-\alpha}{4}.
    \end{equation}
    Hence, by (\ref{eqn: beta-alpha lower bound lemma 1}), we have
    \begin{equation}\label{eqn: b_0-a_0 lower bound}
        b_0-a_0 \geq \sqrt{\frac{A_0}{2L_0}}.
    \end{equation}
    Let $a,b \in \R$ be such that $\alpha<a< a_0 < b_0 <b< \beta$. Now, extend $u$ to a function in $\Omega_h \cup ([a,b] \times \R)$, still denoted by $u$, such that $u \in C^1( [a,b] \times \R; \R^2)$. Let $\varphi \in C_c^\infty (\R)$ such that $\text{supp}(\varphi) \subseteq [a,b]$ and 
    \begin{equation}\label{eqn: EL-phi_integral 0}
     \int_\R \varphi(x) dx=   \int_\alpha^\beta
 \varphi (x) dx =0.
 \end{equation}
As $\varphi \in C_c^\infty ((\alpha, \beta))$, $\varphi \in H^2((\alpha, \beta)) \cap H^1_0((\alpha, \beta))$. Since $\text{Lip} (h) <L_0$ and $h>0$ in $(\alpha, \beta)$ , we can choose $\varepsilon_0 \in \R$ such that for all $\varepsilon \in \R$ with $|\varepsilon|\leq |\varepsilon_0|$, $\text{Lip}(h+\varepsilon \varphi ) \leq L_0$ and $h+\varepsilon \varphi \geq 0$ in $(\alpha, \beta)$. Further, note that by (\ref{eqn: EL-phi_integral 0}),
\begin{equation*}
    \int_\alpha^\beta h(x)+\varepsilon \varphi(x) dx =A_0,
\end{equation*}
and hence $(\alpha, \beta, h+\varepsilon \varphi) \in \mathcal{A}_s$. In turn, $(\alpha, \beta, h+\varepsilon \varphi, u|_{\Omega_{h+\varepsilon \varphi}}) \in \a$ by properties of the extension we chose to have. Now, by (\ref{eqn: define surface energy}),
\begin{gather*}
   \lim_{\varepsilon \rightarrow 0} \frac{1}{\varepsilon} (\s(\alpha, \beta, h+\varepsilon \varphi)-\s(\alpha, \beta, h))=\lim_{\varepsilon \rightarrow 0}\frac{1}{\varepsilon}\gamma \bigg(\int_\alpha^\beta \sqrt{1+(h'(x)+\varepsilon \varphi'(x))^2}-\sqrt{1+(h'(x))^2}dx\bigg)\\
    +\lim_{\varepsilon \rightarrow 0}\frac{1}{\varepsilon}\frac{\nu_0}{2} \bigg(\int_\alpha^\beta \bigg(\frac{(h''(x)+\varepsilon \varphi''(x))^2}{(1+(h'(x)+\varepsilon \varphi'(x))^2)^{5/2}}
       -\frac{(h''(x))^2}{(1+(h'(x))^2)^{5/2}}\bigg)dx\bigg)\\
    = \lim_{\varepsilon \rightarrow 0} \gamma \int_\alpha^\beta \frac{1}{2\sqrt{1+(h'(x)+\varepsilon \varphi'(x))^2}}2 (h'(x)+\varepsilon \varphi'(x))\varphi'(x) dx\\
    + \lim_{\varepsilon \rightarrow 0}\frac{\nu_0}{2}\bigg( \int_\alpha^\beta \frac{2(h''(x)+\varepsilon \varphi''(x))\varphi''(x)}{(1+(h'(x)+\varepsilon \varphi'(x))^2)^{5/2}}dx
     -\frac{5}{2} \int_\alpha^\beta \frac{(h''(x)+\varepsilon \varphi''(x))^2}{(1+(h'(x)+\varepsilon \varphi'(x))^2)^{7/2}} 2 (h'(x)+\varepsilon \varphi'(x))\varphi'(x)dx\bigg)\\
    = \gamma \int_\alpha^\beta \frac{h'(x)\varphi'(x)}{\sqrt{1+(h'(x))^2}} dx+\nu_0\int_\alpha^\beta \frac{h''(x)\varphi''(x)}{(1+(h'(x))^2)^{5/2}}dx
    -\frac{5\nu_0}{2} \int_\alpha^\beta \frac{(h''(x))^2h'(x)\varphi'(x)}{(1+(h'(x)^2)^{7/2}} dx,
\end{gather*}
where we have used mean value theorem and Lebesgue's dominated convergence theorem. By (\ref{eqn: def elastic energy}) and (\ref{eqn: W upper and lower bounds by|.|^2}), 
\begin{eqnarray*}
    \lim_{\varepsilon \rightarrow 0} \frac{1}{\varepsilon} (\e(\alpha, \beta, h+\varepsilon \varphi,u|_{\Omega_{h+\varepsilon \varphi}})-\e(\alpha, \beta, h,u))&=& \int_\alpha^\beta W(Eu(x, h(x))) \varphi(x) dx,
\end{eqnarray*}
and by dominated convergence theorem, 
\begin{eqnarray*}
    \lim_{\varepsilon \rightarrow 0} \frac{1}{\varepsilon} (\t (\alpha, \beta, h+ \varepsilon \varphi; \alpha_0, \beta_0, h_0)-\t (\alpha, \beta, h; \alpha_0, \beta_0, h_0))&=& \frac{1}{\tau} \int_{\alpha_0}^{\beta_0 }\frac{(\htil(x)-h_0(x))\varphi(x)}{(1+(h_0'(x))^2)^{1/2}}dx.
\end{eqnarray*}
Combining all the limits above, it follows from (\ref{eqn: define total energy}) that
\begin{gather}
    \frac{d}{d\varepsilon} (\f^0(\alpha, \beta, h+\varepsilon \varphi, u|_{\Omega_{h+\varepsilon \varphi}})\bigg|_{\varepsilon=0}=\gamma \int_\alpha^\beta \frac{h'(x)\varphi'(x)}{\sqrt{1+(h'(x))^2}} dx+\nu_0\int_\alpha^\beta \frac{h''(x)\varphi''(x)}{(1+(h'(x))^2)^{5/2}}dx\nonumber \\
    -\frac{5\nu_0}{2} \int_\alpha^\beta \frac{(h''(x))^2h'(x)\varphi'(x)}{(1+(h'(x)^2)^{7/2}} dx+\int_\alpha^\beta W(Eu(x, h(x))) \varphi(x) dx 
    +\frac{1}{\tau} \int_{\alpha_0}^{\beta_0 }\frac{(\htil(x)-h_0(x))\varphi(x)}{(1+(h_0'(x))^2)^{1/2}}dx 
    =0. \label{eqn: F derivative 0}
\end{gather}
We rewrite (\ref{eqn: F derivative 0}) as 
\begin{equation}\label{eqn: Euler lagrange eqn with A,B,f}
    \int_\R (A\htil''\varphi''+B\htil'\varphi' +f\varphi )dx=0
\end{equation}
for all $\varphi \in C_c^\infty (\R)$ with $\text{supp}(\varphi) \subseteq [a,b]$ satisfying (\ref{eqn: EL-phi_integral 0}), where
\begin{eqnarray}
    A(x)&:=& \frac{\nu_0}{(1+(\htil'(x))^2)^{5/2}}, \label{eqn: def A(X)} \\
    B(x)&:=& \gamma  \frac{1}{(1+(\htil'(x))^2)^{1/2}}-\frac{5}{2} \nu_0  \frac{(\htil''(x))^2 }{(1+(\htil'(x))^2)^{7/2}}, \label{eqn: def B(x)}\\
    f(x)&:=& W (Eu(x,\htil(x)))+\frac{1}{\tau}\frac{(\htil(x)-h_0(x))}{(1+(h_0'(x))^2)^{1/2}}\chi_{[\alpha_0, \beta_0]}.\label{eqn: def f(x)_EL}
\end{eqnarray}
Note that $\chi_{[\alpha_0, \beta_0]}$ is the characteristic function of the interval $[\alpha_0, \beta_0]$ in $\R$. Further, $A \in C^{0,1/2}([\alpha, \beta])$ and $B \in L^1((\alpha, \beta))$. Define 
\begin{eqnarray*}
   \albar:= \max\{\alpha_0, \alpha\}, \text{ and }
   \bebar:=\min\{\beta_0, \beta\}.
\end{eqnarray*}
By the regularity of $u$, $f \in L^1_{\text{loc}}((\albar, \bebar))$.\\
\paragraph{}Let $\varphi \in C_c^{\infty}([a,b])$ and let $c_\varphi := \int_a^b \varphi dx$. Let $\psi_0 \in C_c^\infty(\R)$ such that $\operatorname{supp}(\psi_0) \subset (0,1)$ and 
\begin{gather}\label{eqn: def psi_0}
    \int_0^1 \psi_0(x) dx =1, \quad \|\psi_0\|_{L^\infty(\R)} \leq 1, \quad \|\psi_0'\|_{L^\infty(\R)} \leq C_0, \quad  \|\psi_0''\|_{L^\infty(\R)} \leq C_0,
\end{gather}
for $C_0>0$. Now, define $\phi_0 \in C_c^\infty (\R)$ with $\operatorname{supp}(\phi_0) \subset (a_0, b_0)$ as 
\begin{equation}\label{eqn: def phi_0 from psi_0}
\phi_0(x):= \frac{1}{b_0-a_0} \psi_0 \bigg( \frac{x-a_0}{b_0-a_0}\bigg)
\end{equation}
for every $x \in \R$. Then, from (\ref{eqn: b_0-a_0 lower bound}) and (\ref{eqn: def psi_0}), we have 
\begin{gather}
    \int_{a_0}^{b_0} \phi_0 dx =1, \quad \|\phi_0\|_{L^\infty(\R)} \leq \sqrt{\frac{2L_0}{A_0}}, \label{eqn: phi_0 integral and L^infinity bound} \\
    \|\phi_0'\|_{L^\infty(\R)} \leq\frac{2C_0L_0}{A_0}, \quad
    \|\phi_0''\|_{L^\infty(\R)} \leq C_0 \bigg(\frac{2L_0}{A_0}\bigg)^{3/2}. \label{eqn: phi_0', phi_0'' L^infinity bound}
\end{gather}
Now, we have that $\varphi-c_\varphi \phi_0 \in C_c^\infty ([a,b]) $ with $ \int_a^b \varphi-c_\varphi \phi_0 dx =0$. Hence, from (\ref{eqn: Euler lagrange eqn with A,B,f}), we have 
\begin{eqnarray*}
    \int_a^b (Ah''\varphi''+Bh'\varphi' +f\varphi )dx= m\int_a^b \varphi dx.
\end{eqnarray*}
where 
\begin{equation}\label{eqn: def m EL section}
    m:= \int_\alpha^\beta (A h'' \phi_0''+Bh'\phi_0'+f\phi_0 )dx,
\end{equation}
since $\text{supp}(\phi_0) \subseteq [a_0, b_0] \subset [\alpha, \beta]$. Fix $x_0 \in(\alpha_0, \beta_0)$. Integrate by parts to get 
\begin{equation}\label{eqn: Integral with F_m}
    \int_a^b (Ah''\varphi''+(Bh'-F_m)\varphi' )dx= 0
\end{equation}
where 
\begin{equation}
   F_m(x):= \int_{x_0}^x (f(s)-m) ds, \label{eqn: def F_m(x)}
\end{equation}
for $x \in (\alpha, \beta)$. Let $\psi \in C_c^\infty ([a,b])$ with $\int_a^b \psi dx=0$. If $\bar{\phi}$ is defined as $\bar{\phi}(x):= \int_a^x \psi(s)ds$ for $x \in [a,b]$ we have $\bar{\phi} \in C_c^\infty ([a,b])$, $\bar{\phi}'(x) = \psi (x)$ and $\bar{\phi}''(x)=\psi'(x)$ for $x \in [a,b]$. Further, since $\bar{\phi}$ satisfies (\ref{eqn: Integral with F_m}), $\psi$ satisfies
\begin{equation}\label{eqn: ABF_m psi eqn}
    \int_a^b (Ah''\psi'+(Bh'-F_m)\psi )dx= 0.
\end{equation}
Let $\phi \in C_c^\infty ([a,b])$. Then, we have that $\phi-c_\phi \phi_0 \in C_c^\infty ([a,b]) $ with $ \int_a^b (\phi-c_\phi \phi_0) dx =0$, where $c_\phi:= \int_a^b \phi dx$. Hence, by (\ref{eqn: ABF_m psi eqn}),
\begin{equation*}
    \int_a^b A h'' (\phi'-c_\phi \phi_0' )+(Bh'-F_m)(\phi-c_\phi \phi_0) dx=0. 
\end{equation*}
This implies that
\begin{eqnarray*}
     \int_a^b A h'' \phi'+(Bh'-F_m)\phi dx &=& c_\phi \int_a^b (Ah''+(Bh'-F_m)) \phi_0' dx=m_1 \int_a^b \phi dx,
\end{eqnarray*}
where
\begin{equation}\label{eqn: def m_1}
    m_1:= \int_\alpha^\beta (Ah''+(Bh'-F_m)) \phi_0 dx.
\end{equation} 
Therefore, as we have
\begin{equation}
    \int_a^b A h'' \phi'  +(Bh'-F_m-m_1) \phi dx =0,
\end{equation}
for every $\phi \in C_c^\infty ([a,b])$, we obtain that $Bh'-F_m-m_1$ is the weak derivative of $Ah''$ in $(a,b)$. Then, by Theorem 3.20 in \cite{leoni2024first},
\begin{equation}\label{eqn: Ah''_EL equation}
    Ah'' = \int_{x_0}^x (Bh'-F_m) ds - m_1( x-x_0)+ A(x_0)h''(x_0)
\end{equation}
 for all $x \in (a,b)$. Since $m_1, m$ do not depend on $a,b$, (\ref{eqn: Ah''_EL equation}) holds in $(\alpha, \beta)$.
Note that we have $A \in C^{0,1/2}((\alpha, \beta))$, $B \in L^1((\alpha, \beta))$ and $F_m \in C^0((\alpha, \beta))$ and there exists $c_0>0$ such that $A \geq c_0$. Hence, $h'' \in C^0((\alpha, \beta))$. This implies that $A \in C^1((\alpha, \beta))$ and $B \in C^0((\alpha, \beta))$ which, in turn, implies that $h'' \in C^1((\alpha, \beta))$. Therefore, $h \in C^3((\alpha, \beta))$ and hence by standard elliptic regularity estimates\cite[Theorem 9.3]{agmon1964estimates}, we have $u \in C^{2,1}(\overline{\Omega}_h^{a,b}; \R^2)$ for every $\alpha<a<b< \beta$. Hence, (\ref{eqn: u C^2, 1/2 interior regularity}) is proved.
\paragraph{}This, in turn, implies that $f \in C^{1}((\alpha, \beta) \setminus \{\alpha_0, \beta_0\})$. Therefore, by (\ref{eqn: def F_m(x)}), $F_m \in C^{2}((\alpha, \beta) \setminus \{\alpha_0, \beta_0\})$. Therefore, by (\ref{eqn: def A(X)}) and (\ref{eqn: def B(x)}), $A \in C^2((\alpha, \beta) \setminus \{\alpha_0, \beta_0\})$ and $B \in C^1((\alpha, \beta) \setminus \{\alpha_0, \beta_0\})$. Hence, by (\ref{eqn: Ah''_EL equation}), $h \in C^{4}((\alpha, \beta) \setminus \{\alpha_0, \beta_0\})$. 
 \paragraph{} Now, differentiate (\ref{eqn: Ah''_EL equation}) twice to obtain 
\begin{equation*}
    (Ah'')''= (Bh')'-F_m'.
\end{equation*}
Hence, by (\ref{eqn: def A(X)}), (\ref{eqn: def B(x)}), (\ref{eqn: def f(x)_EL}) and (\ref{eqn: def F_m(x)}), we have 
\begin{eqnarray}\label{eqn: EL strong equation}
    \nu_0\bigg( \frac{h''}{J^5}\bigg)''=  \gamma \bigg( \frac{h'}{J}\bigg)'-\frac{5}{2}\nu_0 \bigg(\frac{h'(h'')^2 }{J^7}\bigg)'-\overline{W} -\frac{1}{\tau}\frac{(\htil(x)-h_0(x))}{J_0}\chi_{[\alpha_0, \beta_0]}+m
\end{eqnarray}
a.e. in $(\alpha, \beta)$, where $m$ is given by (\ref{eqn: def m EL section}). 

\end{proof}

\begin{remark}\label{remark: 1}
Note that the boundary value problem (\ref{eqn: u BVP EL eqn}) takes the form
    \begin{gather} \label{eqn: u BVP alter form}
    \begin{cases}
        -\mu \nabla u - (\lambda+ \mu)\nabla\,{\operatorname*{div}}\,u=0 \text{ in } \Omega_h, \\
        2\mu(Eu)\nu^h+ \lambda({\operatorname*{div}}\,u)\nu^h=0 \text{ on } \Gamma_h, \\
        u(x,0)=(e_0x,0) \text{ for } \, x \in (\alpha, \beta)
    \end{cases}
    \end{gather}
    by definition of $\C$ in (\ref{eqn: def W}), where $\Gamma_h$ is graph of $h$ in $(\alpha, \beta)$ and $\nu^h$ is the outward unit normal on $\Gamma_h$. When we analyze the behavior of $u$ in the corners of $\Omega_h$ by transforming them into triangular domains, (\ref{eqn: u BVP alter form}) would be used.
\end{remark}
Now, define the functional $B_\tau (h,h_0, \alpha, \alpha_0, \beta, \beta_0)$ by
\begin{equation} \label{eqn: def B_tau}
    B_\tau (h,h_0, \alpha, \alpha_0, \beta, \beta_0):= 1+ \frac{1}{\tau} \int_{\alpha_0}^{\beta_0}|\htil-h_0| dx + \frac{1}{\tau} |\alpha- \alpha_0|+ \frac{1}{\tau} |\beta-\beta_0|.
\end{equation}
Note that $B_\tau (h,h_0, \alpha, \alpha_0, \beta, \beta_0)\geq 1$. The following theorem provides $L^1$ and $L^\infty$ bounds on $\overline{W}, h''$ and $h'''$ in terms of $B_\tau (h,h_0, \alpha, \alpha_0, \beta, \beta_0)$. These estimates prove to be crucial in obtaining global regularity of $u$.
\paragraph{}Given $(\alpha, \beta, h) \in \a_s$, $\alpha\leq a\leq b\leq \beta$ and $\eta>0$, define
\begin{equation}\label{eqn: def Omega h,eta,a,b}
    \Omega_{h,\eta}^{a,b} := \{(x,y): a<x<b,\, h(x)- \eta< y< h(x)\}.
\end{equation}
\begin{thm}\label{thm: 6}
    Under the assumptions of Theorem \ref{thm: 4_EL Theorem}, suppose in addition that there exist $0< \eta_0<1$, $0< \eta_1<1$ and $M>1$ such that
    \begin{gather}
        2 \eta_0 \leq h'(\alpha), \quad h'(\beta) \leq -2 \eta_0, \label{eqn: h'(alpha), h'(beta) bounds}\\
        h(x) \geq 2 \eta_1\text{  for all } x \in [\alpha+\delta_0, \beta- \delta_0],\label{eqn: Lower bound on h_EL section-assumption}\\
        \int_\alpha^\beta |h''(x)|^2 dx \leq M, \label{eqn: h'' integral bound_EL section}
    \end{gather}
    where 
    \begin{equation}\label{eqn: def delta_0}
        \delta_0 := \frac{\eta_0^2}{4M}< \frac{1}{4}.
    \end{equation}
    Then, $W(Eu(\cdot, h(\cdot))) \in L^1((\alpha, \beta))$, $h \in W^{3,1}((\alpha, \beta))$, and there exists a constant $c_0= c_0(\eta_0, \eta_1, M)>1$ (independent of $\alpha, \beta, h, h_0, \tau$) such that 
    \begin{eqnarray}
        \int_\alpha^\beta W(Eu(x, h(x))) dx &\leq& c_0 B_\tau (h,h_0, \alpha, \alpha_0, \beta, \beta_0),\label{eqn: Integral of W B_tau bound}\\
        \|h''\|_{L^\infty((\alpha, \beta))}&\leq& c_0 B_\tau (h,h_0, \alpha, \alpha_0, \beta, \beta_0),\label{eqn: h'' B_tau L infinity bound}\\
        \|h'''\|_{L^1((\alpha, \beta))}&\leq& c_0 B_\tau (h,h_0, \alpha, \alpha_0, \beta, \beta_0), \label{eqn: h''' B_tau L 1 bound}
    \end{eqnarray}
where $B_\tau (h,h_0, \alpha, \alpha_0, \beta, \beta_0)$ is defined in (\ref{eqn: def B_tau}).
\end{thm}
Before we begin the proof of the theorem, we state the following lemmas (Theorem 3.3 and Lemma 3.5 in \cite{dal2025motion}) which will be used in the proof.
\begin{lem} \label{lemma: 7}
    Let $\eta_0>0$, $M>0$, $\alpha< \beta$ and $h \in H^2((\alpha, \beta))\cap H^1_0((\alpha, \beta))$. Assume that $h$ satisfies (\ref{eqn: h'(alpha), h'(beta) bounds}) and (\ref{eqn: h'' integral bound_EL section}). Then, 
    \begin{equation}\label{eqn: beta-alpha lower bound}
        \beta-\alpha \geq  \frac{16\eta_0^2}{M}.
    \end{equation}
\end{lem}
\begin{lem}\label{Thm: 5}
    Let $(\alpha, \beta, h) \in \a_s$, let $u \in \a_e(\alpha, \beta, h)$ be the minimizer of the functional $\e(\alpha, \beta, h, \cdot)$ defined in (\ref{eqn: def elastic energy}), let $\delta>0$, $\eta>0$ and $M>1$. Assume that there exist $\alpha< a< b< \beta$ with $b-a>4 \delta$ such that 
    \begin{gather}
        h(x) \geq 2 \eta \text{  for all  } x \in [a,b],\\
        \int_\alpha^\beta |h''(x)|^2 dx \leq M.
    \end{gather}
Then there exists a constant $C=C(\delta, \eta, M)>0$(independent of $\alpha, \beta, a,b, h, h_0,\tau$) such that 
\begin{equation}
    \|u\|_{C^{1, 1/2}(\overline{\Omega}_{h,\eta}^{a+\delta,b-\delta})} \leq C.
\end{equation}
\end{lem}
\begin{proof}[Proof of Theorem \ref{thm: 6}]\textbf{Step 1:} In this step, we will construct a variation of $(\alpha, \beta, h,u)$ in $\a$. Firstly, extend $h$ to $(\alpha-1, \beta)$ (without relabeling) by setting $h(x):= h'(\alpha)(x-\alpha)$ for $x \in (\alpha-1, \alpha]$. Note that $h'(x)=h'(\alpha)\geq 2\eta_0>\eta_0$ for $x \in (\alpha-1, \alpha)$ by (\ref{eqn: h'(alpha), h'(beta) bounds}). Now, for $x \in (\alpha, \beta)$ such that $(x-\alpha) \leq \frac{\eta_0^2}{M}$, using H\"older's inequality, the fundamental theorem of calculus and (\ref{eqn: h'' integral bound_EL section}),
\begin{eqnarray}
    h'(x)&=& h'(\alpha)+ \int_x^\alpha h''(s) ds \geq h'(\alpha)-(x-\alpha)^{1/2}M^{1/2} \nonumber\\
    &\geq& 2\eta_0-\eta_0=\eta_0. \label{eqn: h' geq eta_0}
\end{eqnarray}
Therefore, using (\ref{eqn: def delta_0}),
\begin{equation} \label{eqn: upper, lower bounds for h' in alpha-delta_0, alpha+ 2 delta_0}
    \eta_0 \leq h'(x) \leq L_0 \text{ for all } x \in [\alpha- \delta_0, \alpha+ 2 \delta_0].
\end{equation}
Note that $\alpha+ 4 \delta_0< \beta$ by Lemma \ref{lemma: 7}.  As $h(\alpha)=0$, using the mean value theorem, (\ref{eqn: def delta_0}), (\ref{eqn: h' geq eta_0}) and (\ref{eqn: upper, lower bounds for h' in alpha-delta_0, alpha+ 2 delta_0}),
\begin{gather}
    h(x) \geq \eta_0(x- \alpha) \text{ for every } x \in [\alpha, \alpha+ 4 \delta_0]), \label{eqn: Lower bound on h_EL section}\\
    |h(x)| \leq L_0 |x-\alpha|\text{ for every } x \in [\alpha-\delta_0, \alpha+ 2 \delta_0].\label{eqn: Upper bound on h_EL section}
\end{gather}
Choose $\varphi_0 \in C^\infty(\R)$ such that $\varphi_0(0)=1$, $\varphi_0(x) \geq 1/2$ for every $x \in [-\frac{\delta_0}{2}, \frac{\delta_0}{2} ]$, $\int_0^{\delta_0} \varphi_0 dx =0$ and supp$(\varphi_0) \subset (-\delta_0, \delta_0)$. Define 
\begin{equation}\label{eqn: define varphi_EL section}
    \varphi(x):= \varphi_0(x- \alpha), \text{ for } x \in \R
\end{equation}
and $\varep_0 :=  \min \bigg\{1, \frac{1}{2} \frac{\delta_0 \eta_0}{\|\varphi_0\|_{C^1(\R)}}\bigg\}$. Then, for every $ \varep \in \R$ such that $|\varep| \leq \varep_0$, $-\frac{\eta_0}{2} \leq \varep \varphi'(x) \leq \frac{\eta_0}{2}$ for every $x \in \R$. Define $f(\varep, x):= h(x)+ \varep \varphi(x)$ for $|\varep| \leq \varep_0$ and $x \in \R$. Therefore, as $\varphi(\delta_0)= \varphi(-\delta_0)=0$, by (\ref{eqn: upper, lower bounds for h' in alpha-delta_0, alpha+ 2 delta_0}), 
\begin{gather}
    h'(x)+ \varep \varphi'(x) \geq  \frac{\eta_0}{2} \text{ for every } x \in [\alpha- \delta_0, \alpha+ \delta_0], \label{eqn: h' + varep varphi' geq eta_0/2}\\
    h(\alpha- \delta_0) + \varep \varphi(\alpha- \delta_0)  <0< h(\alpha+ \delta_0) + \varep \varphi(\alpha+\delta_0).
\end{gather}
By the intermediate value theorem, there exists $\alpha_\varep \in (\alpha- \delta_0, \alpha+ \delta_0)$ such that 
\begin{equation}\label{eqn: alpha_varepsilon condition}
    h(\alpha_\varep)+ \varep \varphi(\alpha_\varep)=0.
\end{equation}
Note that when $\varep=0$, $\alpha_0=\alpha$ satisfies the condition above. Now, using the implicit function theorem, the function $\varep \mapsto \alpha_\varep$ is of class $C^1$ and 
\begin{equation*}
    h'(\alpha_\varep) \frac{d \alpha_\varep}{d\varep}+\varep \varphi'(\alpha_\varep) \frac{d \alpha_\varep}{d\varep}+ \varphi(\alpha_\varep)=0,
\end{equation*}
that is, 
\begin{equation}\label{eqn: derivative of alpha_epsilon wrt epsilon}
    \frac{d \alpha_\varep}{d\varep} = \frac{-\varphi(\alpha_\varep)}{ h'(\alpha_\varep)+\varep \varphi'(\alpha_\varep)}.
\end{equation}
By (\ref{eqn: h' + varep varphi' geq eta_0/2}), we get
\begin{equation}
    \bigg|\frac{d \alpha_\varep}{d\varep} \bigg| \leq \frac{2 \|\varphi_0\|_{L^\infty(\R)}}{\eta_0} \text{ for all } \varep \in [-\varep_0, \varep_0].
\end{equation}
Therefore, by the mean value theorem applied to the function $\varep \mapsto \alpha_\varep$, we have 
\begin{equation}\label{eqn: alpha_epsilon-alpha bound}
    |\alpha_\varep-\alpha| \leq \frac{2 \|\varphi_0\|_{L^\infty(\R)}|\varep|}{\eta_0}.
 \end{equation}
Define $\varep_1 := \min \bigg\{\frac{ \delta_0 \eta_0}{4\|\varphi_0\|_{L^\infty(\R)}}, \varep_0\bigg\}>0$. Then, for $|\varep| \leq \varep_1$, by
(\ref{eqn: alpha_epsilon-alpha bound}), $|\alpha_\varep-\alpha| \leq \frac{\delta_0}{2}$ and hence by construction, $\varphi(\alpha_\varep) \geq \frac{1}{2}$. Therefore, by (\ref{eqn: h' + varep varphi' geq eta_0/2}) and (\ref{eqn: derivative of alpha_epsilon wrt epsilon}), 
\begin{equation}\label{eqn: derivative of alpha_epsilon <0}
    \frac{d \alpha_\varep}{d\varep} < 0 \text{ for } |\varep| \leq \varep_1.
\end{equation}
Now, by (\ref{eqn: h' + varep varphi' geq eta_0/2}) and (\ref{eqn: alpha_varepsilon condition}), we have $h+ \varep \varphi >0$ in $(\alpha_\varep, \alpha+ \delta_0)$.
\paragraph{} We need a variation of $h$. Fix a function $\psi_0 \in C^\infty(\R)$ such that supp$(\psi_0) \subset (\delta_0, 2\delta_0)$ and $\int_{\delta_0}^{2 \delta_0} \psi_0 dx =1$. Define $\psi(x):= \psi_0(x-\alpha)$ for $x \in \R
$. Now, define $h_\varep$ as 
\begin{equation}\label{eqn: define h_epsilon}
    h_\varep(x):= h(x)+ \varep \varphi(x)+ \omega_\varep \psi(x),
\end{equation}
where $\omega_\varep$ is chosen such that 
\begin{equation}\label{eqn: h_epsilon area constraint}
    \int_{\alpha_\varep}^\beta h_\varep dx =A_0.
\end{equation}
Note that $h_\varep \in H^1_0((\alpha_\varep, \beta))\cap H^2((\alpha_\varep, \beta))$ by definition. 
\paragraph{} We need to understand the growth of $\omega_\varep$ with respect to $\varep$ in order to prove that $h_\varep >0$ in $(\alpha_\varep, \beta)$ for $\varep $ sufficiently small. Since supp$(\psi) \subset (\alpha+ \delta_0, \alpha+ 2 \delta_0)$, $\alpha_\varep < \alpha+ \delta_0$, and $\int_{ \delta_0}^{2 \delta_0} \psi_0 dx =1$,
\begin{eqnarray*}
    \int_{\alpha_\varep}^\beta \psi dx = \int_{\alpha_\varep}^{\alpha+ \delta_0} \psi dx+\int_{\alpha+\delta_0}^{\alpha+2\delta_0} \psi dx+\int_{\alpha+ 2 \delta_0}^\beta \psi dx=1.
\end{eqnarray*}
Therefore, using $\int_\alpha^\beta \varphi dx =0$, by (\ref{eqn: area constraint}) and (\ref{eqn: h_epsilon area constraint}), we have
\begin{eqnarray*}
   A_0 &=& \int_\alpha^ \beta(h + \varep \varphi )dx = \int_{\alpha_\varep}^\beta( h+ \varep \varphi+ \omega_\varep \psi) dx\\
    &=& \int_{\alpha_\varep}^\alpha (h + \varep \varphi )dx +  \int_\alpha^ \beta(h + \varep \varphi )dx+ \omega_\varep
\end{eqnarray*}
and hence 
\begin{equation*}
    \frac{\omega_\varep}{\varep} = -\frac{1}{\varep} \int_{\alpha_\varep}^\alpha (h + \varep \varphi )dx.
\end{equation*}
By (\ref{eqn: Upper bound on h_EL section}) and (\ref{eqn: alpha_epsilon-alpha bound}), we have 
\begin{gather*}
   \bigg| \int_{\alpha_\varep}^\alpha h dx\bigg| \leq \frac{4 L_0\|\varphi_0\|^2_{L^\infty(\R)}\varep^2}{\eta_0^2}, \text{ and }
    \bigg|\int_{\alpha_\varep}^\alpha \varphi dx \bigg|\leq \frac{2 \|\varphi_0\|^2_{L^\infty(\R)}|\varep|}{\eta_0}.
\end{gather*}
Using these estimates in the equation above, we obtain
\begin{equation}\label{eqn: omega_epsilon estimate}
    \bigg|\frac{\omega_\varep}{\varep}\bigg|\leq C|\varep|.
\end{equation}
Then, by (\ref{eqn: upper, lower bounds for h' in alpha-delta_0, alpha+ 2 delta_0}),
\begin{eqnarray*}
    h_\varep' &=& h'+ \varep \varphi' + \omega_\varep \psi'\\
    &\geq& \eta_0- |\varep|\|\varphi_0'\|_{L^\infty(\R)}- C\varep^2 \|\psi_0'\|_{L^\infty(\R)},
\end{eqnarray*}
in $(\alpha_\varep, \alpha+ 2 \delta_0)$. Choose $\varep_2 := \min\bigg\{\frac{\eta_0}{2 \|\varphi_0'\|_{L^\infty(\R)}+2C\|\psi_0'\|_{L^\infty(\R)} }, \varep_1\bigg\}$. Then, if $|\varep| \leq \varep_2$, 
\begin{equation}\label{eqn: Lower bound on h_epsilon'}
    h_\varep' \geq \frac{\eta_0}{2}.
\end{equation}
Now, as $\alpha_\varep< \alpha+ \delta_0$ and supp$(\psi) \in (\alpha+ \delta_0, \alpha+ 2\delta_0)$, $\psi(\alpha_\varep)=0$. Hence, by (\ref{eqn: alpha_varepsilon condition}), $h_\varep(\alpha_\varep)=0$. Therefore, by (\ref{eqn: Lower bound on h_epsilon'}), $h_\varep>0$ in $(\alpha_\varep, \alpha+ 2 \delta_0)$. Now, since supp$(\varphi) \in (\alpha-\delta_0, \alpha+ \delta_0)$ and supp$(\psi) \in (\alpha+ \delta_0, \alpha+ 2\delta_0)$, \begin{equation*}
    h_\varep=h \text{ in } [\alpha+ 2 \delta_0, \beta],
\end{equation*}
that is, $h_\varep(\alpha_\varep)= h_\varep(\beta)=0$ and $h_\varep>0$ in $(\alpha_\varep, \beta)$. Combined with (\ref{eqn: h_epsilon area constraint}), this implies that $(\alpha_\varep, \beta, h_\varep) \in \a_s$.
\paragraph{} Now, we construct a variation in the set $\a$. Let $U$ be the interior of the set $\overline{\Omega}_h \cup ([\alpha- \delta_0, \beta] \times [-1, 0])$ and let $\hat{u} : \overline{U} \rightarrow \R^2$ be the function defined by 
\begin{equation}
    \hat{u} (x,y):= \begin{cases}
        u(x,y) \text{ if } (x,y) \in \overline{\Omega}_h , \\
        (e_0x, 0) \text{ if } (x,y) \in [\alpha-\delta_0, \beta] \times [-1,0].
    \end{cases}
\end{equation}
Since both definitions match on $[\alpha, \beta] \times \{0\}$, $\hat{u} \in H^1(U; \R^2)$. Since $\partial U$ has Lipschitz boundary, we can extend $\hat{u}$ to $\R^2$, without relabeling, to obtain $\hat{u} \in H^1(\R^2; \R^2)$. Now, define $u_\varep$ as the restriction of $\hat{u}$ to $\Omega_{h_\varep}:= \{(x,y) \in \R^2: \alpha_\varep < x < \beta, 0< y <h_\varep\}$. Note that $u_\varep \in H^1(\Omega_{h_\varep};\R^2)$ and $u_\varep(x,0)=(e_0x,0)$ for a.e. $x \in (\alpha_\varep, \beta)$. Therefore, $(\alpha_\varep, \beta, h_\varep, u_\varep) \in \a$.\\\\
\textbf{Step 2:} We claim that 
 \begin{gather}
     \liminf_{ \varep \rightarrow 0}\frac{S(\alpha_\varep, \beta, h_\varep)-\s(\alpha, \beta,h)}{\varep}\geq  \gamma \int_\alpha^\beta \frac{h'\varphi'}{\sqrt{1+(h')^2}}dx + \gamma \sqrt{1+h'(\alpha))^2}\frac{1}{h'(\alpha)}-\gamma_0 \frac{1}{h'(\alpha)}\nonumber\\
    +\nu_0 \int_\alpha^\beta \frac{h''\varphi''}{(1+(h')^2)^{5/2}}dx-\frac{5}{2} \nu_0 \int_\alpha^\beta \frac{h'(h'')^2\varphi'}{(1+(h')^2)^{7/2}}dx. \label{eqn: surface energy derivative claim}
 \end{gather}
 By (\ref{eqn: define surface energy}), 
 \begin{eqnarray}
     \frac{S(\alpha_\varep, \beta, h_\varep)-\s(\alpha, \beta,h)}{\varep}&=& \frac{1}{\varep}\bigg(\gamma\int_{\alpha_\varep}^\beta \sqrt{1+(h_\varep')^2} dx- \gamma \int_{\alpha}^\beta \sqrt{1+(h')^2} dx\bigg)\\& &-\frac{1}{\varep}\gamma_0((\beta-\alpha_\varep)-(\beta-\alpha))\nonumber\\
     & & +\frac{1}{\varep}\frac{\nu_0}{2} \bigg( \int_{\alpha_\varep}^\beta \frac{(h_\varep'')^2}{(1+(h_\varep')^2)^{5/2}}dx-\int_\alpha^\beta \frac{(h'')^2}{(1+(h')^2)^{5/2}}dx \bigg)\nonumber \\
     &=:& I_{\varep,1}+I_{\varep,2}+I_{\varep,3}. \label{eqn: surface energy difference quotient break down}
 \end{eqnarray}
 By (\ref{eqn: derivative of alpha_epsilon wrt epsilon}), we have 
 \begin{eqnarray}
   \lim_{\varep \rightarrow 0} I_{\varep,2}=   \lim_{\varep \rightarrow 0}\gamma_0 \frac{(\alpha_\varep-\alpha)}{\varep} = \gamma_0 \frac{d \alpha_\varep}{d \varep}\Bigg\vert_{\varep=0} = \gamma_0 \frac{-\varphi(\alpha)}{h'(\alpha)}=-\gamma_0 \frac{1}{h'(\alpha)}, \label{eqn: I_epsilon, 2 limit_s energy}
 \end{eqnarray}
 as $\varphi(\alpha)=1$. Further,
 \begin{eqnarray}
      \lim_{\varep \rightarrow 0} I_{\varep,1} &=& \bigg( \gamma \int_{\alpha_\varep}^\beta \frac{1}{2\sqrt{1+(h_\varep')^2}} 2 h_{\varep}' (\varphi'+ o(\varep)\psi')dx \nonumber\\
     & & -\gamma \sqrt{1+(h_\varep'(\alpha_\varep))^2} \bigg(-\frac{\varphi(\alpha_\varep)}{h'(\alpha_\varep)+\varep \varphi'(\alpha_\varep)}\bigg)\bigg)\Bigg \vert_{\varep=0}\nonumber\\
     &=& \gamma \int_\alpha^\beta \frac{h'\varphi'}{\sqrt{1+(h')^2}}dx + \gamma \sqrt{1+h'(\alpha))^2}\frac{1}{h'(\alpha)}, \label{eqn: I_epsilon, 1 limit_s energy}
 \end{eqnarray}
 by (\ref{eqn: derivative of alpha_epsilon wrt epsilon}), (\ref{eqn: define h_epsilon}) and (\ref{eqn: omega_epsilon estimate}). Write 
 \begin{eqnarray}
     I_{\varep,3} &=& \frac{1}{\varep}\frac{\nu_0}{2} \bigg( \int_{\alpha_\varep}^\beta \frac{(h_\varep'')^2}{(1+(h_\varep')^2)^{5/2}}dx-\int_{\alpha_\varep}^\beta \frac{(h'')^2}{(1+(h')^2)^{5/2}}dx \bigg)\nonumber\\
     & & + \frac{1}{\varep}\frac{\nu_0}{2} \int_{\alpha_\varep}^\alpha \frac{(h'')^2}{(1+(h')^2)^{5/2}}dx =: I_{\varep,3,1}+I_{\varep,3,2}. \label{eqn: splitting of I_(epsilon,3)}
 \end{eqnarray}
 By (\ref{eqn: derivative of alpha_epsilon wrt epsilon}), (\ref{eqn: define h_epsilon}) and (\ref{eqn: omega_epsilon estimate}) again,
\begin{eqnarray}
     \lim_{\varep \rightarrow 0} I_{\varep,3,1} = \nu_0 \int_\alpha^\beta \frac{h''\varphi''}{(1+(h')^2)^{5/2}}dx-\frac{5}{2} \nu_0 \int_\alpha^\beta \frac{h'(h'')^2\varphi'}{(1+(h')^2)^{7/2}}dx. \label{eqn: I_epsilon,3,1 limit_s energy}
\end{eqnarray}
We infer from (\ref{eqn: derivative of alpha_epsilon <0}) that $\alpha_\varep< \alpha$ when $\varep>0$ and $\alpha_\varep>\alpha$ when $\varep<0$ for $|\varep|\leq \varep_2$. Hence, as $\frac{(h'')^2}{(1+(h')^2)^{5/2}} \geq 0$ everywhere, we have
\begin{equation*}
    \liminf_{\varep \rightarrow 0} I_{\varep,3,2} \geq 0.
\end{equation*}
This proves the claim (\ref{eqn: surface energy derivative claim}). 
\paragraph{}Since $h\geq 0$ in $[\alpha, \beta]$, by (\ref{eqn: area constraint}) and (\ref{eqn: Lower bound on h_EL section-assumption}), 
\begin{eqnarray*}
    A_0 \geq \int_{\alpha+ \delta_0}^{\beta- \delta_0} h dx \geq 2 \eta_1 ( \beta- \alpha - 2 \delta_0).
\end{eqnarray*}
Therefore,
\begin{equation}\label{eqn: Upper bound on beta-alpha}
    \beta- \alpha \leq \frac{A_0}{2 \eta_1}+ 2 \delta_0=:C_0.
\end{equation}
Hence, by (\ref{eqn: h'(alpha), h'(beta) bounds}), (\ref{eqn: h'' integral bound_EL section}), (\ref{eqn: def delta_0}), (\ref{eqn: beta-alpha lower bound}), (\ref{eqn: surface energy derivative claim}), (\ref{eqn: Upper bound on beta-alpha}) and H\"older's inequality,
\begin{gather}
    \liminf_{ \varep \rightarrow 0}\frac{\s(\alpha_\varep, \beta, h_\varep)-\s(\alpha, \beta,h)}{\varep}  \geq-\bigg(\gamma L_0 C_0+ \nu_0 M^{1/2} C_0^{1/2}+ \frac{5}{2} \nu_0 L_0 M\bigg) \|\varphi_0\|_{C^2_b(\R)}-\frac{\gamma_0}{2 \eta_0}\nonumber \\
   =: -C_1 \geq -C_1 B_\tau, \label{eqn: Lower limit on liminf of surface energy}
\end{gather}
where $C_1>0$, as $B_\tau=B_\tau(h,h_0, \alpha, \alpha_0, \beta, \beta_0)\geq 1$ .
\paragraph{}Define $\t^0(\cdot,\cdot,\cdot):= \t_\tau (\cdot,\cdot,\cdot,\alpha_0, \beta_0, h_0)$. 
 \begin{equation*}
     \t^0 (\alpha_\varep, \beta, h_\varep)= \frac{1}{2 \tau} \int_{\alpha_0}^{\beta_0} \frac{(\htil_\varep -h_0)^2}{\sqrt{1+(h_0')^2}} dx + \frac{\sigma_0}{2 \tau } (\alpha_\varep-\alpha_0)^2+\frac{\sigma_0}{2 \tau } (\beta-\beta_0)^2,
 \end{equation*}
 where $\htil_\varep$ is the extension of $h_\varep$ by 0 outside $(\alpha_\varep, \beta)$. We need to evaluate $\frac{\partial}{\pa \varep } \t^0 (\alpha_\varep, \beta, h_\varep)\Bigg\vert_{\varep=0}$. Clearly, $\frac{\partial}{\pa \varep }\bigg(\frac{\sigma_0}{2 \tau } (\beta-\beta_0)^2\bigg)=0$ for every $\varep$. By (\ref{eqn: derivative of alpha_epsilon wrt epsilon}),
 \begin{gather*}
      \frac{\partial}{\pa \varep }\bigg( \frac{\sigma_0}{2 \tau } (\alpha_\varep-\alpha_0)^2\bigg)\Bigg\vert_{\varep=0}= \frac{\sigma_0}{ \tau }(\alpha_\varep-\alpha_0)\frac{\partial \alpha_\varep}{\pa \varep }\Bigg\vert_{\varep=0}\\
      = -\frac{\sigma_0}{ \tau }(\alpha_\varep-\alpha_0)\frac{\varphi(\alpha_\varep)}{h'(\alpha_\varep)+ \varep \varphi'(\alpha_\varep)}\Bigg\vert_{\varep=0} = -\frac{\sigma_0}{ \tau }(\alpha-\alpha_0)\frac{1}{h'(\alpha)}.
 \end{gather*}
 Further, by Lebesgue's dominated convergence theorem,
    \begin{eqnarray*}
       \frac{\partial}{\pa \varep }\frac{1}{2 \tau}\bigg(\int_{\alpha_0}^{\beta_0} \frac{(\htil_\varep -h_0)^2}{\sqrt{1+(h_0')^2}} dx \bigg) \Bigg\vert_{\varep=0}&=& \frac{1}{2 \tau} \bigg(\int_{\alpha_0}^{\beta_0} \frac{2(\htil_\varep -h_0)}{\sqrt{1+(h_0')^2}} (\varphi+o(\varep) \psi)dx \bigg) \Bigg\vert_{\varep=0}\\
       &=& \frac{1}{ \tau} \int_{\alpha_0}^{\beta_0} \frac{(\htil-h_0)\varphi}{\sqrt{1+(h_0')^2}} dx .
    \end{eqnarray*}
Therefore, using (\ref{eqn: h'(alpha), h'(beta) bounds}),
\begin{eqnarray}
    \frac{\partial}{\pa \varep } \t^0 (\alpha_\varep, \beta, h_\varep)\Bigg\vert_{\varep=0} &=& \frac{1}{ \tau} \int_{\alpha_0}^{\beta_0} \frac{(\htil-h_0)\varphi}{\sqrt{1+(h_0')^2}} dx- \frac{\sigma_0}{ \tau }\frac{(\alpha-\alpha_0)}{h'(\alpha)} \nonumber \\
    &\geq& -\frac{\|\varphi\|_{L^\infty(\R)}}{ \tau} \int_{\alpha_0}^{\beta_0}|\htil-h_0| dx  -\frac{\sigma_0}{ 2\tau }\frac{|\alpha-\alpha_0|}{\eta_0}. \label{eqn: derivative of incremental energy term}
\end{eqnarray}
Therefore, we have $C_2>0$ such that 
\begin{eqnarray}\label{eqn: Lower limit on liminf of incremental energy}
    \liminf_{\varep \rightarrow 0}\frac{\t^0 (\alpha_\varep, \beta, h_\varep)- \t^0 (\alpha, \beta, h)}{\varep} \geq -C_2 \geq -C_2 B_\tau.
\end{eqnarray}\\\\
\textbf{Step 3:} In this step, we prove (\ref{eqn: Integral of W B_tau bound}).
\paragraph{} Since $(\alpha_\varep, \beta, h_\varep, u_\varep) \in \a$, and $(\alpha, \beta, h,u) \in \a$ minimizes the total energy $\f^0$, we have
\begin{equation*}
    \f^0(\alpha_\varep, \beta, h_\varep, u_\varep)-\f^0 (\alpha, \beta, h,u) \geq 0
\end{equation*}
and hence
\begin{equation*}
    \limsup_{\varep \rightarrow 0^-}\frac{ \f^0(\alpha_\varep, \beta, h_\varep, u_\varep)-\f^0 (\alpha, \beta, h,u)}{\varep} \leq 0. 
\end{equation*}
Therefore, using (\ref{eqn: define total energy}), (\ref{eqn: Lower limit on liminf of surface energy}) and (\ref{eqn: Lower limit on liminf of incremental energy}), we can find $C>0$ such that 
\begin{equation} \label{eqn: Elastic energy liminf bound_EL section}
     \limsup_{ \varep \rightarrow 0^-}\frac{\e(\alpha_\varep, \beta, h_\varep, u_\varep)-\e(\alpha, \beta,h,u)}{\varep} \leq CB_\tau(h,h_0, \alpha, \alpha_0, \beta, \beta_0). 
\end{equation}
By choice of $\varep_1$, (\ref{eqn: alpha_epsilon-alpha bound}) and (\ref{eqn: derivative of alpha_epsilon <0}), $\alpha< \alpha_\varep < \alpha+ \frac{\delta_0}{2}$ when $-\varep_1< \varep<0$. Since supp$(\varphi) \in (\alpha-\delta_0, \alpha+\delta_0)$ and supp$(\psi) \in (\alpha+\delta_0, \alpha+2\delta_0)$, $h_\varep=h$ in $(\alpha+2 \delta_0, \beta)$ and hence
\begin{gather}
    \frac{\e(\alpha_\varep, \beta, h_\varep, u_\varep)-\e(\alpha, \beta,h,u)}{\varep}= \frac{1}{\varep}\int_{\alpha_\varep}^\beta \bigg(\int_0^{h_\varep(x)} W(E\hat{u}(x,y))dy\bigg) dx \nonumber\\
    -\frac{1}{\varep}\int_\alpha^\beta \bigg(\int_0^{h(x)} W(Eu(x,y))dy\bigg) dx \nonumber \\
    = -\frac{1}{\varep} \int_{\alpha_\varep}^{\alpha+2 \delta_0}\bigg(\int_{h_\varep(x)}^{h(x)} W(E\hat{u}(x,y))dy \bigg)dx  -\frac{1}{\varep} \int_{\alpha}^{\alpha_\varep} \bigg(\int_0^{h(x)} W(Eu(x,y))dy \bigg)dx \nonumber\\
    \geq -\frac{1}{\varep}\int_{\alpha_\varep}^{\alpha+2 \delta_0}\bigg(\int_{h_\varep(x)}^{h(x)} W(E\hat{u}(x,y))dy \bigg)dx \label{eqn: Elastic energy lim inf inequality_EL section}
\end{gather}
where we used the facts that $W(Eu)\geq 0$ everywhere and $\varep<0$. 
\paragraph{} Note that $h(x) \geq \frac{\eta_0 \delta_0}{4}$ for $x \in [\alpha+ \frac{\delta_0}{4}, \alpha+ 4 \delta_0]$ by (\ref{eqn: Lower bound on h_EL section}). In the statement of Lemma \ref{Thm: 5}, take $\delta= \frac{\delta_0}{4}, \eta=\frac{\eta_0 \delta_0}{8},  a= \alpha+ \frac{\delta_0}{4}$ and $b=\alpha+ \frac{9\delta_0}{4}$. Therefore, $\alpha< a< b< \beta$ and $b-a > 4 \delta$ by (\ref{eqn: def delta_0}) and (\ref{eqn: beta-alpha lower bound}). Using (\ref{eqn: h'' integral bound_EL section}), by Lemma \ref{Thm: 5}, we have $C>0$ such that
\begin{gather}
    \|W(Eu(\cdot, h(\cdot)))\|_{L^\infty([\alpha+ \frac{\delta_0}{4}, \alpha+ 2 \delta_0])} \leq \|u\|_{C^{1,1/2}(\overline{\Omega}^{a+ \delta, b-\delta}_{h,\eta})}^2\nonumber\\ \leq C \leq CB_\tau (h,h_0, \alpha, \alpha_0, \beta, \beta_0).
\end{gather}
Therefore, by Lebesgue's dominated convergence theorem, definition of $h_\varep$ and (\ref{eqn: W upper and lower bounds by|.|^2}), we have 
\begin{equation}\label{eqn: W integral limit_EL section}
    \lim_{\varep \rightarrow 0^-} \frac{1}{\varep} \int_{\alpha+ \frac{\delta_0}{2}}^{\alpha+ 2 \delta_0} \bigg(\int_{h_\varep(x)}^{h(x)} W (E \hat{u}(x,y)) dy\bigg) dx = -\int_{\alpha+ \frac{\delta_0}{2}}^{\alpha+ 2 \delta_0}W(Eu(x, h(x))) dx \leq 0.
\end{equation}
Hence, by (\ref{eqn: Elastic energy liminf bound_EL section}), (\ref{eqn: Elastic energy lim inf inequality_EL section}) and (\ref{eqn: W integral limit_EL section}), we have 
\begin{gather}
    \limsup_{\varep \rightarrow 0^-}\bigg( -\frac{1}{\varep} \int_{\alpha_\varep}^{\alpha + \frac{\delta_0}{2}} \bigg(\int_{h_\varep(x)}^{h(x)} W (E \hat{u}(x,y)) dy\bigg) dx\bigg) \nonumber\\\leq \limsup_{\varep \rightarrow 0^-}\bigg( -\frac{1}{\varep} \int_{\alpha_\varep}^{\alpha+  2 \delta_0} \bigg(\int_{h_\varep(x)}^{h(x)} W (E \hat{u}(x,y)) dy\bigg) dx\bigg) \nonumber \\
    + \limsup_{\varep \rightarrow 0^-}\bigg( \frac{1}{\varep} \int_{\alpha + \frac{\delta_0}{2}}^{\alpha+ 2 \delta_0} \bigg(\int_{h_\varep(x)}^{h(x)} W (E \hat{u}(x,y)) dy\bigg) dx\bigg) \nonumber \\
   \leq  \limsup_{ \varep \rightarrow 0}\frac{\e(\alpha_\varep, \beta, h_\varep, u_\varep)-\e(\alpha, \beta,h,u)}{\varep} 
\leq  CB_\tau(h,h_0, \alpha, \alpha_0, \beta, \beta_0). \label{eqn: Elastic energy limsup bound_ EL section}
\end{gather}
Let $x \in (\alpha, \alpha+ \frac{\delta_0}{2})$. Since $\varep\mapsto \alpha_\varep$ is a $C^1$ function and $\alpha_\varep=0$ at $\varep=0$, 
\begin{equation*}
    \lim_{\varep \rightarrow 0^-} \chi_{(\alpha_\varep, \alpha+ \frac{\delta_0}{2})}(x)=1.
\end{equation*}
Further, since $\psi(x)=0$ in $(\alpha, \alpha+ \frac{\delta_0}{2})$, by (\ref{eqn: define h_epsilon}),
\begin{eqnarray*}
    \lim_{\varep \rightarrow 0^-} -\frac{1}{\varep} \int_{h_\varep(x)}^{h(x)} W (E \hat{u}(x,y)) dy = W (Eu(x, h(x)))\varphi(x).
\end{eqnarray*}
Also note that when $-\varep<\varep<0$, by (\ref{eqn: derivative of alpha_epsilon <0}), (\ref{eqn: define h_epsilon}) and definitions of $\varphi, \psi$, 
\begin{equation*}
    -\frac{1}{\varep} \int_{h_\varep(x)}^{h(x)} W (E \hat{u}(x,y)) dy \geq 0.
\end{equation*}
Since $\varphi \geq 1/2$ in $[\alpha, \alpha+ \frac{\delta_0}{2}]$, using Fatou's lemma and (\ref{eqn: Elastic energy limsup bound_ EL section}), we have
\begin{eqnarray*}
    \frac{1}{2} \int_{\alpha}^{\alpha+ \frac{\delta_0}{2}} W(Eu(x, h(x)))dx &\leq& \limsup_{\varep \rightarrow 0^-}\bigg( -\frac{1}{\varep} \int_{\alpha_\varep}^{\alpha+ \frac{\delta_0}{2}} \bigg(\int_{h_\varep(x)}^{h(x)} W (E \hat{u}(x,y)) dy\bigg) dx\bigg)\\
    &\leq&  CB_\tau(h,h_0, \alpha, \alpha_0, \beta, \beta_0). 
\end{eqnarray*}
The corresponding inequality in the interval $[\beta-\frac{\delta_0}{2}, \beta]$ can be obtained similarly. In the interval $[\alpha+\frac{\delta_0}{2}, \beta-\frac{\delta_0}{2}]$, the inequality can be obtained using (\ref{eqn: Upper bound on beta-alpha}) and the result in Lemma \ref{Thm: 5}. Note that the inequality $1\leq B_\tau(h,h_0, \alpha, \alpha_0, \beta, \beta_0)$ is used here. Combining all the inequalities above, we obtain (\ref{eqn: Integral of W B_tau bound}).\\\\
\textbf{Step 4:} We prove the rest of the claims in the theorem. Fix $x_0 \in (\alpha+ \delta_0, \beta-\delta_0)$. By (\ref{eqn: Ah''_EL equation}),
\begin{eqnarray}\label{eqn: Thm 6 h'' equation}
     Ah''(x) = \int_{x_0}^x (Bh'(s)-F_m(s)) ds -m_1(x-x_0)+ A(x_0)h''(x_0)
\end{eqnarray}
where $A, B, F_m$ are defined in (\ref{eqn: def A(X)}), (\ref{eqn: def B(x)}) and (\ref{eqn: def F_m(x)}). Define 
\begin{equation}\label{eqn: def F(x)}
    F(x) := \int_{x_0}^x f(s) ds.
\end{equation}
Then, $F_m(x)= F(x) - m(x-x_0)$. 
\paragraph{} From the proof of Theorem \ref{thm: 4_EL Theorem}, we know that $A \in C^1((\alpha, \beta)), B \in C^0((\alpha, \beta))$ and $h'' \in C^1((\alpha, \beta))$. Differentiating (\ref{eqn: Thm 6 h'' equation}), we have 
\begin{equation}\label{eqn: Ah''' equation}
    A(x) h'''(x) = -A'(x) h''(x) + B(x)h'(x) -F_m(x) - m_1
\end{equation}
By (\ref{eqn: def A(X)}),
\begin{equation}\label{eqn: A'(x)}
    A'(x) = -5 \nu_0 \frac{h''(x)h'(x)}{(1+(h'(x))^2)^{7/2}}
\end{equation} 
Hence, by (\ref{eqn: def A(X)}), (\ref{eqn: def B(x)}), (\ref{eqn: h'' integral bound_EL section}) and (\ref{eqn: A'(x)}), there exists $C>0$ such that 
\begin{equation}\label{eqn: upper bound on H^1 norm of A and L^1 norm of B}
    \|A\|_{H^1((\alpha, \beta))}+\|B\|_{L^1((\alpha, \beta))} \leq C \leq C B_\tau (h,h_0, \alpha, \alpha_0, \beta, \beta_0). 
\end{equation}
Now, we need to estimate the $L^\infty$ bound on $F_m$. For this, we obtain $L^1$ bounds on $f$ and then estimate $m$ in (\ref{eqn: def m EL section}). By (\ref{eqn: def B_tau}),
\begin{equation} \label{eqn: h-h_0 integral C B_tau bound}
    \frac{1}{\tau} \int_\alpha^\beta |\htil -h_0| \chi_{[\alpha_0, \beta_0]} dx \leq \frac{1}{\tau} \int_{\alpha_0}^{\beta_0} |\htil -h_0| dx \leq B_\tau (h,h_0, \alpha, \alpha_0, \beta, \beta_0).
\end{equation}
Hence, by (\ref{eqn: def f(x)_EL}), (\ref{eqn: Integral of W B_tau bound}) and (\ref{eqn: def F(x)}), we have 
\begin{equation}\label{eqn: upper bound on L^infinity bound of F}
    \|F\|_{L^\infty((\alpha, \beta))} \leq C B_\tau (h,h_0, \alpha, \alpha_0, \beta, \beta_0).
\end{equation}
By (\ref{assumption: Lip h< L_0 and h>0}), (\ref{eqn: b_0-a_0 lower bound}), (\ref{eqn: phi_0 integral and L^infinity bound})-(\ref{eqn: def m EL section}), (\ref{eqn: h'' integral bound_EL section}), (\ref{eqn: Integral of W B_tau bound}), (\ref{eqn: upper bound on H^1 norm of A and L^1 norm of B}), (\ref{eqn: h-h_0 integral C B_tau bound})
 and H\"older's inequality, 
 \begin{gather}
     |m| \leq  C_0 \bigg(\frac{2L_0}{A_0}\bigg)^{3/2}\|A\|_{L^2((\alpha, \beta))}\|h''\|_{L^2((\alpha, \beta))}+\frac{2C_0L_0}{A_0}\|B\|_{L^1((\alpha, \beta))} 
    \nonumber \\ + \sqrt{\frac{2L_0}{A_0}} \frac{1}{\tau} \int_\alpha^\beta |\htil -h_0| \chi_{[\alpha_0, \beta_0]} dx + \sqrt{\frac{2L_0}{A_0}}\int_\alpha^\beta W(Eu(x, h(x))) dx
    \nonumber\\ \leq CB_\tau (h,h_0, \alpha, \alpha_0, \beta, \beta_0).\label{eqn: m B_tau bound}
 \end{gather}
Therefore, by (\ref{eqn: Upper bound on beta-alpha}), (\ref{eqn: upper bound on L^infinity bound of F}) and (\ref{eqn: m B_tau bound}), we have 
\begin{eqnarray}
    \|F_m\|_{L^\infty((\alpha, \beta))} \leq C B_\tau (h,h_0, \alpha, \alpha_0, \beta, \beta_0), \label{eqn: F_m L^infinity bound} \\
    \|F_m\|_{L^1((\alpha, \beta))} \leq C B_\tau (h,h_0, \alpha, \alpha_0, \beta, \beta_0).\label{eqn: F_m L^1 bound}
\end{eqnarray}

 Now, by (\ref{eqn: phi_0 integral and L^infinity bound}), (\ref{eqn: def m_1}), (\ref{eqn: h'' integral bound_EL section}), (\ref{eqn: upper bound on H^1 norm of A and L^1 norm of B}), (\ref{eqn: F_m L^1 bound}) and H\"older's inequality, we have
\begin{eqnarray}
    |m_1| &\leq& \|\phi_0\|_{L^\infty(\R)}\int_\alpha^\beta |A(x)h''(x) + B(x) h'(x)-F_m(x)| dx \nonumber\\
    & \leq & \sqrt{\frac{2L_0}{A_0}}(\|A\|_{L^2((\alpha, \beta))}\|h''\|_{L^2((\alpha, \beta))}+ L_0\|B\|_{L^1((\alpha, \beta))}+ \|F_m\|_{L^1((\alpha, \beta))}) \nonumber\\
    &\leq& CB_\tau (h,h_0, \alpha, \alpha_0, \beta, \beta_0), \label{eqn: m_1 bound}
\end{eqnarray}
where $\phi_0$ is defined in (\ref{eqn: def phi_0 from psi_0}).
\paragraph{} Note that $A\geq \frac{\nu_0}{(1+L_0^2)^{1/2}}$ by (\ref{eqn: def A(X)}). Therefore, by (\ref{eqn: Ah''' equation}) and (\ref{eqn: A'(x)}), we have 
\begin{equation}\label{eqn: |h'''(x)| bound}
    |h'''(x)| \leq C |h''(x)|^2+|B(x)|+|F_m(x)|+m_1 
\end{equation}
for $x \in (\alpha, \beta)$. Hence, by (\ref{eqn: h'' integral bound_EL section}), (\ref{eqn: Upper bound on beta-alpha}), (\ref{eqn: upper bound on H^1 norm of A and L^1 norm of B}), (\ref{eqn: F_m L^1 bound}) and (\ref{eqn: m_1 bound}),
\begin{equation*}
    \|h'''\|_{L^1((\alpha, \beta))} \leq CB_\tau (h,h_0, \alpha, \alpha_0, \beta, \beta_0),
\end{equation*}
thus proving (\ref{eqn: h''' B_tau L 1 bound}). Recall that $h''$ is continuous in $(\alpha, \beta)$. For every $ x \in (\alpha, \beta)$, by the fundamental theorem of calculus and mean value theorem,
\begin{equation}
    |h''(x)| \leq |\beta-\alpha|^{-1} \|h''\|_{L^1((\alpha, \beta))}+ \|h'''\|_{L^1((\alpha, \beta))}.
\end{equation}
Hence, by (\ref{eqn: beta-alpha lower bound lemma 1}), (\ref{eqn: h'' integral bound_EL section}), (\ref{eqn: h''' B_tau L 1 bound}) and H\"older's inequality, 
\begin{equation*}
    \|h''\|_{L^\infty((\alpha, \beta))} \leq CB_\tau (h,h_0, \alpha, \alpha_0, \beta, \beta_0),
\end{equation*}
thus proving (\ref{eqn: h'' B_tau L infinity bound}).

\end{proof}
\section{Corner of $\Omega_h$ to Triangular domain}
Consider the triangle 
\begin{equation}\label{eqn: def triangle A^l_r}
    A^l_r := \{(x,y) \in \R^2: 0<x<r,\, 0<y<lx\}
\end{equation}
for some $r>0$ and $l>0$. The aim of this section is to flatten the graph of $h$ in the corners of $\Omega_h$ near $(\alpha,0)$ to $A^l_r$. The estimates developed in this section are crucial in dealing with the elasticity boundary value problem in the corners of $\Omega_h$, when combined with the existing regularity results in domains of the form of $A^l_r$. Fix $0<\eta_0<1$ and $M>1$. Assume that
\begin{gather}
    \alpha<\beta, \quad h \in W^{3,1}((\alpha, \beta))\label{eqnar:assumptions1},\\
h>0 \text{ in } (\alpha, \beta),\quad  h(\alpha)=
h(\beta)=0,\label{eqnar:assumptions2}\\
h'(\alpha)\geq 2 \eta_0,\quad h'(\beta)\leq -2 \eta_0,\quad\text{ Lip }h \leq L_0,\label{eqnar:assumptions3}\\
\int_\alpha^\beta |h''(x)|^2 dx \leq M.\label{eqnar:assumptions4}
\end{gather}

Without loss of generality, assume $\alpha=0$. Define
\begin{equation}\label{eqn: def h_0', h_0''}
    h_0':=h'(0), \quad h_0'':= h''(0).
\end{equation}
Given $r>0$, define $I_r := (0,r)$ and $\Omega_h^{0,r}:= \Omega_h \cap (I_r \times \R)$. Choose $r$ such that
\begin{equation}\label{eqn: assume bound on r}
    0<r\leq \delta_0 = \frac{\eta_0^2}{4M}< \frac{1}{4}.
\end{equation}
For $x \in I_r$, define 
\begin{equation}\label{eqn: def sigma}
    \sigma (x):= \frac{h_0'x}{h(x)}.
\end{equation}
Using $\sigma,$ we define the diffeomorphisms $\Phi: I_r \times \R \rightarrow I_r \times \R$ and $\Psi: I_r \times \R \rightarrow I_r \times \R$ by
\begin{equation}
    \Phi(x,y): (x, y/\sigma(x))\quad \Psi(x,y):= (x, y \sigma(x)).
\end{equation}
For the rest of this section, set $l:=h_0'$. Then,
\begin{equation}\label{eqn: define Phi, Psi}
    \Phi(A^l_r)= \Omega_h^{0,r}\, \text{and}\, \Psi (\Omega_h^{0,r})= A^l_r.
\end{equation}
By (\ref{eqnar:assumptions1}), (\ref{eqn: def sigma}) and direct computation, we have that $\sigma \in C^2(\overline{I}_r)$ and 
\begin{eqnarray}
    \sigma'(x)&=& \frac{h_0'h(x)-h_0'xh'(x)}{(h(x))^2}, \label{eqn: sigma'}\\
    \sigma''(x)&=& - \frac{h_0'xh''(x) h(x)+2 h_0'h(x)h'(x)-2h_0'x(h'(x))^2}{(h(x))^3}\label{eqn: sigma''}
\end{eqnarray}
for $x \in I_r$, and 
\begin{equation}
    \sigma(0)=1,\, \sigma'(0)=-\frac{h_0''}{2h_0'},\, \sigma''(0)= \frac{(h_0'')^2}{2(h_0')^2}-\frac{h_0'''}{3h_0'}.
\end{equation}
This implies that
\begin{equation}\label{eqn: C^2 regularity of Phi, Psi}
    \Phi \in C^2(\overline{A}^l_r; \R^2) \text { and } \Psi \in C^2(\overline{\Omega}_h^{0,r}; \R^2).
\end{equation}
Now, we state a result that follows from Lemma 4.3 and Remark 4.4 in \cite{dal2025motion}.
\begin{lem} \label{lem: 9}
    Under the assumptions (\ref{eqnar:assumptions1}) to (\ref{eqnar:assumptions4}), let $r$ be as in $\ref{eqn: assume bound on r}$ and let $p\geq1$.
    \begin{itemize}
        \item[(i)]If $f \in \Omega_h^{0,r}$ and $w \in W^{2,p}(\Omega_h^{0,r})$, then $f \circ \Phi \in L^p(A^l_r)$ and $w \circ \Phi \in W^{2,p}(A^l_r)$. Moreover, the following estimates hold:
    \begin{eqnarray}
        \|f\circ \Phi\|_{L^p(A^l_r)} &\leq& C_p \|f\|_{L^p(\Omega_h^{0,r})}, \label{eqn: f o Phi L^p_0 bound independent of h}\\
        \|\nabla (w \circ \Phi)\|_{L^p(A^l_r)} &\leq& C_p(1+ \sup_{(x,y)\in A^l_r}|y \sigma'(x)|) \|\nabla w\|_{L^p(\Omega_h^{0,r})}, \label{eqn: nabla w o Phi L^p_0 bound independent of h}\\
        \|\nabla^2( w \circ \Phi)\|_{L^p(A^l_r)} &\leq& C_p(1+ \sup_{(x,y)\in A^l_r}|y \sigma'(x)|^2) \|\nabla^2 w\|_{L^p(\Omega_h^{0,r})} \nonumber \\
        & &+C_p(1+\sup_{x \in I_r}| \sigma'(x)|+\sup_{(x,y)\in A^l_r}|y (\sigma'(x))^2|\nonumber \\
        & &+\sup_{(x,y)\in A^l_r}|y \sigma''(x)| )\|\nabla w\|_{L^p(\Omega_h^{0,r})}
    \end{eqnarray}
    where $C_p$ is a positive constant that depends only on $p$. 
    \item[(ii)] If $f \in L^p(A^l_r)$ and $w \in W^{2,p}(A^l_r)$, then $f \circ \Psi \in \Omega_h^{0,r}$ and $w \circ \Psi \in W^{2,p}(\Omega_h^{0,r})$ and the following estimates hold:
    \begin{eqnarray}
        \|f\circ \Psi\|_{L^p(\Omega_h^{0,r})} &\leq& C_p \|f\|_{L^p(A^l_r)}, \label{eqn: f o Psi L^p_0 bound independent of h}\\
        \|\nabla (w \circ \Psi)\|_{L^p(\Omega_h^{0,r})} &\leq& C_p(1+ \sup_{(x,y)\in \Omega_h^{0,r}}|y \sigma'(x)|) \|\nabla w\|_{L^p(A^l_r)}, \label{eqn: nabla w o Psi L^p_0 bound independent of h}\\
        \|\nabla^2( w \circ \Psi)\|_{L^p(\Omega_h^{0,r})} &\leq& C_p(1+ \sup_{(x,y)\in \Omega_h^{0,r}}|y \sigma'(x)|^2) \|\nabla^2 w\|_{L^p(A^l_r)} \nonumber \\
        & &+C_p(1+\sup_{x \in I_r}| \sigma'(x)|+\sup_{(x,y)\in \Omega_h^{0,r}}|y (\sigma'(x))^2|\nonumber \\
        & &+\sup_{(x,y)\in \Omega_h^{0,r}}|y \sigma''(x)| )\|\nabla w\|_{L^p(A^l_r)},
    \end{eqnarray}
     where $C_p$ is a positive constant that depends only on $p$. 
    \end{itemize}
\end{lem}
From the lemma above, it is clear that regularity estimates on $\Omega_h^{0,r}$ are dependent on those in $A^l_r$ and vice versa, and the bounds on $\sigma$ and its derivatives are crucial in those estimates. The following lemma develops all these bounds in terms of $h$ and its derivatives, where (\ref{eqnar:assumptions1}) is assumed, unlike \cite{dal2025motion}, where $W^{4,1}$ regularity of $h$ is used.
\begin{lem}\label{lem: 8} Under the assumptions (\ref{eqnar:assumptions1})-(\ref{eqnar:assumptions4}), let $r$ be as in (\ref{eqn: assume bound on r}). Then, there exists a constant $C=C(\eta_0, M)>0$ independent of $r$ such that
  
    \begin{eqnarray}
     \|\sigma-1\|_{L^\infty(I_r)} &\leq& \frac{r^{1/2}}{\eta_0}\|h''\|_{L^2(I_r)} \leq Cr\|h''\|_{L^\infty(I_r)}, \label{eqn: sigma-1 bound}\\
        \|\sigma'\|_{L^\infty(I_r)} &\leq& C(|h_0''|+\|h'''\|_{L^1(I_r)}),\label{eqn: sigma' bound}\\
        \sup_{(x,y)\in \Omega_h^{0,r}}|y \sigma'(x)| &\leq& Cr(|h_0''|+\|h'''\|_{L^1(I_r)}),\label{eqn: ysigma' bound_omega} \\ 
        \sup_{(x,y)\in A^l_r}|y \sigma'(x)| &\leq& Cr(|h_0''|+\|h'''\|_{L^1(I_r)}),\label{eqn: ysigma' bound_triangle}\\ 
        \sup_{(x,y)\in \Omega_h^{0,r}}|y \sigma''(x)|&\leq& C(r|h_0''|^2+\|h'''\|_{L^1(I_r)}+r\|h'''\|_{L^1(I_r)}^2), \label{eqn: ysigma'' bound_omega}\\
        \sup_{(x,y)\in A^l_r}|y \sigma''(x)|&\leq& C(r|h_0''|^2+\|h'''\|_{L^1(I_r)}+r\|h'''\|_{L^1(I_r)}^2),\label{eqn: ysigma'' bound_triangle}\\
        \|h\sigma''\|_{L^\infty(I_r)}&\leq& C(r|h_0''|^2+\|h'''\|_{L^1(I_r)}+r\|h'''\|_{L^1(I_r)}^2).\label{eqn: hsigma'' bound}
   \end{eqnarray}
\end{lem}
\begin{proof}
    For $x \in I_r$, using Taylor's formula with integral remainder, we have
    \begin{eqnarray*}
        \sigma(x)-1&=& \frac{h_0'x-h(x)}{h(x)}= \frac{h_0'x-(h_0'x+\int_0^x h''(s)(x-s)ds)}{h(x)}\\
        &=& -\frac{\int_0^x h''(s)(x-s)ds}{h(x)}.
    \end{eqnarray*}
Now, using the theorem hypotheses, we obtain a lower bound for $h(x)$ in terms of $x$. By H{\"o}lder's inequality, (\ref{eqnar:assumptions4}) and (\ref{eqn: assume bound on r}), for $x \in I_r$,
\begin{eqnarray*}
    \bigg| \int_0^x h''(s) ds \bigg| &\leq & \bigg(\int_0^x |h''(s)|^2ds\bigg)^{1/2} \bigg(\int_0^x1 ds\bigg)^{1/2}\\
    &\leq& M^{1/2}x^{1/2}\leq M^{1/2}\frac{\eta_0}{2M^{1/2}} =\frac{\eta_0}{2}< \eta_0.
\end{eqnarray*}
Therefore, by the fundamental theorem of calculus and (\ref{eqnar:assumptions3}),
\begin{equation*}
    h'(x)=h_0'+\int_0^x h''(s) ds \geq 2\eta_0-\eta_0 \geq \eta_0.
\end{equation*}
By using the fundamental theorem of calculus again and the fact that $h_0=0$,
\begin{equation}\label{eqn:Lower bound on h}
    h(x)= \int_0^x h'(s) ds \geq \eta_0 x.
\end{equation}
Therefore, by (\ref{eqn:Lower bound on h}) and H{\"o}lder's inequality, for every $x \in I_r$, 
\begin{eqnarray*}
    | \sigma(x)-1| &=& \bigg| \frac{h_0'x-h(x)}{h(x)}\bigg| \leq \frac{1}{\eta_0 x}\bigg(\int_0^x |h''(s)|^2 ds\bigg)^{1/2}\bigg(\int_0^x (x-s)ds\bigg)^{1/2}\\
    &\leq& \frac{1}{\eta_0 x}\frac{x^{3/2}}{\sqrt{3}}\bigg(\int_0^x |h''(s)|^2 ds\bigg)^{1/2}  \leq\frac{r^{1/2}}{\eta_0}\|h''\|_{L^2(I_r)} \leq \frac{r}{\eta_0}\|h''\|_{L^\infty(I_r)}.
\end{eqnarray*}
Hence (\ref{eqn: sigma-1 bound}) follows.
\paragraph{} To estimate the numerator in (\ref{eqn: sigma'}), we use Taylor's formula with integral remainder again. As $h_0=h(0)=0$, for $x \in I_r$,
\begin{eqnarray*}
    h_0'h(x)-h_0'xh'(x)&=& h_0'\bigg(h_0'x+h_0'' \frac{x^2}{2}+\int_0^x h'''(s)\frac{(x-s)^2}{2} ds\bigg)\\
    & &-h_0'x\bigg(h_0'+h_0''x+\int_0^x h'''(s)(x-s)ds\bigg)\\
    &=& (h_0')^2 x + h_0'h_0''\frac{x^2}{2}+ \frac{h_0'}{2}\int_0^x h'''(s)(x-s)^2ds\\
    & & -(h_0')^2 x-h_0'h_0''x^2-h_0'x\int_0^x h'''(s)(x-s)ds\\
    &=& -h_0'h_0''\frac{x^2}{2}+\frac{h_0'}{2}\int_0^x h'''(s)(x-s)^2ds-h_0'x\int_0^x h'''(s)(x-s)ds.
\end{eqnarray*}
Therefore,
\begin{eqnarray*}
    | h_0'h(x)-h_0'xh'(x)| &\leq& \frac{1}{2}h_0' |h_0''|x^2+\frac{3}{2} h_0' x^2\int_0^x |h'''(s)| ds.
\end{eqnarray*}
Hence, by (\ref{eqn: sigma'}) and (\ref{eqn:Lower bound on h}),
\begin{eqnarray*}
    \|\sigma'\|_{L^\infty(I_r)} &\leq& \frac{1}{(h(x))^2} \bigg(\frac{1}{2}h_0' |h_0''|x^2+\frac{3}{2} h_0'x^2 \int_0^x |h'''(s)| ds \bigg)\\
    &\leq& \frac{L_0}{\eta_0^2}\bigg(\frac{1}{2} |h_0''|+\frac{3}{2} \|h'''\|_{L^1(I_r)} \bigg),
\end{eqnarray*}
which proves (\ref{eqn: sigma' bound}).
\paragraph{}Note that $0< y < L_0r$ when $(x,y) \in \Omega_h^{0,r}$ as $h$ is Lipschitz continuous with Lipschitz constant less than or equal to $L_0$ and $0<y<lr \leq L_0r$ when $(x,y) \in A^l_r$ as $l=h_0' \leq L_0$. Hence, (\ref{eqn: ysigma' bound_omega}) and (\ref{eqn: ysigma' bound_triangle}) follows from (\ref{eqn: sigma' bound}).
\paragraph{} To estimate the numerator in (\ref{eqn: sigma''}), use Taylor's formula with integral remainder again to obtain
\begin{eqnarray*}
    h_0'xh''(x)h(x)&=&h_0'x\bigg(h_0''+\int_0^x h'''(s)ds\bigg)\bigg(h_0'x+\frac{1}{2}h_0''x^2+\frac{1}{2}\int_0^x h'''(s)(x-s)^2 ds\bigg),\\
    2h_0'h(x)h'(x) &=& 2h_0'\bigg(h_0'x+\frac{1}{2}h_0''x^2+\frac{1}{2}\int_0^x h'''(s)(x-s)^2 \bigg)\bigg(h_0'+h_0''x+\int_0^x h'''(s)(x-s)ds\bigg),\\
    2h_0'x(h'(x))^2&=&2h_0'x\bigg(h_0'+h_0''x+\int_0^x h'''(s)(x-s)ds\bigg)^2.
\end{eqnarray*}
Hence,
\begin{eqnarray*}
     -\sigma''(x)(h(x))^3 &=&-\frac{1}{2}h_0'(h_0'')^2x^3+\frac{3}{2}h_0'h_0''x\int_0^x h'''(s)(x-s)^2+(h_0')^2x^2 \int_0^x h'''(s)ds\\
     &&+\frac{1}{2}h_0'h_0''x^3\int_0^x h'''(s)ds+\frac{1}{2}h_0'x\int_0^x h'''(s)ds\int_0^x h'''(s)(x-s)^2ds\\
     &&-2(h_0')^2x\int_0^x h'''(s)(x-s)ds-3h_0'h_0''x^2\int_0^x h'''(s)(x-s)ds\\
     &&+(h_0')^2\int_0^x h'''(s)(x-s)^2ds+ h_0'\int_0^x h'''(s)(x-s)^2ds\int_0^x h'''(s)(x-s)ds\\
     &&-2h_0'x\int_0^x h'''(s)(x-s)ds.
\end{eqnarray*}
Therefore, using Young's inequality for $p=q=2$ and (\ref{eqn:Lower bound on h}),
\begin{eqnarray*}
    \sup_{(x,y) \in \Omega_h^{0,r}} |y \sigma''(x)|&\leq& \frac{C}{(h(x))^3}(r|h_0''|^2 x^3+r|h_0''|\|h'''\|_{L^1(I_r)}x^3+\|h'''\|_{L^1(I_r)}x^3+r\|h'''\|_{L^1(I_r)}^2x^3)\\
    &\leq&\frac{C}{\eta_0^3 x^3}(r|h_0''|^2 x^3+\|h'''\|_{L^1(I_r)}x^3+r\|h'''\|_{L^1(I_r)}^2x^3)\\
    &=& C(r|h_0''|^2 +\|h'''\|_{L^1(I_r)}+r\|h'''\|_{L^1(I_r)}^2),
\end{eqnarray*}
which proves (\ref{eqn: ysigma'' bound_omega}). Following the steps above, (\ref{eqn: ysigma'' bound_triangle}) and (\ref{eqn: hsigma'' bound}) can also be proved.

\end{proof}
\begin{remark}
    By (\ref{eqnar:assumptions4}), (\ref{eqn: assume bound on r}) and (\ref{eqn: sigma-1 bound}),
    \begin{equation*}
        \|\sigma-1\|_{L^\infty(I_r)} \leq 1/2.
    \end{equation*}
    This, in turn, would imply that 
    \begin{gather}
       1/2 \leq \|\sigma \|_{L^\infty(I_r)} \leq 3/2, \text{ and } \label{eqn: sigma upper and lower l^infinity bound}\\
        \|\sigma +1\|_{L^\infty(I_r)} \leq 5/2. \label{eqn: sigma +1 L^infinity bound}
    \end{gather} 
\end{remark}

Now, we will focus on the boundaries of $\Omega_h^{0,r}$ and $\A^l_r$. Let $\Gamma_h^{0,r}:=\Gamma_h \cap (I_r \times \R)$. Define
\begin{equation}\label{eqn: def hypotenuse of triangle A^l_r}
    \Gamma^l_r := \{(x,lx): 0<x<r\} \subset \partial A^l_r, 
\end{equation}
where $r>0$ is as in (\ref{eqn: assume bound on r}). We denote the outer unit normal to $\Omega_h^{0,r}$ at $\Gamma_h^{0,r}$ as $\nu^h:=(\nu^h_1, \nu^h_2)$ and that to $A^l_r$ at $\Gamma^l_r$ as $\nu^0:=(\nu^0_1, \nu^0_2)$. That is,
\begin{eqnarray*}
   \nu^0= \frac{(-h_0',1)}{\sqrt{1+(h_0')^2}}, \quad \nu^h(x)=\frac{(-h'(x),1)}{\sqrt{1+(h'(x))^2}},\quad \text{for } x\in I_r.
\end{eqnarray*}
We define the functions
\begin{eqnarray}
    \omega_1(x)&:=&\nu^0_1-\nu^h_1(x)=-\frac{h_0'}{\sqrt{1+(h_0')^2}}+\frac{h'(x)}{\sqrt{1+(h'(x))^2}},\nonumber\\
    \omega_2(x)&:=& \nu^0_2-\nu^h_2(x)=\frac{1}{\sqrt{1+(h_0')^2}}-\frac{1}{\sqrt{1+(h'(x))^2}},\nonumber\\
    \omega_3(x)&:=& -\sigma'(x)h(x)\nu^h_1(x)=\sigma'(x)h(x)\frac{h'(x)}{\sqrt{1+(h'(x))^2}},\nonumber\\
    \omega_4(x)&:=& -(\sigma(x)-1)\nu^h_1(x)=(\sigma(x)-1)\frac{h'(x)}{\sqrt{1+(h'(x))^2}},\label{eqn: def omega_i}\\
    \omega_5(x)&:=& \sigma'(x)h(x)\nu^h_2(x) =\sigma'(x)h(x)\frac{1}{\sqrt{1+(h'(x))^2}},\nonumber\\
    \omega_6(x)&:=& (\sigma(x)-1)\nu^h_2(x)=(\sigma(x)-1)\frac{1}{\sqrt{1+(h'(x))^2}}.\nonumber
\end{eqnarray}
As we will see in the lemma below (Lemma 4.5, \cite{dal2025motion}), these functions are useful in formalizing the change in Dirichlet conditions when we change variables to move from $\Gamma_h^{0,r}$ to $\Gamma^l_r$.
\begin{lem}\label{lem: 10}
  Under the assumptions (\ref{eqnar:assumptions1})-(\ref{eqnar:assumptions4}), let $r$ be as in (\ref{eqn: assume bound on r}), $1<p< \infty$ and $u \in W^{2, p}_{\text{loc}}(\Omega_h^{0,r}; \R^2)$ with $u \in W^{2, p}({\overline{\Omega}_h^{\rho,r}}; \R^2)$ for every $0<\rho<r$. Let $l=h_0'$, $v: A^l_r \rightarrow \R^2$ be defined by $v(x,y):= u\circ \Phi(x,y)$ and let $g \in W^{1,p}(\Omega_h^{0,r}; \R^2)$. Assume that
  \begin{equation}
      (\C Eu)\nu^h= 2 \mu (E u)\nu^h+ \lambda({\operatorname*{div}}\,u )\nu^h = g \text{ on } \Gamma_h^{0,r}.
  \end{equation}
  Then, if $\nu^0$ is the outer unit normal to $\Gamma^l_r$, 
  \begin{equation}
      (\C Ev)\nu^0= 2 \mu (E v)\nu^0+ \lambda({\operatorname*{div}}\,v )\nu^0 = g \circ \Phi+ \hat{g}^v + \check{g}^v\text{ on } \Gamma^l_r,
  \end{equation}
  where $\hat{g}^v=(\hat{g}^v_1, \hat{g}^v_2)$ and $\check{g}^v=(\check{g}^v_1, \check{g}^v_2)$ are defined by
  \begin{eqnarray}
      \hat{g}^v_1&:=& (2 \mu \omega_1+ \lambda \omega_1) \pa_x v_1+ \mu \omega_2 \pa_x v_2, \label{eqn: g hat g check def 1}\\
      \hat{g}^v_2 &:=& \lambda \omega_2 \pa_x v_1+ \mu\omega_1 \pa_x v_2,\\
      \check{g}^v_1&:=&(\mu \omega_2 + 2 \mu \omega_3 + \lambda \omega_3-\mu\omega_6)\pa_y v_1 \nonumber\\
      & & +(\lambda \omega_1+\lambda \omega_4- \mu\omega_5)\pa_y v_2,\\
      \check{g}^v_2&:=& (\mu\omega_1+ \mu \omega_4- \lambda\omega_5) \pa_y v_1 \nonumber\\
      & &+(2 \mu \omega_2+ \lambda \omega_2+ \mu\omega_3- 2 \mu \omega_6-\lambda \omega_6)\pa_y v_2.\label{eqn: g hat g check def end}
  \end{eqnarray}
\end{lem}
\begin{lem}\label{lem: 11}
Under the assumptions (\ref{eqnar:assumptions1})-(\ref{eqnar:assumptions4}), let $r$ be as in (\ref{eqn: assume bound on r}) and let $\omega_i$, $i=1,...,6$ be defined as in (\ref{eqn: def omega_i}). Then, there exists a constant $C=C(\eta_0, M)>0$ such that 
    \begin{eqnarray}
        \|\omega_i\|_{L^{\infty}(I_r)} &\leq& Cr(\|h''\|_{L^{\infty}(I_r)}+\|h'''\|_{L^{1}(I_r)}),\label{eqn: omega_i bound}\\
        \|\omega_i'\|_{L^{\infty}(I_r)} &\leq& C(\|h''\|_{L^{\infty}(I_r)}+\|h'''\|_{L^{1}(I_r)})+Cr(\|h''\|_{L^{\infty}(I_r)}^2+\|h'''\|_{L^{1}(I_r)}^2), \label{eqn: omega'_i bound}
    \end{eqnarray}
    for $i=1,...,6$.
\end{lem}
\begin{proof}
    For $t \in \R$, define $f_1(t):= \frac{t}{\sqrt{1+t^2}}$ and $f_2(t):= \frac{1}{\sqrt{1+t^2}}$. Note that $|f_i(t)| \leq 1$ for $i=1,2$. Now,
    \begin{eqnarray*}
        f_1'(t)=\frac{1}{(1+t^{2})^{3/2}}, \quad  f_2'(t)=\frac{-t}{(1+t^{2})^{3/2}}.
    \end{eqnarray*}
  Clearly, $|f_i'(t)| \leq 1$ for all $t \in \R$ and $i=1,2$. hence, both $f_1$ and $f_2$ are Lipschitz continuous with Lipschitz constant less than or equal to 1. Now, since by (\ref{eqn: def omega_i}), we have
    \begin{eqnarray*}
        \omega_1(x)=f_1(h'(x))-f_1(h_0'),\quad \omega_2(x)=f_2(h_0')-f_2(h'(x))
    \end{eqnarray*}
    for every $x \in I_r$, 
    \begin{eqnarray*}
        |\omega_i(x)| \leq |h'(x)-h_0'|= \bigg| \int_{0}^x h''(s) ds\bigg| \leq r\|h''\|_{L^{\infty}(I_r)},
    \end{eqnarray*}
    for $i=1,2$, where we used the fundamental theorem of calculus. Therefore, (\ref{eqn: omega_i bound}) follows for $i=1,2$. 
    \paragraph{}On the other hand, for every $x \in I_r$, for $i=1,2$,
    \begin{eqnarray*}
        |\omega_1'(x)| = |f_1'(h'(x))h''(x)| \leq |h''(x)|, \quad |\omega_2'(x)| = |f_2'(h'(x))h''(x)| \leq |h''(x)|.
    \end{eqnarray*}
    Therefore, for $i=1,2$, $\|\omega_i'\|_{L^\infty(I_r)}\leq \|h''\|_{L^{\infty}(I_r)}$ and hence (\ref{eqn: omega'_i bound}) follows.
 \paragraph{}Since $h(x) \leq L_0 x \leq L_0 r$ for every $x \in I_r$ by (\ref{eqnar:assumptions2}) and (\ref{eqnar:assumptions3}),
\begin{eqnarray*}
    |\omega_3(x)|=|\sigma'(x) h(x) f_1(h'(x))| \leq L_0 r |\sigma'(x)| \leq Cr(\|h''\|_{L^\infty(I_r)}+\|h'''\|_{L^1(I_r)}),
\end{eqnarray*}
 using (\ref{eqn: sigma' bound}). Hence, (\ref{eqn: omega_i bound}) follows for $i=3$. Also, since $\omega_3=\sigma'hf_1(h')$,
 \begin{eqnarray*}
    | \omega_3'(x)|&=&|\sigma''(x) h(x) f_1(h'(x))+\sigma'(x) h'(x) f_1(h'(x))+\sigma'(x) h(x) f_1'(h'(x))h''(x)|\\
    &\leq& \|\sigma'' h\|_{L^\infty(I_r)}+L_0 \|\sigma'\|_{L^\infty(I_r)}+L_0r\|\sigma'\|_{L^\infty(I_r)}\| h''\|_{L^\infty(I_r)}\\
    &\leq& C(r|h_0''|^2+\|h'''\|_{L^1(I_r)}+r\|h'''\|_{L^1(I_r)}^2+|h_0''|+\|h'''\|_{L^1(I_r)}+\\
    & &r\|h''\|_{L^\infty(I_r)}(|h_0''|+\|h'''\|_{L^1(I_r)}))\\
    &\leq& C(\|h''\|_{L^\infty(I_r)}+\|h'''\|_{L^1(I_r)})+ Cr(\|h''\|_{L^\infty(I_r)}^2+\|h'''\|_{L^1(I_r)}^2),
 \end{eqnarray*}
 for $C>0$, where we used (\ref{eqn: sigma' bound}), (\ref{eqn: hsigma'' bound}) and Young's inequality. This proves (\ref{eqn: omega'_i bound}) for $i=3$. Since $\omega_5(x)=\sigma'(x)h(x) f_2(h'(x))$, the estimates for $\omega_5$ are similar.\\
 \paragraph{} Now, for $x \in I_r$,
 \begin{eqnarray*}
     |\omega_4(x)|=|(\sigma(x)-1)f_1(h'(x))| \leq |\sigma(x)-1|,
 \end{eqnarray*}
 and hence (\ref{eqn: omega_i bound}) follows from (\ref{eqn: sigma-1 bound}). Further, using (\ref{eqn: sigma-1 bound}) and (\ref{eqn: sigma' bound}), we have
 \begin{eqnarray*}
     |\omega_4'(x)|&=&|\sigma'(x)f_1(h(x)+(\sigma-1)f_1'(h'(x))h''(x)| \leq |\sigma'(x)|+|(\sigma(x)-1)h''(x)|\\
     &\leq& C(\|h''\|_{L^\infty(I_r)}+\|h'''\|_{L^1(I_r)})+ Cr\|h''\|_{L^\infty(I_r)}^2,
 \end{eqnarray*}
 which implies (\ref{eqn: omega'_i bound}) for $i=4$. Since $\omega_6(x)=(\sigma(x)-1)f_2(h'(x))$, the estimates follow similarly.
 
\end{proof}
\section{Regularity at Corners}
In this section, we study the regularity of solutions to the boundary value problem \ref{eqn: u BVP EL eqn} near the corners of $\Omega_h$. We consider the following boundary value problem: 
\begin{equation}\label{eqn: regularity BVP}
      \begin{cases}
          -{\operatorname*{div}}\,(\C Ew)= f  \text{ in } A^l_r,\\
        (\C Ew )\nu^l=g  \text{ on } \Gamma^l_r,\\
        w=0 \text{ on } \pa A^l_r \setminus \Gamma^l_r,
      \end{cases}
  \end{equation}
We will rely on the following theorem which was proved in \cite[Theorem 5.2]{dal2025motion}. We recall that $A^l_r, \Gamma^l_r$ are as defined in (\ref{eqn: def triangle A^l_r}), (\ref{eqn: def hypotenuse of triangle A^l_r}) respectively.
\begin{thm}\label{thm: 13, elliptic regularity in triangle}
    There exists $p_0\in (4/3, 2)$ depending on $\lambda$ and $\mu$ with the property that if $r>0$, $0< \eta_0\leq l\leq L_0$, $f \in L^{p_0}(A^l_r ; \R^2)$, and $g \in W^{1,p_0}(A^l_r ; \R^2)$ with $g=0$ on one of the sides of the triangle $A^l_r$ different from $\Gamma^l_r$, and if $w \in H^1(A^l_r; \R^2)$ is the unique weak solution to the boundary value problem (\ref{eqn: regularity BVP}), then $w \in W^{2,p_0}(A^l_r; \R^2)$ and 
  \begin{equation}\label{eqn: elliptic regularity bound on second derivative}
      \|\nabla^2 w\|_{L^{p_0}(A^l_r; \R^2)} \leq \kappa (\|f\|_{L^{p_0}(A^l_r; \R^2)}+ \|\nabla g\|_{L^{p_0}(A^l_r; \R^2)}),
  \end{equation}
  for a constant $\kappa>0$ depending on $\lambda, \mu, \eta_0$ and $L_0$ but independent of $r,l,f$ and $g$.
    
\end{thm}

Using the result above, we prove the $W^{2,p_0}$ regularity of the elastic displacement $u$ throughout $\Omega_h$.
\begin{thm}\label{thm: 14}
    Under the assumptions of Theorem \ref{thm: 6}, we have that $u \in W^{2, p_0}(\Omega_h)$.
\end{thm}
To prove this theorem, we use the following Poincar\'e's inequality from Lemma 6.2 and Remark 6.3 in \cite{dal2025motion}.
\begin{lem}\label{lemma: 15, Poincare's inequality}
Let $0< r< \beta- \alpha$, $0< \eta_0 \leq m \leq L_0$ and let $p\geq1$.
\begin{itemize}
    \item[(i)] Then,
    \begin{equation}\label{eqn: Pincare Omega}
        \|v\|_{L^p(\Omega_h^{\alpha, \alpha+r})} \leq L_0r\|\nabla v\|_{L^p(\Omega_h^{\alpha, \alpha+r})}
    \end{equation}
    for every $v \in W^{1,p}(\Omega_h^{\alpha, \alpha+r})$ such that $v(x,0)=0$ for $x \in (\alpha, \alpha+r)$ (in the sense of traces).
    \item[(ii)] Also, 
    \begin{equation}\label{eqn: Poincare triangle}
        \|v\|_{L^p(A^l_r)} \leq \max \{L_0, 1/\eta_0\}r \|\nabla v\|_{L^p(A^l_r)}
    \end{equation}
    for every $v \in W^{1,p}(A^l_r)$ such that $v(x,0)=0$ for $x \in (0,r)$ or $v(r,y)=0$ for $y \in (0, lr)$ (in the sense of traces).
\end{itemize}
\end{lem}
Now, we prove theorem \ref{thm: 14}.
\begin{proof}[Proof of Theorem \ref{thm: 14}]
    By Theorem \ref{thm: 4_EL Theorem}, we have that $u \in C^{2, 1/2}(\overline{\Omega}_h^{a,b};\R^2)$ for every $\alpha<a<b< \beta$. Hence, it suffices to prove that there exists $r>0$ sufficiently small such that $u \in W^{2, p_0}(\Omega_h^{\alpha, \alpha+r})\cup W^{2, p_0}(\Omega_h^{\beta-r, \beta})$
    where $\Omega_h^{a,b}$ is defined in (\ref{eqn: Part of Omega_h}). We prove that $u \in W^{2, p_0}(\Omega_h^{\alpha, \alpha+r})$ and the result in the other corner follows in the same way.
    \paragraph{} Without loss of generality, we assume that $\alpha=0$.\\\\
    \textbf{Step 1: Localization.}  For $(x,y) \in \R^2$, define 
    \begin{equation}\label{eqn: def w_0}
        w_0(x,y):= (e_0x,0).
    \end{equation}
 For every $r>0$, let $\varphi_r \in C^\infty(\R^2)$ be such that $0\leq \varphi_r\leq 1$, $\varphi(x,y)=1$ for $x \leq 5r/8$, $\varphi(x,y)=0$ for $x \geq 7r/8$, $\|\nabla \varphi_r\|_{L^\infty(\R^2)}\leq C/r$ and $\|\nabla^2 \varphi_r\|_{L^\infty(\R^2)}\leq C/r^2$, where $C>0$ is a constant independent of $r$. 
 \paragraph{} Define $\tilde{u}$ as 
 \begin{equation}\label{eqn: def u-tilde}
     \tilde{u}:= (u-w_0)\varphi_r
 \end{equation}
 in $\Omega_h$. Let $0<r<\beta$. Then, by (\ref{eqn: def w_0}) and (\ref{eqn: def u-tilde}), $\util(x,0)=0$ for $0\leq x\leq r$ and $\util(r,y)=0$ for $0\leq y\leq h(r)$, that is, $\util=0$ on $\pa \Omega_h^{0,r} \setminus \Gamma_h^{0,r}$, where $\Gamma_h^{0,r}:=\{(x,h(x)): 0<x<r\}$. Further, 
 \begin{equation*}
     \nabla \util = \varphi_r \nabla (u-w_0)+ (u-w_0)\otimes \nabla \varphi_r.
 \end{equation*}
 Since  $\|\nabla \varphi_r\|_{L^\infty(\R^2)}\leq C/r$, by Lemma \ref{lemma: 15, Poincare's inequality}, we have 
 \begin{equation}
     \|\nabla \util\|_{L^2(\Omega_h^{0,r})} \leq C\|\nabla(u-w_0)\|_{L^2(\Omega_h^{0,r})}.
 \end{equation}
By Remark \ref{remark: 1}, we have that
\begin{equation}\label{eqn: elliptic PDE, regularity at corners}
    \begin{cases}
        -{\operatorname*{div}}\, \C E \util= -\mu \Delta \util -(\lambda+\mu) \nabla{\operatorname*{div}}\,\util=f \text{ in } \Omega_h^{0,r},\\
        \C E\util \nu^h= 2 \mu (E \util)\nu^h+ \lambda({\operatorname*{div}}\,\util )\nu^h=g  \text{ on } \Gamma_h^{0,r},\\
        \util=0 \text{ on } \pa \Omega_h^{0,r}\setminus \Gamma_h^{0,r},
    \end{cases}
\end{equation}
where $\nu^h(x)$ is the outward unit normal to $\Gamma_h$ and 
\begin{eqnarray}
     f&:= &-\mu(2 (\nabla u-\nabla w_0)\nabla \varphi_r+(u-w_0)\Delta \varphi_r) \nonumber\\
     & &-(\lambda+\mu)(\text{div }u- \text{ div }w_0)\nabla \varphi_r + (\nabla u-\nabla w_0)^T \nabla \varphi_r+ \nabla^2 \varphi_r(u-w_0), \label{eqn: def f_u regularity proof}\\
     g&:=& \mu((u-w_0)\otimes \nabla \varphi_r+ \nabla \varphi_r \otimes (u-w_0))\nu^h+ \lambda \,\text{trace}((u-w_0)\otimes \nabla \varphi_r)\nu^h. \label{eqn: def g_u regularity proof}
\end{eqnarray}
 Note that $f \in L^2(\Omega_h^{0,r};\R^2)$ by Theorem \ref{thm: 2_existence of minimizer}. Since $h \in C^2((0, \beta))$, $\nu^h$ is a $C^1$ function. Therefore, $g \in H^1(\Omega_h^{0,r};\R^2)$. Note that the trace of $g$ vanishes on $\pa \Omega_h^{0,r}\setminus \Gamma_h^{0,r}$, since $u(x,0)-w_0(x,0)=0$ for a.e. $x \in [0,r]$ and $\varphi_r(r,y)=0$ for all $y \in \R$.  Hence, for every $f\in L^2(\Omega_h^{0,r};\R^2)$ and $g \in H^1(\Omega_h^{0,r};\R^2)$, problem (\ref{eqn: elliptic PDE, regularity at corners}) has a unique weak solution in $H^1(\Omega_h^{0,r};\R^2)$ by Theorem \ref{thm: 13, elliptic regularity in triangle}.\\
 \vspace{0.2 cm}\\
 \textbf{Step 2: Straightening the boundary.} Define the function 
 \begin{equation}\label{eqn: define v using utilde}
     v:= \util \circ \Phi \in H^1(A^l_r; \R^2),
 \end{equation}
where $r$ satisfies (\ref{eqn: assume bound on r}), $m =h'(0)$, $A^l_r$ and $ \Phi$ are defined in (\ref{eqn: def triangle A^l_r}) and (\ref{eqn: define Phi, Psi}) respectively.
\paragraph{} By (\ref{eqn: define v using utilde}), $\util (x,y)= v(x ,\sigma(x)y)$ where $\sigma$ is defined in (\ref{eqn: def sigma}).  Using direct computation, it can be shown that $v$ satisfies the boundary value problem 
\begin{equation}\label{eqn: BVP satisfied by v=utilde o Phi}
    \begin{cases}
        -\operatorname{div }\C Ev= - \mu\Delta v-(\lambda+\mu) \nabla \operatorname{ div }v = f \circ \Phi + f^v \text{ in } A^l_r,\\
        \C Ev \nu^0=2 \mu (E v)\nu^0+ \lambda(\operatorname{ div }v )\nu^0 = g \circ \Phi + \hat{g}^v+ \check{g}^v \text{ on } \Gamma^l_r,\\
        v= 0 \text{ on } \pa A^l_r \setminus \Gamma^l_r,
    \end{cases}
\end{equation}
 where $\nu^0=(\frac{(-h_0',1)}{\sqrt{1+(h_0')^2}})$ is the outward unit normal vector on $\Gamma^l_r$, where $\Gamma^l_r$ is defined in (\ref{eqn: def hypotenuse of triangle A^l_r}), and $f^v:=(f^v_1, f^v_2) \in H^{-1}(A^l_r;\R^2)$ is defined by 
 \begin{eqnarray}
     f^v_1 &:=& \mu \bigg[ \bigg(\sigma^2-1+y^2\frac{(\sigma')^2}{\sigma^2}\bigg)\pa_{yy}^2 v_1+2y \frac{\sigma'}{\sigma}\pa_{xy}^2v_1\bigg]\nonumber \\
     & & +(\lambda+\mu) \bigg[ y^2\frac{(\sigma')^2}{\sigma^2})\pa_{yy}^2 v_1 +2y \frac{\sigma'}{\sigma}\pa_{xy}^2v_1+(\sigma-1)\pa_{xy}^2v_2+y \sigma' \pa_{yy}^2v_2\bigg]\nonumber\\
     & & +\mu y \frac{\sigma''}{\sigma}\pa_{y}v_1+ (\lambda+\mu)\bigg[ y \frac{\sigma''}{\sigma}\pa_{y}v_1+ \sigma' \pa_y v_2\bigg] \label{eqn: def f^v_1}\\
     f^v_2&:=& \mu \bigg[ \bigg(\sigma^2-1+y^2\frac{(\sigma')^2}{\sigma^2}\bigg)\pa_{yy}^2 v_2+2y \frac{\sigma'}{\sigma}\pa_{xy}^2v_2\bigg]\nonumber \\
     & & +(\lambda+\mu) [ (\sigma-1)\pa_{xy}^2v_1+y \sigma' \pa_{yy}^2v_1+(\sigma^2-1)\pa_{yy}^2 v_2 ]\nonumber\\
     & & +\mu y \frac{\sigma''}{\sigma}\pa_{y}v_2 + (\lambda+\mu)\sigma' \pa_y v_1, \label{eqn: def f^v_2}
 \end{eqnarray}
 and $\hat{g}^v:= (\hat{g}^v_1, \hat{g}^v_2) \in W^{1, p_0}(A^l_r; \R^2)$ and $\check{g}^v:= (\check{g}^v_1, \check{g}^v_2) \in W^{1, p_0}(A^l_r; \R^2)$ are defined in (\ref{eqn: g hat g check def 1})- (\ref{eqn: g hat g check def end}).\\\\
 \textbf{Step 3: Fixed Point Argument.} Define the space
 \begin{equation}\label{eqn: def X^l_r}
     X^l_r := \{ z \in W^{2, p_0}(A^l_r; \R^2) : z=0 \text{ on } \pa A^l_r \setminus \Gamma^l_r\}.
 \end{equation}
 For $z \in X^l_r$, define the norm $\|\cdot\|_{X^l_r}$ as
 \begin{equation}
    \|z\|_{X^l_r}:=\|\nabla^2 z \|_{L^{p_0}(A^l_r;\R^2)}. 
 \end{equation}
 We claim that $(X^l_r, \|\cdot\|_{X^l_r})$ is a Banach space where $\|\cdot\|_{X^l_r}$ is equivalent to $\|\cdot\|_{ W^{2, p_0}(A^l_r; \R^2)}$. Note that the existence and completeness of $\|\cdot\|_{X^l_r}$ follows from the equivalence with $\|\cdot\|_{ W^{2, p_0}}$ as $X^l_r$ is a closed set of the Banach space $(W^{2, p_0}(A^l_r; \R^2), \|\cdot\|_{W^{2, p_0}(A^l_r; \R^2)})$. By Lemma \ref{lemma: 15, Poincare's inequality}, there exists a constant $C>0$, depending only on $\eta_0$ and $L_0$, such that for every $\varphi \in W^{1,p_0}(A^l_r ; \R^2)$ with $\varphi=(0,0)$ on one of the sides of the triangle $A^l_r$ other that $\Gamma^l_r$,
 \begin{equation}\label{eqn: Poincare in X^l_r}
     \|\varphi\|_{L^p_0(A^l_r; \R^2)} \leq C r  \|\nabla \varphi\|_{L^p_0(A^l_r; \R^2)}.
 \end{equation}
If $z \in X^l_r$, $\pa_x z (x,0)=(0,0)$ for all $x \in [0,r]$ and $\pa_y z (r,y)=(0,0)$ for all $y \in [0, lr]$. Then, by (\ref{eqn: def X^l_r}) and (\ref{eqn: Poincare in X^l_r}), 
\begin{eqnarray*}
    \|z\|_{L^p_0(A^l_r; \R^2)} &\leq& C r  \|\nabla z\|_{L^p_0(A^l_r; \R^2)},\\
    \|\pa_x z\|_{L^p_0(A^l_r; \R^2)} &\leq& C r  \|\nabla \pa_x z\|_{L^p_0(A^l_r; \R^2)},\\
    \|\pa_y z\|_{L^p_0(A^l_r; \R^2)} &\leq& C r  \|\nabla \pa_y z\|_{L^p_0(A^l_r; \R^2)}.
\end{eqnarray*}
Therefore, there exists a constant $K_r$ dependent on $r$ such that
\begin{equation}
    \|z\|_{W^{2,p_0}(A^l_r;\R^2)} \leq K_r  \|\nabla^2 z\|_{L^p_0(A^l_r; \R^2)}
\end{equation}
for every $z \in X^l_r$. Hence, the claim is proved.
\paragraph{} Let $r$ satisfy (\ref{eqn: assume bound on r}) and $z \in W^{2, p_0}(A^l_r; \R^2)$. Let $w=T_1(z)$ be the solution of the boundary value problem 
\begin{equation}\label{eqn: BVP satisfied by T_1(z)}
    \begin{cases}
        -\text{div }(\C Ew)= f \circ \Phi + f^z \text{ in } A^l_r,\\
        (\C Ew )\nu^0=g \circ \Phi + \hat{g}^z+ \check{g}^z \text{ on } \Gamma^l_r,\\
        w=0 \text{ on } \pa A^l_r \setminus \Gamma^l_r.
    \end{cases}
\end{equation}
Now, we separate the boundary value problem above into two parts such that one depends on $z$ and the other does not, and obtain $T_1(z)$ as the sum of their respective solutions using linearity. Let $v_0$ be the solution to the boundary value problem
\begin{equation}
    \begin{cases}
        -\text{div }(\C Ew)= f \circ \Phi \text{ in } A^l_r,\\
        (\C Ew )\nu^0=g \circ \Phi  \text{ on } \Gamma^l_r,\\
        w=0 \text{ on } \pa A^l_r \setminus \Gamma^l_r,
    \end{cases}
\end{equation}
and $T_2(z)$ be the solution of the boundary value problem
\begin{equation}
    \begin{cases}
        -\text{div }(\C Ew)= f^z \text{ in } A^l_r,\\
        (\C Ew )\nu^0= \hat{g}^z+ \check{g}^z \text{ on } \Gamma^l_r,\\
        w=0 \text{ on } \pa A^l_r \setminus \Gamma^l_r.
    \end{cases}
\end{equation}
Clearly, $T_1(z)=v_0+T_2(z)$. Note that $T_1, T_2$ are well defined and $T_1(z), T_2(z), v_0 \in X^l_r$ for $z \in X^l_r$ by Theorem \ref{thm: 13, elliptic regularity in triangle}. Also, $T_1, T_2$ are linear operators on $X^l_r$ and $T_2(0)=0$.
\paragraph{} We claim that $T_1: X^l_r \rightarrow X^l_r$ is a contraction for sufficiently small $r>0$, that is, $\|T_1(z)-T_1(w)\|_{X^l_r} \leq \|z-w\|_{X^l_r}$. The claim follows directly if $T_2: X^l_r \rightarrow X^l_r$ is a contraction. By linearity, it is sufficient to show that 
\begin{equation}\label{eqn: sufficient condition for contraction map}
    \|T_2(z)\|_{X^l_r} \leq 1/2 \text{ for every } z \in X^l_r \text{ with } \|z\|_{X^l_r} \leq 1.
\end{equation}
By (\ref{eqn: elliptic regularity bound on second derivative}) from Theorem \ref{thm: 13, elliptic regularity in triangle}, we have 
\begin{equation}\label{eqn: sufficient condition for contraction map 2}
    \|T_2(z)\|_{X^l_r} \leq \kappa (\|f^z\|_{L^{p_0}(A^l_r; \R^2)}+ \|\nabla \hat{g}^z\|_{L^{p_0}(A^l_r; \R^2)} + \|\nabla \check{g}^z\|_{L^{p_0}(A^l_r; \R^2)})
\end{equation}
for $\kappa$ independent of $r,h,f,g$ and $z$. Hence, to prove (\ref{eqn: sufficient condition for contraction map}), we will prove that
\begin{equation}
    \kappa (\|f^z\|_{L^{p_0}(A^l_r; \R^2)}+ \|\nabla \hat{g}^z\|_{L^{p_0}(A^l_r; \R^2)} + \|\nabla \check{g}^z\|_{L^{p_0}(A^l_r; \R^2)}) \leq 1/2.
\end{equation}
\paragraph{} In the rest of the proof, let $C_h$ denote a constant that depends only on $h$ and is independent of $r,m$. Let $C_{\lambda, \mu}$ denote a constant that depends only on the Lam\'e coefficients $\lambda, \mu$.
\paragraph{} Let $r>0$ satisfy (\ref{eqn: assume bound on r}). Recall that $I_r := (0,r)$. By (\ref{eqn: sigma upper and lower l^infinity bound}), (\ref{eqn: sigma +1 L^infinity bound}), (\ref{eqn: def f^v_1}) and (\ref{eqn: def f^v_2}), 
\begin{eqnarray}
    |f^z(x,y)| &\leq& C_{\lambda, \mu}(\|\sigma-1\|_{L^\infty(I_r)}+ \sup_{(x,y)\in A^l_r}|y \sigma'(x)|+ \sup_{(x,y)\in A^l_r}|y \sigma'(x)|^2)|\nabla^2 z(x,y)|\nonumber\\
    & &+ C_{\lambda, \mu}(\|\sigma'\|_{L^\infty(I_r)}+ \sup_{(x,y)\in A^l_r}|y \sigma''(x)|)|\nabla z(x,y)|, \label{eqn: f^z pointwise estimate}
\end{eqnarray}
for a.e. $(x,y) \in A^l_r$. Hence, by (\ref{eqn: sigma-1 bound}), (\ref{eqn: ysigma' bound_omega}), (\ref{eqn: ysigma' bound_triangle}) and (\ref{eqn: ysigma'' bound_triangle}), we have 
\begin{equation}
    |f^z(x,y)| \leq C_h r |\nabla^2 z(x,y)|+ C_h |\nabla z(x,y)|,
\end{equation}
for a.e. $(x,y) \in A^l_r$. Since $\|z\|_{X^l_r}= \|\nabla ^2 z\|_{L^{p_0}(A^l_r; \R^2)} \leq 1$, using Poincar\'e's inequality from Lemma \ref{lemma: 15, Poincare's inequality} in the equation above, we have
\begin{eqnarray}\label{eqn: f^z C_h r bound}
    \|f^z\|_{L^{p_0}(A^l_r;\R^2)} &\leq & C_h r \|\nabla ^2 z\|_{L^{p_0}(A^l_r; \R^2)} + C_h \|\nabla  z\|_{L^{p_0}(A^l_r; \R^2)} \leq C_h r.
\end{eqnarray}
From (\ref{eqn: g hat g check def 1})-(\ref{eqn: g hat g check def end}), we can compute $\nabla \hat{g}$ and $\nabla \check{g}$ and obtain 
\begin{eqnarray}
    \|\nabla \hat{g}^z\|_{L^{p_0}(A^l_r;\R^2)} &\leq & C_{\lambda, \mu}\bigg( \sum_{i=1}^6 \|\omega_i\|_{L^\infty(I_r)}\|\nabla ^2 z\|_{L^{p_0}(A^l_r; \R^2)}\bigg)\nonumber\\
    & &+C_{\lambda, \mu}\bigg(\sum_{i=1}^6 \|\omega_i'\|_{L^\infty(I_r)}\|\nabla  z\|_{L^{p_0}(A^l_r; \R^2)} \bigg), \label{eqn: nabla g hat L^p_0 bound}\\
    \|\nabla \check{g}^z\|_{L^{p_0}(A^l_r;\R^2)} &\leq & C_{\lambda, \mu}\bigg( \sum_{i=1}^6 \|\omega_i\|_{L^\infty(I_r)}\|\nabla ^2 z\|_{L^{p_0}(A^l_r; \R^2)}\bigg)\nonumber\\
    & &+C_{\lambda, \mu}\bigg(\sum_{i=1}^6 \|\omega_i'\|_{L^\infty(I_r)}\|\nabla  z\|_{L^{p_0}(A^l_r; \R^2)} \bigg).\label{eqn: nabla g check L^p_0 bound}
\end{eqnarray}
Hence, by (\ref{eqn: omega_i bound}), (\ref{eqn: omega'_i bound}) from Lemma \ref{lem: 11} and Poincar\'e's inequality from Lemma \ref{lemma: 15, Poincare's inequality}, we have 
\begin{equation}\label{eqn: g hat g check C_h r bound}
     \|\nabla \hat{g}^z\|_{L^{p_0}(A^l_r;\R^2)}+ \|\nabla \check{g}^z\|_{L^{p_0}(A^l_r;\R^2)} \leq C_h r.
\end{equation}
Choose $C_h$ large enough such that it satisfies (\ref{eqn: f^z C_h r bound}) and (\ref{eqn: g hat g check C_h r bound}). Then, if $0<r< \min\bigg\{\frac{1}{4\kappa C_h}, \frac{\eta_0^2}{4M}\bigg\}$, by (\ref{eqn: sufficient condition for contraction map 2}), (\ref{eqn: f^z C_h r bound}) and (\ref{eqn: g hat g check C_h r bound}), $\|T_2(z)\|_{X^l_r} \leq 1/2$ and hence $T_2$ is a contraction. This implies that $T_1$ is a contraction as well.\\\\
\textbf{Step 4: Conclusion.} Apply the Banach fixed point theorem to $T_1$ to obtain $z _0 \in X^l_r$ such that $T_1(z_0)=z_0$. Hence, by (\ref{eqn: BVP satisfied by T_1(z)}), $z_0$ solves (\ref{eqn: BVP satisfied by v=utilde o Phi}). Then, by (\ref{eqn: C^2 regularity of Phi, Psi}), the function $u_0:= z_0 \circ \Phi \in W^{2,p_0}(\Omega_h^{0,r}; \R^2)$ and solves (\ref{eqn: elliptic PDE, regularity at corners}). Note that $W^{2,p_0}(\Omega_h^{0,r}; \R^2) \subset H^1(\Omega_h^{0,r}; \R^2)$ since $p_0 \in (4/3,2)$. Therefore, by uniqueness of solution of (\ref{eqn: elliptic PDE, regularity at corners}), $\util=u_0$ and hence $\util \in W^{2,p_0}(\Omega_h^{0,r}; \R^2)$. Since $\util =u-w_0$ in $\Omega_h^{0, r/2}$ and $w_0 \in W^{2,p_0}(\Omega_h^{0,r}; \R^2)$, we obtain that $u \in W^{2,p_0}(\Omega_h^{0,r/2}; \R^2)$.

\end{proof}
In the following theorem, we derive the exact relation between the endpoints $\alpha, \beta$ of the film and their counterparts at time $t=0$. This would be used later in obtaining the differential equation satisfied by the endpoints as functions of time.
\begin{thm}\label{thm: 16}
    Under the hypotheses of Theorem \ref{thm: 6}, we have 
    \begin{equation}\label{eqn: h''(alpha)=h''(beta)=0}
        h''(\alpha)=0, \quad h''(\beta)=0.
    \end{equation}
    In addition to this, if we assume that 
    \begin{equation} \label{eqn: alpha, alpha_0 are close, beta, beta_0 are close}
        |\alpha-\alpha_0| < \delta_0 \quad \text{and} \quad |\beta-\beta_0|< \delta_0,
    \end{equation}
    then, we have
    \begin{gather}
        \sigma_0 \frac{\alpha-\alpha_0}{\tau} = 
        \frac{\gamma}{J(\alpha)}-\gamma_0 + \nu_0 \frac{h'(\alpha)}{J(\alpha)^2}\bigg(\frac{h''}{J^3}\bigg)'(\alpha), \label{eqn: alpha ode precursor}\\
        \sigma_0 \frac{\beta-\beta_0}{\tau} = 
        -\frac{\gamma}{J(\beta)}+\gamma_0 -\nu_0 \frac{h'(\beta)}{J(\beta)^2}\bigg(\frac{h''}{J^3}\bigg)'(\beta),\label{eqn: beta ode precursor}
    \end{gather}
    where $\htil$ is defined in (\ref{eqn: def h tilde}) and $J$ is defined in (\ref{eqn: Def J, J_0}).
\end{thm}

\begin{remark}\label{remark: additional term in ODEs} Note that $\alpha, \beta$ are dependent on $\tau$. The assumption (\ref{eqn: alpha, alpha_0 are close, beta, beta_0 are close}) is reasonable since we prove later that $|\alpha-\alpha_0|$ and $|\beta-\beta_0|$ converge to $0$ as we converge the time step size $\tau$ to $0$.
\end{remark}
To prove the theorem, we use the following lemma \cite[Lemma 6.5]{dal2025motion}.
\begin{lem}\label{lem: 17}
    Under the hypotheses of Theorem \ref{thm: 6}, let $1\leq p< \infty$ and let $v \in W^{2,p}(\Omega_h;\R^2)$ such that $v(x,0)=(0,0)$ for a.e. $x \in (\alpha, \beta)$. Then, there exists $\hat{v} \in W^{2,p}(\R^2;\R^2)$ such that $\hat{v}=v$ in $\Omega_h$ and $\hat{v}(x,0)=(0,0)$ for a.e. $x \in \R$.
\end{lem}

\begin{proof}[Proof of Theorem \ref{thm: 16}]\textbf{Step 1:} We know that $v= u-w_0 \in W^{2,p_0}(\Omega_h;\R^2)$ by Theorem \ref{thm: 14} and that $v(x,0)=(0,0)$ for a.e. $x \in (\alpha, \beta)$, where $w_0$ is as defined in (\ref{eqn: def w_0}). Then, by Lemma \ref{lem: 17}, extend $v$ to $\hat{v} \in W^{2,p_0}(\R^2;\R^2)$ such that $\hat{v}=v$ in $\Omega_h$ and $\hat{v}(x,0)=(0,0)$ for a.e. $x \in \R$. Then the function $\hat{u}:= \hat{v}+ w_0 \in W^{2,p_0}_{\text{loc}}(\R^2;\R^2)$ and $\hat{u}(x,0)=(e_0x,0)$ for a.e. $x \in \R$.
\paragraph{} Now, let $\varphi, \alpha_\varep, h_\varep$ be as defined in (\ref{eqn: define varphi_EL section}), (\ref{eqn: alpha_varepsilon condition}), (\ref{eqn: define h_epsilon}) and let $u_\varep$ be the restriction of $\hat{u}$ to $\Omega_{h_\varep}$. Note that $u_\varep(x,0)=(e_0x,0)$ for a.e. $x \in (\alpha_\varep, \beta)$. We claim that 
\begin{equation}\label{eqn: derivative of elastic energy term}
    \frac{d}{d \varep} \e(\alpha_\varep, \beta, h_\varep, u_\varep)\bigg|_{\varep=0} = \int_\alpha^\beta W(Eu(x,h(x)))\varphi(x) dx.
\end{equation}
By (\ref{eqn: Elastic energy lim inf inequality_EL section}), there exists $\varep_1>0$ such that
\begin{eqnarray*}
     \frac{\e(\alpha_\varep, \beta, h_\varep, u_\varep)-\e(\alpha, \beta,h,u)}{\varep}&=& -\frac{1}{\varep} \int_{\alpha_\varep}^{\alpha+2 \delta_0}\bigg(\int_{h_\varep(x)}^{h(x)} W(E\hat{u}(x,y))dy \bigg)dx \nonumber\\
    & & -\frac{1}{\varep} \int_{\alpha}^{\alpha_\varep} \bigg(\int_0^{h(x)} W(Eu(x,y))dy \bigg)dx= :I_\varep+ II_\varep
\end{eqnarray*}
for all $-\varep_1< \varep< 0$. Note that $\hat{u}, u_\varep$ used here are different from those used in (\ref{eqn: Elastic energy lim inf inequality_EL section}).
\paragraph{}Since $4/3<p_0< 2$ and $\hat{u} \in W^{2,p_0}_{\text{loc}}(\R^2;\R^2)$, by the Sobolev-Gagliardo-Nirenberg embedding theorem, $\nabla \hat{u} \in L^{p_0^*}_{\text{loc}}(\R^2;\R^{2\times 2})$ where $4<p_0^*= \frac{2 p_0}{2-p_0}$ is the critical Sobolev exponent. Let $p=2+\delta$ where $\delta>0$ is chosen such that $4+2 \delta < p_0^\ast$. Then, $\nabla \hat{u} \in L^{4+2\delta}_{\text{loc}}(\R^2;\R^2)$. Let $q$ be the conjugate of $p$. Then, there exists $\delta'>0$ such that $\frac{1}{q}= \frac{1}{2}+\delta'$. Now, we use (\ref{eqn: W upper and lower bounds by|.|^2}), (\ref{eqn: alpha_epsilon-alpha bound}) and H\"older's inequality with $p,q$ to get
\begin{eqnarray*}
    \int_{\alpha}^{\alpha_\varep} \bigg(\int_0^{h(x)} W(Eu(x,y))dy \bigg)dx &\leq& C \int_{\alpha}^{\alpha_\varep} \bigg(\int_0^{h(x)} |\nabla u(x,y)|^2 dy \bigg) dx \\
    &\leq& C \bigg(\int_{\alpha}^{\alpha_\varep} \bigg(\int_0^{h(x)} |\nabla u(x,y)|^{4+2\delta} dy \bigg) dx\bigg)^{1/p} \bigg( \int_{\alpha}^{\alpha_\varep} h(x)dx\bigg)^{1/q}\\
    &\leq& C \bigg(\int_{\alpha}^{\alpha_\varep} \bigg(\int_0^{h(x)} |\nabla u(x,y)|^{4+2\delta} dy \bigg) dx\bigg)^{1/p} |\varep|^{1+2\delta'},
\end{eqnarray*}
where we used the fact that $|h(x)| \leq L_0|\alpha_\varep-\alpha|$ when $x \in (\alpha, \alpha_\varep)$. This shows that $II_\varep \rightarrow 0$ as $\varep \rightarrow 0^-$.
\paragraph{}Let $\zeta = E\hat{u} \in W^{1,p_0}_{\text{loc}}(\R^2;\R^{2\times 2})$. We claim that $W\circ\zeta \in W^{1,1}_{\text{loc}}(\R^{2})$. Since $p_0>4/3$, $p_0'=\frac{p_0}{p_0-1}< \frac{2 p_0}{2-p_0}=p_0^*$ and hence $\zeta \in L^{p_0'}(\R^2; \R^{2 \times 2})$. By (\ref{eqn: def W}), we have the pointwise estimate
\begin{equation*}
    |\nabla W\circ \zeta| \leq C |\zeta||\nabla \zeta|
\end{equation*}
for $C>0$ independent of $\zeta$. Therefore, for every compact set $K \subset \R^2$,
\begin{equation}\label{eqn: W(zeta) L^1 bound}
    \int_K |\nabla W\circ \zeta| dx dy \leq C \|\zeta\|_{L^{p_0'}(K)}\|\nabla \zeta\|_{L^{p_0}(K)}< \infty.
\end{equation}
By considering the representative of $W \circ \zeta$ that is locally absolutely continuous on the line parallel to the axes, we obtain 
\begin{eqnarray*}
    W(\zeta(x,y))= W(\zeta(x, h(x)))-\int_y^{h(x)} \pa_s (W \circ \zeta)(x,s) ds
\end{eqnarray*}
for a.e. $x \in (\alpha_\varep, \alpha+ 2 \delta_0)$. Therefore, $I_\varep$ decomposes as $A_\varep+ B_\varep$, where
\begin{eqnarray*}
    A_\varep&:=& -\frac{1}{\varep} \int_{\alpha_\varep}^{\alpha+2 \delta_0}\bigg(\int_{h_\varep(x)}^{h(x)} W(\zeta(x,h(x))dy \bigg)dx, \text{ and}\\
    B_\varep&:=& -\frac{1}{\varep} \int_{\alpha_\varep}^{\alpha+2 \delta_0}\bigg(\int_{h_\varep(x)}^{h(x)} \bigg(\int_y^{h(x)}\pa_s (W \circ \zeta)(x,s) ds\bigg) dy \bigg)dx.
\end{eqnarray*}
By (\ref{eqn: alpha_epsilon-alpha bound}), (\ref{eqn: define h_epsilon}) and (\ref{eqn: omega_epsilon estimate}),
\begin{equation*}
    A_\varep = \int_{\alpha_\varep}^{\alpha+2 \delta_0}  W(\zeta(x,h(x)) \bigg(\varphi(x)+ \frac{\omega_\varep}{\varep}\psi(x)\bigg) dx  \rightarrow \int_{\alpha}^{\alpha+2 \delta_0}  W(\zeta(x,h(x)) \varphi(x) dx 
\end{equation*} as $\varep \rightarrow 0^-$. Further, by Fubini's theorem, (\ref{eqn: define h_epsilon}) and (\ref{eqn: W(zeta) L^1 bound}),
\begin{eqnarray*}
    |B_\varep| &=& \bigg| \frac{1}{\varep} \int_{\alpha_\varep}^{\alpha+2 \delta_0}\int_{h_\varep(x)}^{h(x)}(s-h_\varep(x))\pa_s (W \circ \zeta)(x,s) ds dx \bigg| \\
    &\leq& C \bigg\| \varphi+ \frac{\omega_\varep}{\varep}\psi\bigg\|_{L^\infty(\R)}\|\zeta\|_{L^{p_0'}(K)}\bigg( \int_{\alpha_\varep}^{\alpha+2 \delta_0}\int_{h_\varep(x)}^{h(x)} |\pa_y\zeta (x,s)|^{p_0}\bigg)^{1/p_0} \rightarrow 0
\end{eqnarray*}
as $\varep \rightarrow 0^-$, where $K = [\alpha, \beta] \times [0, L_0(\beta-\alpha)]$. Therefore, 
\begin{equation}\label{eqn: left derivative of elastic energy}
    \lim_{\varep \rightarrow 0^-} \frac{\e(\alpha_\varep, \beta, h_\varep, u_\varep)-\e(\alpha, \beta,h,u)}{\varep} = \int_{\alpha}^{\beta}  W(\zeta(x,h(x)) \varphi(x) dx 
\end{equation}
as supp$(\varphi) \in (\alpha- \delta_0, \alpha+ \delta_0)$. When $0< \varep< \varep_1$, we can follow similar steps and obtain the same limit as above as $\varep \rightarrow 0^+$.\\

\textbf{Step 2:} By Theorem \ref{thm: 6}, we get continuity of $h''$ up to the endpoints of $[\alpha, \beta]$. Hence, by (\ref{eqn: derivative of alpha_epsilon wrt epsilon}),
\begin{equation}\label{eqn: I_epsilon,3,2 limit}
    \lim\limits_{\varep \rightarrow 0}I_{\varep,3,2}=  \frac{\nu_0}{2} \frac{(h''(\alpha))^2}{(1+(h'(\alpha))^2)^{5/2}}\frac{1}{h'(\alpha)}
\end{equation}
where $I_{\varep,3,2}$ is defined in (\ref{eqn: splitting of I_(epsilon,3)}). Hence, by (\ref{eqn: surface energy difference quotient break down}), (\ref{eqn: I_epsilon, 2 limit_s energy}), (\ref{eqn: I_epsilon, 1 limit_s energy}), (\ref{eqn: I_epsilon,3,1 limit_s energy}) and (\ref{eqn: I_epsilon,3,2 limit}),
\begin{eqnarray}
     \lim_{ \varep \rightarrow 0}\frac{S(\alpha_\varep, \beta, h_\varep)-\s(\alpha, \beta,h)}{\varep}=  \gamma \int_\alpha^\beta \frac{h'\varphi'}{\sqrt{1+(h')^2}}dx + \gamma \sqrt{1+h'(\alpha))^2}\frac{1}{h'(\alpha)}-\gamma_0 \frac{1}{h'(\alpha)}\nonumber\\
     + \nu_0 \int_\alpha^\beta \frac{h''\varphi''}{(1+(h')^2)^{5/2}}dx-\frac{5}{2} \nu_0 \int_\alpha^\beta \frac{h'(h'')^2\varphi'}{(1+(h')^2)^{7/2}}dx
     +\frac{\nu_0}{2} \frac{(h''(\alpha))^2}{(1+(h'(\alpha))^2)^{5/2}}\frac{1}{h'(\alpha)}. \label{eqn: derivative of surface energy term}
 \end{eqnarray}
Then, by (\ref{eqn: derivative of incremental energy term}), (\ref{eqn: derivative of elastic energy term}) and (\ref{eqn: derivative of surface energy term}), $\f^0(\alpha_\varepsilon, \beta, h_\varepsilon, u_\varepsilon)$ is differentiable with respect to $\varep$ and since $(\alpha, \beta, h, u) \in \a$ minimizes $\f_0$,
\begin{equation}\label{eqn: derivative of total energy}
    \frac{d}{d \varepsilon} \f^0(\alpha_\varepsilon, \beta, h_\varepsilon, u_\varepsilon)\bigg|_{\varepsilon=0}=0.
\end{equation}
  Hence, by (\ref{eqn: left derivative of elastic energy}), (\ref{eqn: derivative of surface energy term}), and (\ref{eqn: derivative of total energy}),
  \begin{eqnarray*}
     \gamma \int_\alpha^\beta\frac{h'\varphi'}{J} dx + \gamma \frac{J(\alpha)}{h'(\alpha)} -\gamma_0 \frac{1}{h'(\alpha)}
     +\nu_0 \int_\alpha^\beta \frac{h''\varphi''}{J^5}dx
     -\frac{5}{2}\nu_0 \int_\alpha^\beta \frac{h'(h'')^2 \varphi'}{J^7}dx\\+\frac{\nu_0}{2}\frac{(h''(\alpha))^2}{(J(\alpha)^5)}\frac{1}{h'(\alpha)} 
     +\int_\alpha^\beta \overline{W}\varphi dx 
     +\frac{1}{\tau}\int_{\alpha_0}^{\beta_0} \frac{(\htil(x)-h_0(x))\varphi(x)}{J_0}-\sigma_0 \frac{\alpha-\alpha_0}{\tau}\frac{1}{h'(\alpha)} =0.
  \end{eqnarray*}
Integrating by parts,
\begin{eqnarray*}
 -  \gamma \int_\alpha^\beta\bigg( \frac{h'}{J}\bigg)'\varphi dx-\gamma \frac{h'(\alpha)}{J(\alpha)}\varphi(\alpha)+\gamma \frac{h'(\beta)}{J(\beta)}\varphi(\beta)
 + \gamma \frac{J(\alpha)}{h'(\alpha)} -\gamma_0 \frac{1}{h'(\alpha)}\\+ \nu_0 \bigg(\int_\alpha^\beta\bigg( \frac{h''}{J^5}\bigg)''\varphi dx
 -\frac{h''(\alpha)}{J^5(\alpha)} \varphi'(\alpha)+\frac{h''(\beta)}{J^5(\beta)}\varphi'(\beta) +  \bigg(\frac{h''}{J^5}\bigg)' (\alpha) \varphi(\alpha)
 - \bigg(\frac{h''}{J^5}\bigg)' (\beta) \varphi(\beta)\bigg)\\
  +\frac{5\nu_0}{2} \bigg(\int_\alpha^\beta \bigg(\frac{h'(h'')^2 }{J^7}\bigg)'\varphi dx
  +\frac{h'(\alpha)(h''(\alpha))^2}{J^7(\alpha)}\varphi(\alpha)-\frac{h'(\beta)(h''(\beta))^2}{J^7(\beta)}\varphi(\beta)\bigg)\\+\frac{\nu_0}{2}\frac{(h''(\alpha))^2}{(J(\alpha)^5)}\frac{1}{h'(\alpha)} 
     +\int_\alpha^\beta \overline{W}\varphi dx +\frac{1}{\tau}\int_{\alpha_0}^{\beta_0}\frac{(\htil-h_0)\varphi}{J_0}dx -\sigma_0 \frac{\alpha-\alpha_0}{\tau}\frac{1}{h'(\alpha)} =0.
\end{eqnarray*}
As $\varphi(\alpha)=1$, $\varphi(\beta)=0$ and $\varphi'(\beta)=0$,
\begin{gather}
    \nu_0 \int_\alpha^\beta\bigg( \frac{h''}{J^5}\bigg)''\varphi dx =  \gamma \int_\alpha^\beta\bigg( \frac{h'}{J}\bigg)'\varphi dx -\frac{5}{2}\nu_0 \int_\alpha^\beta \bigg(\frac{h'(h'')^2 }{J^7}\bigg)'\varphi dx\nonumber
    -\int_\alpha^\beta \overline{W}\varphi dx \\-\frac{1}{\tau}\int_{\alpha_0}^{\beta_0}\frac{(\htil-h_0)\varphi}{J_0}dx\nonumber
      -\gamma \frac{J(\alpha)}{h'(\alpha)}+\gamma \frac{h'(\alpha)}{J(\alpha)} +\gamma_0 \frac{1}{h'(\alpha)}+\nu_0 \frac{h''(\alpha)}{J^5(\alpha)} \varphi'(\alpha)-\nu_0 \bigg(\frac{h''}{J^5}\bigg)' (\alpha) \nonumber\\
     -\frac{5 \nu_0}{2}\frac{h'(\alpha)(h''(\alpha))^2}{J^7(\alpha)} -\frac{\nu_0}{2}\frac{(h''(\alpha))^2}{(J(\alpha)^5)}\frac{1}{h'(\alpha)} +\sigma_0 \frac{\alpha-\alpha_0}{\tau}\frac{1}{h'(\alpha)}. \label{eqn: before dividing by varphi'(alpha)}
\end{gather}
Since $\|\varphi\|_{L^{\infty}(\R)} \leq 2$, the left hand side of (\ref{eqn: before dividing by varphi'(alpha)}) is bounded by a number independent of $\varphi$. Hence, if we divide (\ref{eqn: before dividing by varphi'(alpha)}) by $\varphi'(\alpha)$ and take $\varphi'(\alpha) \rightarrow \infty$, we have 
\begin{equation*}
    \nu_0 \frac{h''(\alpha)}{J^5(\alpha)} =0
\end{equation*}
Hence, we obtain
\begin{equation}\label{eqn:h'' alpha =0}
    h''(\alpha)=0.
\end{equation}
Using (\ref{eqn: EL strong equation}), (\ref{eqn: before dividing by varphi'(alpha)}), (\ref{eqn:h'' alpha =0}) and the fact that $m\int_\alpha^\beta \varphi dx =0$, we have
\begin{eqnarray}
    0&=&\frac{1}{\tau}\bigg(\int_{\alpha}^\beta \frac{\htil-h_0}{J_0}\chi_{[\alpha_0, \beta_0]}\varphi dx-\int_{\alpha_0}^{\beta_0} \frac{\htil-h_0}{J_0}\varphi dx \bigg) -\gamma \frac{J(\alpha)}{h'(\alpha)}+\gamma \frac{h'(\alpha)}{J(\alpha)} +\gamma_0 \frac{1}{h'(\alpha)}\nonumber \\
    & & -\nu_0 \bigg(\frac{h''}{J^5}\bigg)' (\alpha) +\sigma_0 \frac{\alpha-\alpha_0}{\tau}\frac{1}{h'(\alpha)}. \label{eqn: ODE precursor with integral term}
\end{eqnarray}
Recall the definition of $\varphi$ from (\ref{eqn: define varphi_EL section}). Let $0<\delta< \delta_0$ and choose $\varphi_0^\delta \in C^\infty(\R)$ such that $\varphi_0^\delta(0)=1$, $\varphi_0^\delta(x) \geq 1/2$ for every $x \in [-\frac{\delta}{2}, \frac{\delta}{2} ]$, $\int_0^{\delta} \varphi_0 dx =0$ and supp$(\varphi_0) \subset (-\delta, \delta)$. Define $\varphi^\delta:= \varphi_0^\delta(x-\alpha)$ for $x \in \R$. Note that the analysis above would not be affected if $\varphi^\delta$ was used instead of $\varphi$. By (\ref{eqn: beta-alpha lower bound}) and (\ref{eqn: alpha, alpha_0 are close, beta, beta_0 are close}), we have  $\beta, \beta_0 \notin \text{supp}(\varphi^\delta)$ since
\begin{eqnarray*}
    |\alpha-\beta_0|\geq |\alpha-\beta|-|\beta-\beta_0|\geq 63 \delta_0> \delta.
\end{eqnarray*}
 Therefore, 
\begin{equation*}
    \frac{h'(\alpha)}{\tau}\bigg(\int_{\alpha_0}^{\beta_0} \frac{\htil-h_0}{J_0}\varphi^\delta dx -\int_{\alpha}^\beta \frac{\htil-h_0}{J_0}\chi_{[\alpha_0, \beta_0]}\varphi^\delta dx\bigg) = \begin{cases}
        \frac{h'(\alpha)}{\tau} \int_{\alpha_0}^{\alpha}\frac{\htil-h_0}{J_0}\varphi^\delta dx \text{ if } \alpha_0< \alpha,\\
        0, \;\;\;\;\;\;\;\;\;\;\;\;\;\;\;\;\;\;\;\;\; \quad \text{if } \alpha_0 \geq \alpha.
    \end{cases}
\end{equation*}
By letting $\delta \rightarrow 0^+$, it follows from the Lebesgue's dominated convergence theorem that $\int_{\alpha_0}^{\alpha}\frac{\htil-h_0}{J_0}\varphi^\delta dx \rightarrow 0$.
Now, multiplying (\ref{eqn: ODE precursor with integral term}) by $h'(\alpha)$, we have
\begin{eqnarray*}
    \sigma_0 \frac{\alpha-\alpha_0}{\tau} &=& \gamma \bigg(J(\alpha)-\frac{(h'(\alpha))^2}{J(\alpha)}\bigg)-\gamma_0 +\nu_0 h'(\alpha) \bigg(\frac{h''}{J^5}\bigg)' (\alpha)\\
    %&=& \frac{h'(\alpha)}{\tau}\bigg(\int_{\alpha_0}^{\beta_0} \frac{\htil-h}{J_0}\varphi dx -\int_{\alpha}^\beta \frac{\htil-h}{J_0}\chi_{[\alpha_0, \beta_0]}\varphi dx\bigg)\\
    %& & +\gamma \bigg(\frac{(J(\alpha))^2-(h'(\alpha))^2}{J(\alpha)}\bigg)-\gamma_0 +\nu_0 h'(\alpha) \bigg(\frac{h''}{J^5}\bigg)' (\alpha)\\
    &=& \frac{\gamma}{J(\alpha)}- \gamma_0 +\nu_0 \frac{h'(\alpha)}{J(\alpha)^2} \bigg(\frac{h''}{J^3}\bigg)' (\alpha),
\end{eqnarray*}
using (\ref{eqn:h'' alpha =0}). Similarly, we can prove that $h''(\beta)=0$ and that 
\begin{eqnarray*}
    \sigma_0 \frac{\beta-\beta_0}{\tau} %= \frac{h'(\beta)}{\tau}\bigg(\int_{\alpha_0}^{\beta_0} \frac{\htil-h}{J_0}\varphi_1 dx -\int_{\alpha}^\beta \frac{\htil-h}{J_0}\chi_{[\alpha_0, \beta_0]}\varphi_1 dx\bigg)
     %-\frac{\gamma}{J(\beta)}\\+ \gamma_0-\nu_0 h'(\beta) \bigg(\frac{h''}{J^5}\bigg)' (\beta)
    =-\frac{\gamma}{J(\beta)}+ \gamma_0 -\nu_0 \frac{h'(\beta)}{J(\beta)^2} \bigg(\frac{h''}{J^3}\bigg)' (\beta).
\end{eqnarray*}
   
\end{proof}
We proved in Theorem \ref{thm: 14} that $u \in W^{2,p_0}(\Omega_h)$. Now, we obtain an estimate of the $L^{p_0}$ norm of $\nabla^2 u$ in the corner of the domain $\Omega_h$ near $\alpha$. We also obtain an estimate of the $L^p$ norm of $\nabla u$ on the graph of $h$ near $\alpha$ in the sense of traces, where $p=\frac{p_0}{2-p_0}$. Both estimates are calculated in terms of the $L^{p_0}$ norm of $\nabla u$ in $\Omega_h^{0,r}$. Recall that $\alpha$ is taken to be $0$ without loss of generality, $\Omega_h^{0,r}$ is defined in (\ref{eqn: Part of Omega_h}) and that $\Gamma_h^{0,r}=\{(x,h(x): 0<x<r\}$.

\begin{thm}\label{thm: 18}
    Under the assumptions of Theorem \ref{thm: 6}, let $0<r< \delta_0$. Then, there exist two constants $0<c_1=c_1(\eta_0, \eta_1,L_0, M)<1$ and $c_2=c_2(\eta_0, \eta_1,L_0, M)>0$ independent of $r$, such that if 
    \begin{equation}\label{eqn: r h''' L^1 bound by constant}
        r\|h'''\|_{L^1(I_r)} \leq c_1,
    \end{equation}
    then 
    \begin{gather}
        \|\nabla^2 u\|_{L^{p_0}(\Omega_h^{0,r/2})} \leq c_2 + \frac{c_2}{r}\|\nabla u\|_{L^{p_0}(\Omega_h^{0,r})} , \label{eqn: bound on nabla^2 u with nabla u}\\
        \|\nabla u\|_{L^{p_0/(2-p_0)}(\Gamma_h^{0,r/2})}  \leq c_2 + \frac{c_2}{r}\|\nabla u\|_{L^{p_0}(\Omega_h^{0,r})}. \label{eqn: bound on nabla u trace with nabla u}
    \end{gather}
\end{thm}
We state the following lemma \cite[Lemma 6.7]{dal2025motion} to be used in the proof of the theorem.
\begin{lem}\label{lem: 19}
    Let $0< \eta_0\leq l\leq L_0$, $1<p<\infty$ and $\R^2_+:= \R \times (0, +\infty)$. Define
    \begin{equation}\label{eqn: def A^l_infty}
        A^l_\infty := \{(x,y) \in \R^2: x>0 \text{ and } 0<y< lx\},
    \end{equation}
    and let $v \in W^{2,p}(A^l_\infty)$ such that $v(x,0)=0$ for $x>0$. Then, the function 
    \begin{equation}
        \hat{v}(x,y):= \begin{cases}
            3v(2y/m-x,y)-2v(3y/m-2x,y) \text{ if } (x,y) \in \R^2_+ \setminus A^l_\infty,\\
            v(x,y) \hspace{5.45 cm} \text{if } (x,y) \in A^l_\infty
        \end{cases}
    \end{equation}
    belongs to $W^{2,p}(\R^2_+)$ with $\hat{v}(x,0)=0$ for $x >0$. Moreover, 
    \begin{equation}
        \|\nabla^2 \hat{v} \|_{L^p(\R^2_+)} \leq \|\nabla v\|_{L^p(A^l_\infty)},
    \end{equation}
    for some constant $C>0$ depending on $\eta_0, L_0$ and $p$.
\end{lem}

\begin{proof}[Proof of Theorem \ref{thm: 18}.] As in the proof of Theorem \ref{thm: 14}, we transform the corner $\Omega_h^{0,r}$ to the triangle $A^l_r$ and then use the estimates in Lemma \ref{lem: 8} and the results from Lemma \ref{lem: 19} to prove the theorem.
\paragraph{}Let $v$ be as defined in (\ref{eqn: define v using utilde}). Since $\util \in W^{2,p_0}(\Omega_h^{0,r}; \R^2)$ by Theorem \ref{thm: 14}, $v \in W^{2,p_0}(A^l_r; \R^2)$ by (\ref{eqn: C^2 regularity of Phi, Psi}). By (\ref{eqn: elliptic regularity bound on second derivative}) and (\ref{eqn: BVP satisfied by v=utilde o Phi}), 
\begin{eqnarray*}
     \|\nabla^2 v \|_{L^{p_0}(A^l_r)} &\leq&  \kappa ( \|f \circ \Phi\|_{L^{p_0}(A^l_r; \R^2)}+ \|f^v\|_{L^{p_0}(A^l_r; \R^2)}\\
     & &+ \|\nabla g \circ \Phi\|_{L^{p_0}(A^l_r; \R^2)}+\|\nabla \hat{g}^v\|_{L^{p_0}(A^l_r; \R^2)} + \|\nabla \check{g}^v\|_{L^{p_0}(A^l_r; \R^2)}).
\end{eqnarray*}
  Note that by (\ref{eqn: def u-tilde}) and (\ref{eqn: define v using utilde}), $v(x,y)=0$ for every $(x,y) \in A^l_r$ such that $x \geq 7r/8$. Hence, we can extend $v$ to $A^l_\infty$ (without relabeling) by setting $v=0$ in $A^l_\infty \setminus A^l_r$ and have 
  \begin{eqnarray}
      \|\nabla^2 v \|_{L^{p_0}(A^l_\infty)} &\leq&  \kappa ( \|f \circ \Phi\|_{L^{p_0}(A^l_r; \R^2)}+ \|f^v\|_{L^{p_0}(A^l_r; \R^2)}\nonumber \\
     & &+ \|\nabla g \circ \Phi\|_{L^{p_0}(A^l_r; \R^2)}+\|\nabla \hat{g}^v\|_{L^{p_0}(A^l_r; \R^2)} + \|\nabla \check{g}^v\|_{L^{p_0}(A^l_r; \R^2)}). \label{eqn: L^p bound of nabla^2 v in infinite triangle}
  \end{eqnarray}
The constant $\kappa $ is independent of $r,m,f$ and $g$. Now, we will estimate each term in the right hand side of the equation above in terms of $r$ and $h$. 
\paragraph{} By Theorem \ref{thm: 16}, $h''(\alpha)=0$. Hence, by the fundamental theorem of calculus,
\begin{equation}\label{eqn: h'' L^infinity bound using h''' L^1 bound}
    \|h''\|_{L^\infty(I_r)} \leq \|h'''\|_{L^1(I_r)}.
\end{equation}
Now, by using (\ref{eqn: sigma-1 bound}), (\ref{eqn: sigma' bound}), (\ref{eqn: ysigma' bound_triangle}), (\ref{eqn: ysigma'' bound_triangle}), (\ref{eqn: f^z pointwise estimate}) and (\ref{eqn: h'' L^infinity bound using h''' L^1 bound}), we have
\begin{eqnarray*}
    \|f^v\|_{L^{p_0}(A^l_r; \R^2)} &\leq& C(r\|h'''\|_{L^1(I_r)}+r^2\|h'''\|_{L^1(I_r)}^2)\|\nabla^2 v \|_{L^{p_0}(A^l_r)}\\
    & & + C(\|h'''\|_{L^1(I_r)}+r\|h'''\|_{L^1(I_r)}^2)\|\nabla v \|_{L^{p_0}(A^l_r)}.
\end{eqnarray*}
We can use Poincar\'e's inequality in Lemma \ref{lemma: 15, Poincare's inequality} on $\nabla v$ since $\pa_xv=\pa_y v =0 $ 
on the side $\{(r,y): 0\leq y\leq lr\}$ of the triangle $A^l_r$. Then, $\|\nabla v \|_{L^{p_0}(A^l_r)} \leq Cr\|\nabla^2 v \|_{L^{p_0}(A^l_r)}$ for $C=C(L_0, \eta_0)$ and 
\begin{equation}
    \|f^v\|_{L^{p_0}(A^l_r; \R^2)} \leq C(r\|h'''\|_{L^1(I_r)}+r^2\|h'''\|_{L^1(I_r)}^2)\|\nabla^2 v \|_{L^{p_0}(A^l_r)}. \label{eqn: f^v L^p_0 bound for v bound}
\end{equation}
Using the Poincar\'e's inequality in Lemma \ref{lemma: 15, Poincare's inequality}, (\ref{eqn: omega_i bound}), (\ref{eqn: omega'_i bound}) and (\ref{eqn: h'' L^infinity bound using h''' L^1 bound}) in (\ref{eqn: nabla g hat L^p_0 bound}), 
\begin{eqnarray}
    \|\nabla \hat{g}^v\|_{L^{p_0}(A^l_r;\R^2)} &\leq& Cr\|h'''\|_{L^1(I_r)}\|\nabla^2 v \|_{L^{p_0}(A^l_r)} \nonumber\\
    & & + C(\|h'''\|_{L^1(I_r)}+r\|h'''\|_{L^1(I_r)}^2)\|\nabla v \|_{L^{p_0}(A^l_r)}\nonumber\\
    &\leq& C (r\|h'''\|_{L^1(I_r)}+r^2\|h'''\|_{L^1(I_r)}^2)\|\nabla^2 v \|_{L^{p_0}(A^l_r)}.\label{eqn: nabla g hat L^p_0 bound for v bound}
\end{eqnarray}
Similarly, from (\ref{eqn: nabla g check L^p_0 bound}), we have 
\begin{equation}
    \|\nabla \check{g}^v\|_{L^{p_0}(A^l_r;\R^2)} \leq C (r\|h'''\|_{L^1(I_r)}+r^2\|h'''\|_{L^1(I_r)}^2)\|\nabla^2 v \|_{L^{p_0}(A^l_r)}. \label{eqn: nabla g check L^p_0 bound for v bound}
\end{equation}
Recall that in the definition of $f$ and $g$ ((\ref{eqn: def f_u regularity proof}) and (\ref{eqn: def g_u regularity proof})), $u-w_0$ satisfies the conditions of Lemma \ref{lemma: 15, Poincare's inequality} and $\varphi_r$ has the properties $\|\nabla \varphi_r\|_{L^\infty(\R^2)}\leq C/r$ and $\|\nabla^2 \varphi_r\|_{L^\infty(\R^2)}\leq C/r^2$. Therefore, by (\ref{eqn: f o Phi L^p_0 bound independent of h}) from Lemma \ref{lem: 9},
\begin{gather}
    \|f \circ \Phi\|_{L^{p_0}(A^l_r; \R^2)} \leq  C \|f\|_{L^{p_0}(\Omega_h^{0,r}; \R^2)}
    \leq C\|\nabla^2 \varphi_r\|_{L^\infty(\R^2)}\|u-w_0\|_{L^{p_0}(\Omega_h^{0,r})}\nonumber\\
     +C\|\nabla \varphi_r\|_{L^\infty(\R^2)}\|\nabla u-\nabla w_0\|_{L^{p_0}(\Omega_h^{0,r})} \leq \frac{C}{r}\|\nabla u-\nabla w_0\|_{L^{p_0}(\Omega_h^{0,r})}. \label{eqn: f o Phi L^p_0 bound for v bound}
\end{gather}
 Further, by (\ref{eqn: nabla w o Phi L^p_0 bound independent of h}) from Lemma \ref{lem: 9}, (\ref{eqn: ysigma' bound_triangle}) from Lemma \ref{lem: 8} and (\ref{eqn: h'' L^infinity bound using h''' L^1 bound}), since $h''(\alpha)=0$ by Theorem \ref{thm: 16},
\begin{eqnarray}
    \|\nabla (g \circ \Phi)\|_{L^p(A^l_r)} &\leq& C(1+ \sup_{(x,y)\in A^l_r}|y \sigma'(x)|) \|\nabla g\|_{L^p(\Omega_h^{0,r})}\nonumber\\
    &\leq& C(1+ \sup_{(x,y)\in A^l_r}|y \sigma'(x)|) (\|\nabla^2 \varphi_r\|_{L^\infty(\R^2)}\|u-w_0\|_{L^{p_0}(\Omega_h^{0,r})}\nonumber\\
      &&+\|\nabla \varphi_r\|_{L^\infty(\R^2)}\|\nabla u-\nabla w_0\|_{L^{p_0}(\Omega_h^{0,r})}+\|\nabla \varphi_r\|_{L^\infty(\R^2)}\| u-w_0\|_{L^{p_0}(\Omega_h^{0,r})}\|h''\|_{L^\infty(I_r)})\nonumber\\
    &\leq& C(1+r\|h'''\|_{L^1(I_r)})\bigg[\frac{1}{r^2}\|u-w_0\|_{L^{p_0}(\Omega_h^{0,r})} \nonumber\\
    &&+\frac{1}{r}\|\nabla u-\nabla w_0\|_{L^{p_0}(\Omega_h^{0,r})}+ 
    \frac{1}{r}\|u-w_0\|_{L^{p_0}(\Omega_h^{0,r})}\|h'''\|_{L^1(I_r)}\bigg] \nonumber\\
    &\leq& \frac{C}{r}(1+r\|h'''\|_{L^1(I_r)}+r^2\|h'''\|_{L^1(I_r)}^2 )\|\nabla u-\nabla w_0\|_{L^{p_0}(\Omega_h^{0,r})} .\label{eqn: nabla g o Phi L^p_0 bound for v bound}
\end{eqnarray}
Therefore, applying (\ref{eqn: f^v L^p_0 bound for v bound})-(\ref{eqn: nabla g o Phi L^p_0 bound for v bound}) in (\ref{eqn: L^p bound of nabla^2 v in infinite triangle}), 
\begin{eqnarray}
     \|\nabla^2 v \|_{L^{p_0}(A^l_r)} &\leq& \frac{C}{r}(1+r\|h'''\|_{L^1(I_r)}+r^2\|h'''\|_{L^1(I_r)}^2 )\|\nabla u-\nabla w_0\|_{L^{p_0}(\Omega_h^{0,r})}\nonumber\\
     & & +C(r\|h'''\|_{L^1(I_r)}+r^2\|h'''\|_{L^1(I_r)}^2)\|\nabla^2 v \|_{L^{p_0}(A^l_r)}. \label{eqn: nabla^2 v L^p upper bound with u,w_0,h}
\end{eqnarray}
Fix $c_1=c_1(\eta_0, \eta_1, M)>0$ such that $c_1<1$ and $Cc_1<1/4$ where $C$ is the constant in (\ref{eqn: nabla^2 v L^p upper bound with u,w_0,h}). If $r$ satisfies (\ref{eqn: r h''' L^1 bound by constant}) for this $c_1$, we have 
\begin{eqnarray*}
    C(r\|h'''\|_{L^1(I_r)}+r^2\|h'''\|_{L^1(I_r)}^2)\|\nabla^2 v \|_{L^{p_0}(A^l_r)} \leq C(c_1+c_1^2)\|\nabla^2 v \|_{L^{p_0}(A^l_r)} \leq \frac{1}{2}\|\nabla^2 v \|_{L^{p_0}(A^l_r)}
\end{eqnarray*}
and hence, from (\ref{eqn: nabla^2 v L^p upper bound with u,w_0,h}), we have 
\begin{eqnarray*}
    \|\nabla^2 v \|_{L^{p_0}(A^l_r)} &\leq& \frac{C}{r}\|\nabla u-\nabla w_0\|_{L^{p_0}(\Omega_h^{0,r})}\\
    &\leq&  \frac{C}{r}\|\nabla u\|_{L^{p_0}(\Omega_h^{0,r})}+  \frac{C}{r} r^{2/p_0}
\end{eqnarray*}
as $\|\nabla w_0\|_{L^{p_0}(\Omega_h^{0,r})} \leq  C\l^N(\Omega_h^{0,r})^{1/p_0}\leq Cr^{2/p_0}$. Since $p_0<2$ and $r<1$, this in turn gives us 
\begin{equation} \label{eqn: nabla^2 v L^p_0 bound}
    \|\nabla^2 v \|_{L^{p_0}(A^l_r)} \leq  \frac{C}{r}\|\nabla u\|_{L^{p_0}(\Omega_h^{0,r})}+ C.
\end{equation}
Using Lemma \ref{lem: 9}, Lemma \ref{lemma: 15, Poincare's inequality}, (\ref{eqn: sigma' bound}), (\ref{eqn: ysigma' bound_omega}), (\ref{eqn: ysigma'' bound_omega}) and the fact that $r<1$, we have 
\begin{eqnarray*}
     \|\nabla^2 \util \|_{L^{p_0}(\Omega_h^{0,r})} &\leq& C(1+r\|h'''\|_{L^1(I_r)}+r^2\|h'''\|_{L^1(I_r)}^2)\|\nabla^2 v \|_{L^{p_0}(A^l_r)}\\
     &\leq& C |\nabla^2 v \|_{L^{p_0}(A^l_r)} \leq \frac{C}{r}\|\nabla u\|_{L^{p_0}(\Omega_h^{0,r})}+  \frac{C}{r} r^{2/p_0},
\end{eqnarray*}
where we used the assumption (\ref{eqn: r h''' L^1 bound by constant}). By the choice of $\varphi_r$ in (\ref{eqn: def u-tilde}), $\util = u-w_0$ in $\Omega_h^{0, r/2}$. Therefore, since $\nabla^2 w_0=0$, we obtain (\ref{eqn: bound on nabla^2 u with nabla u}).\\\\
\textbf{Step 2:} In this step, we prove (\ref{eqn: bound on nabla u trace with nabla u}) by extending the function $v$ to $\R^2$ and then using standard trace theory. Note that $v \in W^{2, p_0}(A^l_r ; \R^2)$ is zero outside $A^l_{7r/8}$. Hence, we can extend $v$ by zero to obtain $v \in  W^{2, p_0}(A^l_\infty ; \R^2)$ (not relabeled) where $A^l_\infty$ is defined in (\ref{eqn: def A^l_infty}).
\paragraph{} By Lemma \ref{lem: 19}, the function $\hat{v}$ defined by
\begin{equation}\label{eqn: def v hat}
    \hat{v}(x,y):= \begin{cases}
        3v(2y/m-x,y)-2v(3y/m-2x,y) \text{ if } (x,y) \in \R^2_+ \setminus A^l_\infty,\\
        v(x,y) \hspace{5.3 cm} \text{ if } (x,y) \in A^l_\infty
    \end{cases}
\end{equation}
belongs to $W^{2, p_0}(\R^2_+ ; \R^2)$ with $\hat{v}(x,0)=(0,0)$ for every $x \in \R$. Now, define $\overline{v}$ that extends $\hat{v}$ to $\R^2$ as 
\begin{equation}\label{eqn: def v bar}
\overline{v}(x,y):= \begin{cases}
        3 \hat{v}(x,-y)-2\hat{v}(x, -2y) \text{ if } (x,y) \in \R^2 \setminus \R^2_+,\\
        \hat{v}(x,y) \hspace{3 cm} \text{ if } (x,y) \in \R^2_+.
    \end{cases}
\end{equation}
Clearly, $\overline{v} \in W^{2, p_0}(\R^2 ; \R^2)$ with 
\begin{equation*}
    \|\nabla^2 \overline{v}\|_{L^{ p_0}(\R^2)}\leq C\|\nabla^2 v \|_{L^{p_0}(A^l_r)},
\end{equation*}
for $C>0$ independent of $v,r$. 
Define $\Gamma^l_\infty:= \{(x,lx): x \geq 0\}$. By Theorem 18.24 in \cite{leoni2024first}, $\nabla \overline{v} \in L^{p_0/2-p_0}(\Gamma^l_\infty; \R^{2 \times 2})$ in the sense of traces and satisfies the inequality
\begin{equation}\label{eqn: nabla v bar trace bound}
    \|\nabla \overline{v}\|_{L^{p_0/(2-p_0)}(\Gamma^l_\infty)} \leq C \|\nabla^2 \overline{v}\|_{L^{ p_0}(\R^2)}\leq C \|\nabla^2 v \|_{L^{p_0}(A^l_r)}.
\end{equation}
By choice of $\varphi_r$, (\ref{eqn: def u-tilde}), (\ref{eqn: define v using utilde}), (\ref{eqn: def v hat}) and (\ref{eqn: def v bar}), $\overline{v} \circ \Psi= u-w_0$ in $A^l_{r/2}$. Hence, for a.e. $x \in I_{r/2}$, 
\begin{equation*}
    |\nabla (u-w_0)(x, h(x))|  \leq  C \sup_{(x,y) \in \Omega_h^{0,r/2}}(1+ |y\sigma'(x)|) |\nabla \overline{v}(x,lx)| 
\end{equation*}
Now, by (\ref{eqn: ysigma' bound_omega}), (\ref{eqn: h''(alpha)=h''(beta)=0}), (\ref{eqn: r h''' L^1 bound by constant}), (\ref{eqn: nabla^2 v L^p_0 bound}) and (\ref{eqn: nabla v bar trace bound}), we have 
\begin{eqnarray*}
    \|\nabla u-\nabla w_0\|_{L^{p_0/(2-p_0)}(\Gamma^{0, r/2}_h)} \leq C \|\nabla \overline{v}\|_{L^{p_0/2-p_0}(\Gamma^l_\infty)}\\
    \leq C  \|\nabla^2 v \|_{L^{p_0}(A^l_r)} \leq \frac{C}{r}\|\nabla u\|_{L^{p_0}(\Omega_h^{0,r})}+ C.
\end{eqnarray*}
Since $r<1$, by direct computation, we have $ \|\nabla w_0\|_{L^{p_0/(2-p_0)}(\Gamma^{0, r/2}_h)} \leq |e_0| (1+L_0^2)^{\frac{2-p_0}{2 p_0}}\leq C$ for some $C>0$. Hence, (\ref{eqn: bound on nabla u trace with nabla u}) follows.

\end{proof}
\begin{cor}\label{cor: 20}
    Under the hypotheses of Theorem \ref{thm: 6}, let $1 \leq p \leq p_0/ (2-p_0)$ and let 
    \begin{equation}\label{eqn: def r_1}
        r_1:= \min\bigg\{\frac{c_1}{2c_0 B_\tau (h,h_0, \alpha, \alpha_0, \beta, \beta_0)}, \frac{\delta_0}{2}, \frac{\eta_1}{2} \bigg\},
    \end{equation}
    where $c_0>0$ is the constant in the statement of Theorem \ref{thm: 6}, $0<c_1<1$ is the constant in (\ref{eqn: r h''' L^1 bound by constant}), $\eta_1$ is from (\ref{eqn: Lower bound on h_EL section-assumption}), $\delta_0$ is defined in (\ref{eqn: def delta_0}) and $B_\tau (h,h_0, \alpha, \alpha_0, \beta, \beta_0)$ is defined in (\ref{eqn: def B_tau}). Then, there exist two constants $c_3= c_3(\eta_0, \eta_1, M)$ and $c_4= c_4(\eta_0, \eta_1, M, p)>0$ such that 
    \begin{gather}
        \|\nabla^2 u\|_{L^{p_0}(\Omega_h^{\alpha, \alpha+r_1})} \leq c_3 B_\tau (h,h_0, \alpha, \alpha_0, \beta, \beta_0)^{2-(2/p_0)} ,\\
        \|\nabla u\|_{L^{p}(\Gamma_h^{\alpha, \alpha+r_1})} \leq c_4 B_\tau (h,h_0, \alpha, \alpha_0, \beta, \beta_0)^{1-(1/p)}. \label{eqn: nabla u trace L^p bound with B_tau}
    \end{gather}
\end{cor}
 \begin{proof}
  Using H\"older's inequality with exponent $q = 2/p_0 >1$, we have 
     \begin{eqnarray*}
          \|\nabla u\|_{L^{p_0}(\Omega_h^{\alpha, \alpha+r})}^{p_0} &=& \int_{\Omega_h^{\alpha, \alpha+r}} |\nabla u|^{p_0} dx dy \leq \|\nabla u\|_{L^2(\Omega_h)}\L^2 (\Omega_h^{\alpha, \alpha+r})^{(2-p_0)/2},
     \end{eqnarray*}
     where $\l^2$ denotes the Lebesgue measure in $\R^2$. Observe that $\L^2 (\Omega_h^{\alpha, \alpha+r}) \leq L_0 r^2$.  Using (\ref{eqn: u H^1 upper bound}) from Theorem \ref{thm: 2_existence of minimizer}, we have 
     \begin{equation}\label{eqn: nabla u corner bound with r}
         \|\nabla u\|_{L^{p_0}(\Omega_h^{\alpha, \alpha+r})} \leq C r^{(2/p_0) -1}.
     \end{equation}
     Choosing $r_1$ as in (\ref{eqn: def r_1}), using (\ref{eqn: h''' B_tau L 1 bound}), we have 
     \begin{equation*}
        2 r_1\|h'''\|_{L^1(I_r)} \leq c_1.
     \end{equation*}
     Therefore, we can apply Theorem \ref{thm: 18}, (\ref{eqn: def r_1}) and (\ref{eqn: nabla u corner bound with r}) to obtain 
     \begin{eqnarray*}
          \|\nabla^2 u\|_{L^{p_0}(\Omega_h^{\alpha, \alpha+r})}  \leq  c_2 + \frac{c_2 C}{r_1}r_1^{(2/p_0)-1} \leq \, C + \frac{C}{r_1^{2-(2/p_0)}}
          \leq CB_\tau (h,h_0, \alpha, \alpha_0, \beta, \beta_0)^{2-(2/p_0)},
     \end{eqnarray*}
     since $B_\tau (h,h_0, \alpha, \alpha_0, \beta, \beta_0) \geq 1$.
\paragraph{} Let $1\leq p< \frac{p_0}{2-p_0}$ and $p_1:= \frac{p_0}{p(2-p_0)}>1$. Then, by H\"older's inequality, 
\begin{eqnarray*}
    \int_{\Gamma_h^{\alpha, \alpha+r_1}} |\nabla u|^{p} d\h^1 &\leq& \bigg(\int_{\Gamma_h^{\alpha, \alpha+r_1}} |\nabla u|^{p_0/(2-p_0)} d\h^1\bigg)^{1/p_1} \bigg(\int_{\Gamma_h^{\alpha, \alpha+r_1}} 1 d\h^1\bigg)^{1/p_1'}.
\end{eqnarray*}
where $p_1'= \frac{p_0}{p_0-2p+pp_0}$ is the conjugate of $p_1$. Therefore, as $r_1<1$, using (\ref{eqn: bound on nabla u trace with nabla u}) and (\ref{eqn: nabla u corner bound with r}),
\begin{gather*}
    \|\nabla u\|_{L^{p}(\Gamma_h^{\alpha, \alpha+r_1})}  \leq  C \|\nabla u\|_{L^{p_0/(2-p_0)}(\Gamma_h^{\alpha, \alpha+r_1})}r_1^{1+\frac{1}{p}-\frac{2}{p_0}}
    \leq C+Cr_1^{\frac{1}{p}-\frac{2}{p_0}}\|\nabla u\|_{L^{p_0}(\Omega_h^{\alpha, \alpha+2r_1})}\\
    \leq C+Cr_1^{\frac{1}{p}-1}
    \leq C B_\tau (h,h_0, \alpha, \alpha_0, \beta, \beta_0)^{1-(1/p)}.
\end{gather*}
When $p=\frac{p_0}{2-p_0}$, $1-\frac{1}{p}= 2-\frac{2}{p_0}$. Hence, (\ref{eqn: nabla u trace L^p bound with B_tau}) follows from (\ref{eqn: bound on nabla u trace with nabla u}) and (\ref{eqn: nabla u corner bound with r}).

 \end{proof} 

For the rest of the section, we set $I_r(x):= (x-r, x+r)$ for $r>0$ and $x \in \R$. Under the assumptions of Theorem \ref{thm: 6}, fix $r$ and $x_0$ such that 
\begin{equation}\label{eqn: Fix r and x_0 for global regularity of nabla^2 u}
    0< 2r\leq \min\{\delta_0, \eta_1\} \text{ and } \alpha + 8r/\eta_0 \leq x_0 \leq \beta- 8r/ \eta_0.
\end{equation}
For $\alpha<a<b<\beta$ and $\eta>0$, recall that $ \Omega_{h,\eta}^{a,b}$ is as defined in (\ref{eqn: def Omega h,eta,a,b}). Define 
\begin{equation}\label{eqn: def Gamma h,a,b}
    \Gamma_h^{a,b}:= \{(x,h(x)): a<x<b\}.
\end{equation}
We state the following theorem (Theorem 5.10, \cite{dal2025motion}) without proof:
\begin{thm}\label{thm: 21}
    Under the assumptions of Theorem \ref{thm: 6}, let $r$ and $x_0$ be as in (\ref{eqn: Fix r and x_0 for global regularity of nabla^2 u}). Then, for every $2 \leq q < \infty$, there exists constants $c_5=c_5(\eta_0, \eta_1, M)>0$ and $c_6=c_6(\eta_0, \eta_1, M,q)>0$ independent of $h,r$ and $x_0$ such that, if 
    \begin{equation}\label{eqn: h'' L^infinity bound with 1/r}
        \|h''\|_{L^\infty(I_{4r}(x_0))} \leq \frac{1}{r},
    \end{equation}
    then $u \in H^2(\Omega_{h, 2r}^{x_0-2r, x_0+2r})$ and we have 
   \begin{eqnarray}
       \|\nabla^2 u\|_{L^2(\Omega_{h, 2r}^{x_0-2r, x_0+2r})} &\leq&\frac{c_5}{r}  \|\nabla u\|_{L^2(\Omega_{h, 4r}^{x_0-4r, x_0+4r})}, \label{eqn: nabla^2 u bound in interior with nabla u bound in bigger set}\\
        \|\nabla u\|_{L^q(\Gamma_{h}^{x_0-2r, x_0+2r})} &\leq& \frac{c_6}{r^{1-1/q}}\|\nabla u\|_{L^2(\Omega_{h, 4r}^{x_0-4r, x_0+4r})}. \label{eqn: nabla u bound on interior of graph of h with nabla u bound in bigger set}
   \end{eqnarray}
\end{thm}
Using the theorem above, we prove the following result.
\begin{thm}\label{thm: 23}
    Under the assumptions of Theorem \ref{thm: 6}, let $2\leq p \leq \frac{p_0}{2-p_0}$. Then, there exists a constant $c_7 = c_7(\eta_0, \eta_1, M, p)>0$ such that 
    \begin{equation}\label{eqn: L^p bound on nabla u in entire graph by B_tau}
        \|\nabla u\|_{L^p(\Gamma_h)} \leq c_7   B_\tau (h,h_0, \alpha, \alpha_0, \beta, \beta_0)^{1-1/p}.
    \end{equation}
\end{thm}
\begin{proof}
    Since $p$ satisfies the conditions of Corollary \ref{cor: 20},
 from (\ref{eqn: nabla u trace L^p bound with B_tau}), we have 
 \begin{equation}
     \int_\alpha^{\alpha+r_1} |\nabla u(x, h(x))|^p (1+(h'(x)^2)^{1/2} dx \leq c_4^p  B_\tau (h,h_0, \alpha, \alpha_0, \beta, \beta_0)^{p-1}.
 \end{equation}
  By following the same arguments, we also have the estimate above in the interval $(\beta-r_1, \beta)$. Now, we need to prove it in the interval $(\alpha+r_1, \beta-r_1)$. Let $r:= \frac{\eta_0 r_1}{4}$ and $x_i:= ir$ for every $i \in \Z$. Let $I_r(x_i)=(x_i-r, x_i+r)$. Define $Z \subseteq \Z$ as 
  \begin{equation*}
      Z := \{i \in \Z: \alpha+r_1<x_i<\beta-r_1\}.
  \end{equation*}
 Note that $(\alpha+r_1, \beta-r_1) \subseteq \bigcup\limits_{i \in Z}I_r(x_i)$. Therefore,
 \begin{equation}\label{eqn: interior boundary L^p estimate bound by sum of L^p estimates of parts}
      \int_{\alpha+r_1}^{\beta-r_1} |\nabla u(x, h(x))|^p  (1+(h'(x)^2)^{1/2} dx \leq \sum_{i \in Z} \bigg(\int_{I_r(x_i)} |\nabla u(x, h(x))|^p (1+(h'(x)^2)^{1/2} dx\bigg) .
 \end{equation}
By the choice of $r$ and (\ref{eqn: def r_1}), $r/2$ satisfies (\ref{eqn: Fix r and x_0 for global regularity of nabla^2 u}). Therefore, by Theorem \ref{thm: 21}, (\ref{eqn: nabla u bound on interior of graph of h with nabla u bound in bigger set}) is obtained in $\Gamma_h^{x_i-r, x_i+r}$ for every $i \in Z$.
\paragraph{} Recall that $c_1<1$ from the proof of Theorem \ref{thm: 18} and that $0< \eta_0<1$ from the assumptions of Theorem \ref{thm: 6}. Hence, by (\ref{eqn: h'' B_tau L infinity bound}), (\ref{eqn: def r_1}) and the choice of $r$,
\begin{gather*}
    \frac{r}{2}\|h''\|_{L^\infty(I_{2r}(x_0))} \leq  \frac{r}{2}c_0B_\tau (h,h_0, \alpha, \alpha_0, \beta, \beta_0)\\
    \leq r_1 \frac{\eta_0}{8}c_0B_\tau (h,h_0, \alpha, \alpha_0, \beta, \beta_0)
    \leq \,c_1  \frac{\eta_0}{8}\leq 1.
\end{gather*}
 Therefore, $ \frac{r}{2}$ satisfies (\ref{eqn: h'' L^infinity bound with 1/r}) and hence we can use (\ref{eqn: nabla u bound on interior of graph of h with nabla u bound in bigger set}) to obtain
 \begin{equation*}
     \int_{I_r(x_i)} |\nabla u(x, h(x))|^p (1+(h'(x)^2)^{1/2} dx\leq \frac{2^{p-1}c_6^p}{r^{p-1 }}\|\nabla u\|_{L^2(\Omega_{h, 2r}^{x_i-2r, x_i+2r})}^p,
 \end{equation*}
for every $i \in Z$. Note that $\Omega_{h, 2r}^{x_i-2r, x_i+2r}$ intersects at most 7 sets of the same form. Therefore, taking sum with respect to $i \in Z$ on both sides, by (\ref{eqn: u H^1 upper bound}) and (\ref{eqn: interior boundary L^p estimate bound by sum of L^p estimates of parts}), we have 
\begin{eqnarray}
     \int_{\alpha+r_1}^{\beta-r_1} |\nabla u(x, h(x))|^p (1+(h'(x)^2)^{1/2} dx &\leq &  \frac{C}{r^{p-1}} \|\nabla u\|_{L^2(\Omega_{h})}^p\nonumber \leq  \frac{C}{r_1^{p-1}}.
\end{eqnarray}
Since we have 
\begin{equation*}
    \frac{1}{r_1}\leq \max \bigg\{\frac{c_0}{c_1}, \frac{2}{\delta_0}, \frac{2}{\eta_1}\bigg\}B_\tau (h,h_0, \alpha, \alpha_0, \beta, \beta_0)
\end{equation*}
by (\ref{eqn: def r_1}), it follows that
\begin{equation*}
    \int_{\alpha+r_1}^{\beta-r_1} |\nabla u(x, h(x))|^p (1+(h'(x)^2)^{1/2} dx \leq C B_\tau (h,h_0, \alpha, \alpha_0, \beta, \beta_0)^{p-1}.
\end{equation*}
Therefore, (\ref{eqn: L^p bound on nabla u in entire graph by B_tau}) is proved.

\end{proof}
The proof of the next result follows the proof of Theorem 7.4 in \cite{dal2025motion}. Note that $B_\tau (h,h_0, \alpha, \alpha_0, \beta, \beta_0)$ is different from the $B_{\tau}(H,H_0, \alpha, \alpha_0, \beta, \beta_0)$ in \cite{dal2025motion}.
\begin{thm}\label{thm: 24}
    Under the hypothesis of Theorem \ref{thm: 6}, let $1 < p\leq \frac{p_0}{4-2p_0} $. Then, we have $h^{(iv)} \in L^{p_0/(4-2p_0)}((\alpha, \beta))$ and there exists a constant $c_8=c_8(\eta_0, \eta_1, M,p)>0$ such that 
    \begin{equation}\label{eqn: h^(iv) L^p bound}
        \|h^{(iv)}\|_{L^p((\alpha, \beta))}^p \leq c_8  B_\tau (h,h_0, \alpha, \alpha_0, \beta, \beta_0)^{q},
    \end{equation}
    where $ B_\tau (h,h_0, \alpha, \alpha_0, \beta, \beta_0)$ is defined in (\ref{eqn: def B_tau}) and 
    \begin{equation*}
        q=\max\bigg\{p, \frac{5p}{2}-1, 3p-2, 2p-1\bigg\}.
    \end{equation*}
\end{thm}
\begin{proof}
    By (\ref{eqn: EL strong equation}), we have 
   \begin{gather*}
       \nu_0\bigg( \frac{h''}{J^5}\bigg)''=  \gamma \bigg( \frac{h'}{J}\bigg)'-\frac{5}{2}\nu_0 \bigg(\frac{h'(h'')^2 }{J^7}\bigg)'-\overline{W} -\frac{1}{\tau}\frac{(\htil(x)-h_0(x))}{J_0}\chi_{[\alpha_0, \beta_0]}+m,
   \end{gather*}
  where $m$ is defined in (\ref{eqn: def m EL section}). Hence, for $p>1$, we have the point-wise estimate 
   \begin{equation} \label{eqn: h^4 L^p bound precursor}
       |h^{(iv)}|^p \leq C(|h''|^p+ |h''|^p|h'''|^p+|h''|^{3p}+ \overline{W}^p+\frac{1}{\tau^p}|\htil-h_0|^p+|m|^p ),
   \end{equation}
   where $C>0$ is a constant. By (\ref{eqn: def B_tau}), we have that 
   \begin{equation*}
       \frac{1}{\tau^p}|\htil-h_0|^p \leq B_\tau (h,h_0, \alpha, \alpha_0, \beta, \beta_0)^p. 
   \end{equation*}
   Following the proof of Theorem 7.4 in \cite{dal2025motion}, we can prove that 
  \begin{gather}
      \int_\alpha^\beta |h''|^p dx \leq C B_\tau (h,h_0, \alpha, \alpha_0, \beta, \beta_0)^p, \label{eqn: h'' L^p bound}\\
       \int_\alpha^\beta |h''|^p|h'''|^p dx \leq C B_\tau (h,h_0, \alpha, \alpha_0, \beta, \beta_0)^{\max\{5p/2-1, 3p-2\}}, \label{eqn: h'' h''' L^p bound} \\
        \int_\alpha^\beta |h''|^{3p} dx \leq C B_\tau (h,h_0, \alpha, \alpha_0, \beta, \beta_0)^{3p-2}, \text{ and} \label{eqn: h'' L^3p bound}\\
        \int_\alpha^\beta \overline{W}^p dx \leq  C B_\tau (h,h_0, \alpha, \alpha_0, \beta, \beta_0)^{2p-1}. \label{eqn: W bar L^p bound}
 \end{gather}
By (\ref{eqn: m B_tau bound}), we have 
 \begin{equation}\label{eqn: m^p bound}
      |m|^p \leq  C B_\tau (h,h_0, \alpha, \alpha_0, \beta, \beta_0)^p.
 \end{equation}
 Hence, (\ref{eqn: h^(iv) L^p bound}) follows from (\ref{eqn: Upper bound on beta-alpha}) and (\ref{eqn: h^4 L^p bound precursor})-(\ref{eqn: m^p bound}).
 
\end{proof}
\section{Discretization in Time}\label{sec: Discretization}
For the rest of the paper, we use the notation 
\begin{eqnarray*}
    \dot{f}(x,t)= \frac{\pa f}{\pa t} (x,t), \quad f'(x,t)= \frac{\pa f}{\pa x} (x,t)
\end{eqnarray*}
for functions $f: \R \times \R \rightarrow \R$. 
In this section, we use the existence results proved in Section \ref{sec: EL} to set up the minimizing movements scheme, by time discretization and recursive minimization. We fix $(\alpha_0, \beta_0, h_0, u_0) \in \a$ with $\text{Lip}(h_0)<L_0$. For $k \in \N$, set $\tau_k := \frac{1}{k}$, $t^i_k := i \tau_k$ for $i \in \N \cup \{0\},\, i\leq k$, and let
\begin{equation}
    (\alpha_k^0, \beta_k^0, h_k^0, u_k^0):= (\alpha_0, \beta_0, h_0, u_0).
\end{equation}
Given $(\alpha_k^{i-1}, \beta_k^{i-1}, h_k^{i-1}, u_k^{i-1}) \in \a$, by Theorem \ref{thm: 2_existence of minimizer}, there exists $(\alpha_k^i, \beta_k^i, h_k^i, u_k^i) \in \a$ which minimizes the functional
\begin{eqnarray}
    \f_k^{i-1}(\alpha, \beta, h, u):= \s(\alpha, \beta, h)+\e(\alpha, \beta, h,u)+\t_{\tau_k}(\alpha, \beta, h; \alpha_k^{i-1}, \beta_k^{i-1}, h_k^{i-1}).
\end{eqnarray}
We interpolate $\alpha, \beta, h$ linearly in time to obtain
\begin{eqnarray}
    \alpha_k(t) &:=& \alpha_k^{i-1}+ \frac{t-t_k^{i-1}}{\tau_k}(\alpha_k^i-\alpha_k^{i-1}),\label{eqn: def alpha_k}\\
    \beta_k(t)&:=& \beta_k^{i-1}+ \frac{t-t_k^{i-1}}{\tau_k}(\beta_k^i-\beta_k^{i-1}),\label{eqn: def beta_k}\\
    h_k(t,x)&:=& \htil_k^{i-1}(x)+ \frac{t-t_k^{i-1}}{\tau_k}(\htil_k^{i}(x)-\htil_k^{i-1}(x)).\label{eqn: def h_k}
\end{eqnarray}
for $t \in [t_k^{i-1}, t_k^i]$, $1\leq i \leq k$ and $x \in \R$, where $\htil_k^i$ is the extension of $h^i_k$ by $0$ outside the interval $[\alpha_k^i, \beta_k^i]$. Now, for $t \in (t_k^{i-1}, t_k^i]$ and $x \in \R$ we define  the constant-in-time interpolations 
\begin{equation}\label{eqn: piecewise constant extension of alpha, beta, h}
    \hat{\alpha}_k(t) := \alpha_k^i, \,\,\hat{\beta}_k(t):= \beta_k^i,\,\, \hat{h}_k(t,x):= \htil_k^i(x).
\end{equation}
For $t=0$, set $\hat{\alpha}_k(0)=\alpha_0$, $\hat{\beta}_k(0)=\beta_0$ and $\hat{h}_k(0,x)= \htil_0(x)$ for $x \in \R$.
Clearly, $(\hat{\alpha}_k(t), \hat{\beta}_k(t), \hat{h}_k(\cdot,t)) \in \a_s$ for every $t\geq 0$, $\hat{h}_k(\cdot,t) \in H^2((\hat{\alpha}_k(t), \hat{\beta}_k(t)))$ and by (\ref{eqn: beta-alpha lower bound lemma 1}),
\begin{equation}\label{eqn: lower bound on beta hat_k -alpha hat_k}
    \hat{\beta}_k(t)- \hat{\alpha}_k(t) \geq \sqrt{\frac{2A_0}{L_0}}.
\end{equation}
Moreover, both  linear and constant-in-time interpolations of $h$ are Lipschitz continuous in the $x$ variable with Lipschitz constant $L_0$ and preserve the area under them as $A_0$, at every fixed time point. 
\paragraph{} Recall the definition of $a_0, \,b_0$, $\phi_0$ and  $m$ from (\ref{eqn: def a_0, b_0}), (\ref{eqn: def phi_0 from psi_0}), and (\ref{eqn: def m EL section}). For every $\alpha^i_k$, $\beta^i_k$, we choose $a^i_{0,k}$, $b^i_{0,k}$ and $\phi^i_{0,k}$ such that $a^i_{0,k}=\alpha^i_k + \frac{\beta^i_k-\alpha^i_k}{4}$ and $b^i_{0,k}=\beta^i_k - \frac{\beta^i_k-\alpha^i_k}{4}$. We define $\phi^i_{0,k} \in C_c^\infty ([a^i_{0,k},b^i_{0,k}])$ as 
\begin{equation}
    \phi_{0,k}^i(x):=\frac{1}{b^i_{0,k}-a^i_{0,k}} \psi_0 \bigg( \frac{x-a^i_{0,k}}{b^i_{0,k}-a^i_{0,k}} \bigg),
\end{equation}
where $\psi_0$ is defined in (\ref{eqn: def psi_0}). Clearly, $\phi^i_{0,k} $ satisfies (\ref{eqn: phi_0 integral and L^infinity bound}) and (\ref{eqn: phi_0', phi_0'' L^infinity bound}). Define
\begin{equation}\label{eqn: def m^i_k}
    m^i_k:= \int_{\alpha^i_k}^{\beta^i_k} A^i_k(h^i_k)''(\phi^i_{0,k})''+ B^i_k (h^i_k)'(\phi^i_{0,k})'+f^i_k\phi^i_{0,k},
\end{equation}
where
\begin{gather*}
    A^i_k(x):= \frac{\nu_0}{(1+((\htil^{i}_{k})'(x))^2)^{5/2}},\\
    B^i_k(x):=  \frac{\gamma}{1+((\htil^{i}_{k})'(x))^2)^{1/2}}-\frac{5}{2} \nu_0  \frac{((\htil^i_k)''(x))^2 }{(1+((\htil^{i}_{k})'(x))^2)^{7/2}}, \text{ and}\\
    f^i_k(x):= W (Eu^i_k(x,\htil^i_k(x)))+\frac{1}{\tau_k}\frac{(\htil^i_k(x)-h^{i-1}_k(x))}{(1+((h^{i-1}_{k})'(x))^2)^{1/2}}\chi_{[\alpha^{i-1}_k, \beta^{i-1}_k]},
\end{gather*}
for $x \in \R$, $k \in \N$ and $1 \leq i \leq k$, $i \in \N$. Define  the constant in time interpolations
\begin{equation}\label{eqn: def a_0k hat, b_0k hat, phi_0,k hat, m_k hat}
   \hat{a}_{0,k}(t):= a^i_{0,k}, \quad \hat{b}_{0,k}(t):= b^i_{0,k}, \quad \hat{\phi}_{0,k}(t,x):= \phi^i_{0,k}(x), \quad \hat{m}_k(t):= m^i_k
\end{equation}
for $t \in (t^{i-1}_k, t^i_k]$ and $x \in \R$. Hence, by (\ref{eqn: lower bound on beta hat_k -alpha hat_k}), we have 
\begin{equation} \label{eqn: b_0k hat-a_0k hat lower bound}
    \hat{b}_{0,k}(t)-\hat{a}_{0,k}(t) \geq \sqrt{\frac{A_0}{2L_0}}
\end{equation}
for every $t>0$. Also, by (\ref{eqn: m B_tau bound}) and (\ref{eqn: lower bound on beta hat_k -alpha hat_k}), we have 
\begin{equation}
    |m^i_k| \leq C B_\tau (h^i_k,h^{i-1}_k, \alpha^i_k, \alpha^{i-1}_k, \beta^i_k, \beta^{i-1}_k)
\end{equation}
for $C$ independent of $\alpha^i_k$, $\beta^i_k$ and $h^i_k$. This is crucial for the uniform regularity estimates on $h^i_k$ and $u^i_k$.
\paragraph{} By Lemma 8.1 and Proposition 8.2 in \cite{dal2025motion}, there exist constants $M_0, M_1>0$ such that 
 \begin{gather}
        \s(\alpha^i_k, \beta^i_k, h^i_k)+ \e(\alpha^i_k, \beta^i_k, h^i_k, u^i_k)+ \sum_{j=1}^i \t_{\tau_k}(\alpha^j_k, \beta^j_k, h^j_k; \alpha^{j-1}_k, \beta^{j-1}_k, h^{j-1}_k) \leq M_0,\label{eqn: discrete energy upper bound}\\
    \frac{1}{\tau_k} \sum_{j=1}^\infty (\alpha^j_k - \alpha^{j-1}_k)^2 \leq M_1,\,\, \frac{1}{\tau_k} \sum_{j=1}^\infty (\beta^j_k - \beta^{j-1}_k)^2 \leq M_1,\label{eqn: proof of derivative of alpha bound}\\
    \int_0^\infty (\dot{\alpha}_k(t))^2 dt \leq M_1, \text{ and }  \int_0^\infty (\dot{\beta}_k(t))^2 dt \leq M_1,\\
    |\alpha_k(t_2)- \alpha_k(t_1)| \leq M_1^{1/2}|t_2-t_1|^{1/2}, \,\, |\beta_k(t_2)- \beta_k(t_1)| \leq M_1^{1/2}|t_2-t_1|^{1/2}, \label{eqn: alpha, beta C^{0, 1/2} bound}\\
        \alpha_0-(tM_1)^{1/2} \leq \alpha_k(t) \leq \beta_k(t) \leq \beta_0 + (tM_1)^{1/2} \label{eqn: upper bound on |beta_k-alpha_k|}
    \end{gather}
  for every $i \leq k \in \N$ and $t,t_1, t_2 \in [0, \infty)$. Using the bounds above, we prove the following proposition.
\begin{prop}\label{prop: 27}
  There exists $M_2>0$ such that for every $k \in \N$, 
  \begin{equation}\label{eqn: bound on time derivative of h_k}
    \int_0^\infty \int_\R |\dot{h}_k(t,x)|^2 dx dt \leq M_2.
\end{equation}
In particular,
\begin{eqnarray} 
    \|h_k(t_2,\cdot)-h_k(t_1,\cdot)\|_{L^2(\R)} &\leq& M_2^{1/2}|t_2-t_1|^{1/2}\label{eqn: bound on L^2 norm of h_k(t_2,.)-h_k(t_1,.)}\\
    \|\hat{h}_k(t_2,\cdot)-\hat{h}_k(t_1,\cdot)\|_{L^2(\R)} &\leq& M_2^{1/2}(|t_2-t_1|^{1/2}+\tau_k^{1/2})\label{eqn: bound on L_2 norm of h_k hat(t_1)-h_k hat(t_2)}\\
    \|h_k(t, \cdot)\|_{L^2(\R)} &\leq& M_2^{1/2} t^{1/2}+\|\tilde{h}^0\|_{L^2(\R)} \label{eqn: bound on L^2 norm of h_k(t,.)}
\end{eqnarray}
for every $t, t_1, t_2 \in [0, \infty)$ and every $k \in \N$.
\end{prop}

\begin{proof}

As $ \dot{h}_k(t,x) = \frac{1}{\tau_k}(\tilde{h}^j_k(x)-\tilde{h}^{j-1}_k(x))$ for $t \in (t^{j-1}_k, t^j_k)$ and $x \in \R$,
\begin{gather*}
    \int_0^\infty \int_\R |\dot{h}_k(t,x)|^2 dx dt = \sum_{j=1}^\infty \int_{t^{j-1}_k}^{t^j_k} \frac{1}{\tau_k^2} \int_\R |\tilde{h}^j_k(x)-\tilde{h}^{j-1}_k(x)|^2 dx dt\\
    = \sum_{j=1}^\infty \frac{1}{\tau_k} \int_{\alpha^{j-1}_k}^{\beta^{j-1}_k} |\tilde{h}^j_k(x)-\tilde{h}^{j-1}_k(x)|^2 dx  
     +\sum_{j \in A_k} \frac{1}{\tau_k} \int_{\alpha^{j}_k}^{\alpha^{j-1}_k} |\tilde{h}^j_k(x)|^2 dx\\
    +\sum_{j \in B_k} \frac{1}{\tau_k} \int_{\beta^{j-1}_k}^{\beta^{j}_k} |\tilde{h}^j_k(x)|^2 dx
    =: I_1+I_2+I_3,
\end{gather*}
where
\begin{gather*}
    A_k:=\{j \in \N: \alpha^j_k<\alpha^{j-1}_k\}, \quad B_k:= \{j \in \N: \beta^j_k>\beta^{j-1}_k\}.
\end{gather*}
Consider $I_2$. Using (\ref{eqn: proof of derivative of alpha bound}), (\ref{eqn: alpha, beta C^{0, 1/2} bound}), the uniform Lipschitz continuity of the functions $h^i_k$ and the fact that $\htil^j_k(\alpha^j_k)=0$, , we have
\begin{gather*}
    I_2 = \sum_{j \in A_k}\frac{1}{\tau_k} \int_{\alpha^{j}_k}^{\alpha^{j-1}_k} |h^j_k(x)-h^j_k(\alpha^j_k)|^2 dx \leq \sum_{j=1}^\infty \frac{1}{\tau_k} L_0^2 |\alpha^j_k - \alpha^{j-1}_k|^2 |\alpha^j_k - \alpha^{j-1}_k|\\
    \leq \sum_{j=1}^\infty \frac{1}{\tau_k} L_0^2 |\alpha^j_k - \alpha^{j-1}_k|^2 M_1^{1/2} \tau_k^{1/2}
    = L_0^2 M_1^{1/2} \tau_k^{1/2} \sum_{j=1}^\infty \frac{1}{\tau_k} |\alpha^j_k - \alpha^{j-1}_k|^2 \\\leq L_0^2 M_1^{3/2} \tau_k^{1/2} \leq L_0^2 M_1^{3/2}=: M_2^0.
\end{gather*}
Clearly, we can prove $I_3 \leq M_2^0$ following the same steps as above.
\paragraph{}Now, consider $I_1$. Using the upper bound in (\ref{eqn: discrete energy upper bound}) and the non-negativity of $\s$, $\e$ and $\t_{\tau_k}$ from (\ref{eqn: surface energy is non-negative}), (\ref{eqn: W upper and lower bounds by|.|^2}) and (\ref{eqn: define incremental energy functional}), we have 
\begin{equation*}
    \sum_{j=1}^i \frac{1}{2 \tau_k} \int_{\alpha^{j-1}_k}^{\beta^{j-1}_k} \frac{(\tilde{h}^{j}_k-h^{j-1}_k)^2}{\sqrt{1+((h^{j-1}_k)')^2}}dx \leq M_0
\end{equation*}
for every $i$. Therefore, as $\|(\tilde{h}^{j-1}_k)'\|_{L^{\infty}(\R)} \leq L_0$ for every $j, k$, 
\begin{eqnarray*}
   \frac{1}{ \tau_k}  \sum_{j=1}^\infty \int_{\alpha^{j-1}_k}^{\beta^{j-1}_k} \frac{(\tilde{h}^{j}_k-h^{j-1}_k)^2}{\sqrt{1+L_0^2}} dx  & \leq & \frac{1}{\tau_k} \sum_{j=1}^\infty \int_{\alpha^{j-1}_k}^{\beta^{j-1}_k} \frac{(\tilde{h}^{j}_k-h^{j-1}_k)^2}{\sqrt{1+((h^{j-1}_k)')^2}}dx \\
   &\leq & 2M_0.
\end{eqnarray*}
Hence, $I_1= \frac{1}{ \tau_k}  \sum_{j=1}^\infty \int_{\alpha^{j-1}_k}^{\beta^{j-1}_k} (\tilde{h}^{j}_k-h^{j-1}_k)^2 dx \leq 2M_0 \sqrt{1+ L_0^2}=: M_2^1$. Therefore, we have
\begin{equation*}
    \int_0^\infty \int_\R |\dot{h}_k(t,x)|^2 dx dt \leq M_2 =: 2M^0_2+ M^1_2.
\end{equation*}

Now, using the fundamental theorem of calculus, H\"older's inequality and Fubini's theorem, given $t_1, t_2 \in [0, \infty)$,
\begin{gather*}
    \int_\R |h_k(t_2,x)- h_k(t_1,x)|^2 dx = \int_\R \bigg|\int_{t_1}^{t_2} \dot{h}_k(t,x) dt \bigg|^2 dx 
    \leq \int_\R \int_{t_1}^{t_2} |\dot{h}_k(t,x)|^2 dt \,\, |t_2-t_1| dx\\
    \leq\int_{t_1}^{t_2}  \int_\R |\dot{h}_k(t,x)|^2  dx dt \,\,|t_2-t_1|
    \leq M_2 |t_2-t_1|.
\end{gather*}
This proves (\ref{eqn: bound on L^2 norm of h_k(t_2,.)-h_k(t_1,.)}). Now, consider $t^{i_1}_k, t^{i_2}_k$ such that $t_1 \in (t^{i_1-1}_k, t^{i_1}_k]$ and $t_2 \in (t^{i_2-1}_k, t^{i_2}_k]$. Then, by (\ref{eqn: def h_k}), (\ref{eqn: piecewise constant extension of alpha, beta, h}) and (\ref{eqn: bound on L^2 norm of h_k(t_2,.)-h_k(t_1,.)}), we have 
\begin{eqnarray*}
     \|\hat{h}_k(t_2,\cdot)-\hat{h}_k(t_1,\cdot)\|_{L^2(\R)} &\leq& \|h_k(t^{i_2}_k,\cdot)-h_k(t^{i_1}_k,\cdot)\|_{L^2(\R)} \\
     &\leq& M_2^{1/2}|t^{i_2}_k-t^{i_1}_k|\leq M_2^{1/2}(|t_2-t_1|^{1/2}+\tau_k^{1/2}),
\end{eqnarray*}
thus proving (\ref{eqn: bound on L_2 norm of h_k hat(t_1)-h_k hat(t_2)}). Further, as $h_k(0,x)= \tilde{h}_0(x)$ for every $ x \in \R$ and $k \in \N$, we have $\|h_k(t,\cdot)\|_{L^2(\R)} \leq M_2^{1/2} t^{1/2}+ \|\tilde{h}_0\|_{L^2(\R)}$ for every $t \in [0, \infty)$, proving (\ref{eqn: bound on L^2 norm of h_k(t,.)}).

\end{proof}
For $t \in (t^{i-1}_k, t^i_k]$, define
\begin{equation}\label{eqn: def alpha_k bar, beta_k bar}
    \bar{\alpha}_k(t):= \max \{\alpha^{i-1}_k, \alpha^i_k\}, \quad \bar{\beta}_k(t):= \min \{\beta^{i-1}_k, \beta^i_k\}
\end{equation}
for $i,k \in \N$ and $1 \leq i \leq k$. Let $\bar{\alpha}(0)=\alpha_0$ and $\bar{\beta}(0)=\beta_0$.  
\begin{prop}\label{prop: H^2 regularity of h_k(t,.)}
    For every $k \in \N$ and $t \in [0, \infty)$, $h_k(t, \cdot) \in H^2((\bar{\alpha}_k(t), \bar{\beta}_k(t)))$. In particular, there exists $M_4>0$ such that 
    \begin{equation} \label{eqn: upper bound on H^2 norm of h_k(t,.)}
        \|h_k(t, \cdot)\|_{H^2((\bar{\alpha}_k(t), \bar{\beta}_k(t)))} \leq M_4(1+t^{1/4}+ t^{1/2})
    \end{equation}
for every $k \in \N$ and  $t \in [0, \infty)$.
\end{prop}
\begin{proof}
    By (\ref{eqn: def h_k}), 
    \begin{equation*}
        h_k'(t,x)= (\htil_k^{i-1})' (x) + \frac{t-t^{i-1}_k}{\tau_k}(\htil^i_k)'(x)
    \end{equation*}
for $x \in \R$ and $t \in (t^{i-1}_k, t^i_k)$. Therefore, since $|(h^j_k)'(x)|\leq L_0$ for every $j \leq k$, $k \in \N$ and $x \in \R$, using (\ref{eqn: upper bound on |beta_k-alpha_k|}),
\begin{eqnarray}
   \bigg( \int_{\alpha_k(t)}^{\beta_k(t)}|h_k'(t,x)|^2 dx\bigg)^{1/2} & \leq & \bigg( \int_{\alpha^{i-1}_k}^{\beta^{i-1}_k} |(\htil_k^{i-1})' (x)|^2 dx\bigg)^{1/2} +\bigg(\int_{\alpha^{i}_k}^{\beta^{i}_k} |(\htil_k^{i})' (x)|^2 dx\bigg)^{1/2}\nonumber\\
    &\leq& 2L_0(\beta_0- \alpha_0 + 2t^{1/2}M_1^{1/2})^{1/2} \nonumber\\
    &\leq& 2 L_0 ((\beta_0- \alpha_0 )^{1/2}+ 2^{1/2}M_1^{1/4}t^{1/4}),
\end{eqnarray}
where we used the fact that $\htil^i_k =0 $ outside $(\alpha^i_k, \beta^i_k)$. 
Using the same ideas as in the proof of Proposition 8.5 in \cite{dal2025motion}, there exists a constant $M_3>0$ such that for every $i, k$
 we have
  \begin{equation} \label{eqn:uni L^2 bound h^i_k''}
       \int_{\alpha^i_k}^{\beta^i_k} |(h^i_k)''(x)|^2 dx \leq M_3.
   \end{equation}
    Hence, by (\ref{eqn: def h_k}),
\begin{eqnarray}
    \bigg( \int_{\bar{\alpha}_k(t)}^{\bar{\beta}_k(t)} |h_k''(t,x)|^2 dx \bigg)^{1/2} & \leq & \bigg( \int_{\alpha^{i-1}_k}^{\beta^{i-1}_k} |(\htil_k^{i-1})'' (x)|^2 dx\bigg)^{1/2} +\bigg(\int_{\alpha^{i}_k}^{\beta^{i}_k} |(\htil_k^{i})'' (x)|^2 dx\bigg)^{1/2}\nonumber\\
    &\leq& 2M_3^{1/2},
\end{eqnarray}
where we used the fact that $\htil^i_k =0 $ outside $(\alpha^i_k, \beta^i_k)$ again. By (\ref{eqn: bound on L^2 norm of h_k(t,.)}), we already have $ \|h_k(t, \cdot)\|_{L^2(\R)} \leq M_2^{1/2} t^{1/2}+\|\tilde{h}^0\|_{L^2(\R)}$. Choose $M_4:= \max\{2M_3^{1/2}+2 L_0 (\beta_0- \alpha_0 )^{1/2}+\|\tilde{h}^0\|_{L^2(\R)}, M_2^{1/2}, 2^{3/2}L_0M_1^{1/4}\}$ to prove the claim (\ref{eqn: upper bound on H^2 norm of h_k(t,.)}).

\end{proof}
\begin{remark}
   Note that by (\ref{eqn:uni L^2 bound h^i_k''}),  
   \begin{equation}
       |(h^i_k)'(x_2)-(h^i_k)'(x_1)| \leq M_3^{1/2}(x_2-x_1)^{1/2}
   \end{equation}
   for all $\alpha^i_k \leq x_1 \leq x_2 \leq \beta^i_k$.
\end{remark}

\section{Convergence}
In this section, we prove the convergence of the sequences of functions obtained in the previous section through interpolations. These convergence results lead us to the proof of Theorem \ref{thm: main thm}. We assume the following throughout this section.
 \begin{gather}
       (\alpha_0, \beta_0, h_0) \in \a_s, \label{assumption: alpha_0, beta_0, h_0 are in A_s}\\ 
        h_0(x)>0 \text{ for every } x \in (\alpha_0, \beta_0),\\
        h_0'(\alpha_0)>0 \text{ and } h_0'(\beta_0)<0,\\
        \text{Lip }h_0<L_0. \label{assumption: Lip h_0< L_0}
    \end{gather}
 We state the following result from Proposition 9.1 in \cite{dal2025motion} without proof.
\begin{prop}\label{prop: 30}
  Let $\alpha_k, \beta_k, \hat{\alpha}_k, \hat{\beta}_k$ be as defined in (\ref{eqn: def alpha_k}), (\ref{eqn: def beta_k}) and (\ref{eqn: piecewise constant extension of alpha, beta, h}). Then there exist functions $\alpha, \beta : (0, \infty) \rightarrow \R$ such that $\alpha, \beta \in H^1((0,T))$ for every $T>0$ and up to subsequences (not relabeled),
\begin{gather}
    \alpha_k \rightharpoonup \alpha \text{ and } \beta_k \rightharpoonup \beta \text{ weakly in } H^1((0,T)), \label{eqn: alpha, beta weak conv}\\
    \alpha_k \rightarrow \alpha \text{ and } \beta_k \rightarrow \beta \text{ uniformly in } [0,T], \label{eqn: uniform conv_alpha, beta}\\
    \hat{\alpha}_k \rightarrow \alpha \text{ and } \hat{\beta}_k \rightarrow \beta \text{ uniformly in } [0,T], \label{eqn: uniform conv_alpha hat, beta hat}
\end{gather}
for every $T>0$. Moreover,
\begin{equation}\label{eqn: alpha(0), beta(0), beta-alpha bdd below}
    \alpha(0)= \alpha_0, \,\beta(0)=\beta_0, \text{ and } \beta(t)-\alpha(t) \geq \sqrt{\frac{2A_0}{L_0}}
\end{equation}
    for every $t\geq 0$.
\end{prop} 
Now, we prove a corollary.
\begin{cor} \label{cor: a_0, b_0, phi_0 conv}
    Let $\hat{a}_{0,k}, \hat{b}_{0,k}, \hat{\phi}_{0,k}$ be as defined in (\ref{eqn: def a_0k hat, b_0k hat, phi_0,k hat, m_k hat}). Then, 
    \begin{itemize}
        \item[(i)] there exist functions $a_0, b_0 : (0, \infty) \rightarrow \R$ such that $a_0, b_0 \in H^1((0,T))$ for every $T>0$ and up to a subsequence (not relabeled),
    \begin{gather}\label{eqn: uniform conv_a_0, b_0}
        \hat{a}_{0,k} \rightarrow a_0 \text{ and } \hat{b}_{0,k} \rightarrow b_0 \text{ uniformly in } [0,T],
    \end{gather}
    for every $T>0$. Moreover, 
    \begin{equation}\label{eqn: b_0-a_0 function lower bound}
        b_0(t)-a_0(t) \geq \sqrt{\frac{A_0}{2L_0}}
    \end{equation}
    for every $t\geq 0$.
    \item[(ii)] There exists $\phi_0 : (0, \infty) \times \R \rightarrow \R$ such that $\phi_0(t, \cdot) \in C_c^\infty(\R)$ with $\operatorname{supp}(\phi_0(t, \cdot) \subset (a_0(t), b_0(t)))$ for every $t>0$ and up to a subsequence (not relabeled),
    \begin{gather}
        \hat{\phi}_{0,k} \rightarrow \phi_0, \quad  \hat{\phi}'_{0,k} \rightarrow \phi_0', \quad \hat{\phi}_{0,k}'' \rightarrow \phi_0''  \text{ uniformly in } [0,T] \times \R, \label{eqn: phi_0,k hat with two derivatives uniform convergence}
    \end{gather}
    for every $T>0$.
    \end{itemize}
\end{cor}
\begin{proof}
\begin{itemize}
        \item[(i)] Define $a_0, b_0 : (0, \infty) \rightarrow \R$ as
    \begin{equation}
        a_0(t):= \alpha(t)+ \frac{\beta(t)-\alpha(t)}{4}, \quad b_0(t):=\beta(t)-\frac{\beta(t)-\alpha(t)}{4}
    \end{equation} 
    for $t \in (0, \infty)$. Then, (\ref{eqn: uniform conv_a_0, b_0}) and (\ref{eqn: b_0-a_0 function lower bound}) follow from (\ref{eqn: def a_0k hat, b_0k hat, phi_0,k hat, m_k hat}), (\ref{eqn: b_0k hat-a_0k hat lower bound}), (\ref{eqn: uniform conv_alpha hat, beta hat}) and (\ref{eqn: alpha(0), beta(0), beta-alpha bdd below}). $a_0, b_0 \in H^1((0,T))$ for every $T>0$ since $\alpha, \beta \in H^1((0,T))$ for every $T>0$ by Proposition \ref{prop: 30}.
    \item[(ii)]  Note that
    \begin{gather}
        \hat{\phi}_{0,k}'(t,x)= \frac{1}{(\hat{b}_{0,k}(t)-\hat{a}_{0,k}(t))^2}\psi_0'\bigg(\frac{x-\hat{a}_{0,k}(t)}{\hat{b}_{0,k}(t)-\hat{a}_{0,k}(t)} \bigg),\\
        \hat{\phi}_{0,k}''(t,x)= \frac{1}{(\hat{b}_{0,k}(t)-\hat{a}_{0,k}(t))^3}\psi_0''\bigg(\frac{x-\hat{a}_{0,k}(t)}{\hat{b}_{0,k}(t)-\hat{a}_{0,k}(t)} \bigg),
    \end{gather}
  where $\psi_0$ is as in (\ref{eqn: def psi_0}). Now, define $\phi_0 : (0, \infty) \times \R \rightarrow \R$ as 
    \begin{equation}\label{eqn: def phi_0(t,x)}
        \phi_0(t,x):= \frac{1}{b_0(t)-a_0(t)} \psi_0 \bigg( \frac{x-a_0(t)}{b_0(t)-a_0(t)}\bigg)
    \end{equation}
    for $t \in (0, \infty)$ and $x \in \R$.
    Since $\psi_0 \in C_c^\infty(\R)$ with $\operatorname{supp}(\psi) \subset (0,1)$,  $\phi_0(t, \cdot) \in C_c^\infty(\R)$ with $\operatorname{supp}(\phi_0(t, \cdot)) \subset (a_0(t), b_0(t))$. Further,
    \begin{gather}
        \phi_0'(t,x) =\frac{1}{(b_0(t)-a_0(t))^2} \psi_0' \bigg( \frac{x-a_0(t)}{b_0(t)-a_0(t)}\bigg), \label{eqn: phi_0'(t,.)}\\
        \phi_0''(t,x)= \frac{1}{(b_0(t)-a_0(t))^3} \psi_0'' \bigg( \frac{x-a_0(t)}{b_0(t)-a_0(t)}\bigg). \label{eqn: phi_0''(t,.)}
    \end{gather}
    Then, (\ref{eqn: phi_0,k hat with two derivatives uniform convergence}) follows from (\ref{eqn: uniform conv_a_0, b_0}).
   \end{itemize} 
\end{proof}

\begin{prop}\label{prop:31}
    Let $h_k$ be defined as in (\ref{eqn: def h_k}). Then there exist a subsequence (not relabeled) and a non-negative function $h_* \in C^{0, 1/2}([0, \infty); L^2(\R))$ such that for every $t \in [0, \infty)$, 
    \begin{gather}
        {\operatorname*{Lip}}\, h_* (t,\cdot) \leq L_0 \text{ and } h_* (t,x)=0 \text{ for } x \notin (\alpha(t), \beta(t)), \label{eqn: Lip h* and h* support} \\
        h_k (t,\cdot) \rightarrow h_{*}(t,\cdot) \text{ uniformly in } \R \label{eqn: uniform convergence of h_k in R},\\
        h_k' (t,\cdot) \overset{\ast}{\rightharpoonup} h_{*}'(t,\cdot) \text{ weakly star in } L^\infty(\R). \label{eqn: weak star convergence of h_k'(t,.) to h_* (t,.)}
    \end{gather}
Moreover, if we define $h(t,\cdot)$ as the restriction of $h_*(t,\cdot)$ to $(\alpha(t), \beta(t))$, then $h(t,\cdot) \in H^1_0((\alpha(t), \beta(t)))\cap H^2((\alpha(t), \beta(t)))$ and 
\begin{gather}
    \int_{\alpha(t)}^{\beta(t)} h(t,x) dx = A_0,\label{eqn: integral of h(t,.)=A_0}\\
     \int_{\alpha(t)}^{\beta(t)} |h''(t,x)|^2 dx \leq M_3 \label{eqn: L^2 bound on h''},\\
     |h'(t,x_2)-h'(t,x_1)| \leq M_3^{1/2} (x_2-x_1)^{1/2} \text{ for all } \alpha(t) \leq x_1 \leq x_2 \leq \beta(t),\label{eqn: C^{0,1/2} seminorm bound on h'(t,.) in (alpha(t), beta(t))}
\end{gather}
for every $t \in [0, \infty)$, where $M_3>0$ is the constant in (\ref{eqn:uni L^2 bound h^i_k''}).
\end{prop}
 \begin{proof}
 Fix $T>0$. In Theorem \ref{thm: general Ascoli Arzela}, let $X=[0, T]$ and $Y= B_{L^2(\R)}(0, M_5)$ be the ball around $0$ with radius $M_5>0$  in $L^2(\R)$. Equip $Y$ with the weak topology in $L^2(\R)$ and choose $M_5:= M_2^{1/2} T^{1/2}+ \|\tilde{h}^0\|_{L^2(\R)}$. Note that the topology on $Y$ is metrizable as $L^2(\R)$ is separable, by Theorem 3.29  from \cite{brezis2011functional}. By (\ref{eqn: bound on L^2 norm of h_k(t,.)}) from Proposition \ref{prop: 27}, $h_k(t, \cdot)  \in Y$ for every $t \in X$. By (\ref{eqn: bound on L^2 norm of h_k(t_2,.)-h_k(t_1,.)}) and (\ref{eqn: bound on L^2 norm of h_k(t,.)}), $\{h_k(t, \cdot)\}_{k\in \N}$ is equicontinuous and has compact closure for each $t \in X$. Since $h_k(t, \cdot) \in C^{0,1/2}(X;Y)$ for every $ k \in \N$, by Theorem \ref{thm: general Ascoli Arzela}, we have a subsequence (not relabeled) and $h_* \in C^{0, 1/2}(X;Y)$ such that $h_k (t,\cdot)\rightharpoonup h_\ast(t,\cdot)$ in $L^2(\R)$ for every $t \in [0,T]$.
 \paragraph{} Now, relabel the convergent subsequence for $T=1$ as $\{h_{1,k}\}_{k \in \N}$. For $T=2$, consider the subsequence above to obtain a further subsequence and proceed by induction as above for every $T=N \in \N$. By uniqueness of limits, we obtain $h_\ast \in C^{0, 1/2}([0, \infty); L^2(\R))$ (note that $M_5(T) \rightarrow \infty$ as $T \rightarrow \infty$) and subsequences $\{h_{N,k}\}_{k \in \N}$ such that $h_{N, k}(t,\cdot) \rightharpoonup h_\ast (t,\cdot)$ in $L^2(\R)$. By abuse of notation, define the sequence $h_k:= h_{k,k}$ for every $k \in \N$. Clearly, 
 \begin{equation}\label{eqn: h_k(t,.) weakly conv to h_*(t,.)}
     h_k(t, \cdot) \rightharpoonup h_*(t,\cdot) \text{ in } L^2(\R)
 \end{equation}
 for every $t \in [0, \infty)$. 
 %Now, note that $h_k(t, \cdot) \in H^2(\R)$ for every $t \in [0, \infty)$ by Proposition \ref{prop: H^2 regularity of h_k(t, \cdot)}.Hence, using (\ref{eqn: h_k(t, \cdot) weakly conv to h_*(t,.)}), $h_* (t,.) \in H^2(\R)$ for every $t \in [0, \infty)$. Therefore, by Morrey's theorem, $h_*(t,.) \in C^{1,1/2}(\R)$ for every $t \in [0, \infty)$.
\paragraph{} Since Lip$(h_k(t, \cdot)) \leq L_0$, $\{h_k'(t,\cdot)\}_k \subset L^\infty(\R)$ for every $k \in \N$ and $t \in [0, \infty)$. Hence, for every $\varphi \in C_c^\infty(\R)$ such that $\|\varphi\|_{L^1(\R)} \leq 1$, 
\begin{equation*}
    \bigg|\int_\R h_k'(t,x) \varphi(x) dx\bigg| \leq L_0,
\end{equation*}
for every $t \in [0, \infty)$ and every $k \in \N$. By the Bolzano-Weierstrass theorem, there exists a subsequence of  $\{h_k'(t,\cdot)\}_k$ (not relabeled) that converges weakly star to some $\bar{h}(t,\cdot)$ in $L^\infty(\R)$. But, by (\ref{eqn: h_k(t,.) weakly conv to h_*(t,.)}), we have that $\bar{h}(t, \cdot)=h_\ast'(t, \cdot)$ for a.e. $x \in \R$ and that $h_k'(t,\cdot)\overset{\ast}{\rightharpoonup} h_\ast'(t,\cdot)$, for every $t>0$. Hence, (\ref{eqn: weak star convergence of h_k'(t,.) to h_* (t,.)}) is proved. 
\paragraph{} For $0<T < \infty$, fix $t \in [0, T]$, $k \in \N$ and find $i$ such that $t^{i-1}_k \leq t \leq t^i_k$. By (\ref{eqn: def alpha_k}), (\ref{eqn: def beta_k}) and (\ref{eqn: alpha, beta C^{0, 1/2} bound}), we have 
\begin{gather}
\alpha_k(t) -(M_1 \tau_k)^{1/2} \leq \min \{\alpha^{i-1}_k, \alpha^i_k\}\leq \max \{\alpha^{i-1}_k, \alpha^i_k\} \leq \alpha_k(t) +(M_1 \tau_k)^{1/2}, \label{eqn: bounds on min max of alpha_ks}\\
\beta_k(t) -(M_1 \tau_k)^{1/2} \leq \min \{\beta^{i-1}_k, \beta^i_k\}\leq \max \{\beta^{i-1}_k, \beta^i_k\} \leq \beta_k(t) +(M_1 \tau_k)^{1/2}. \label{eqn: bounds on min max of beta_ks}
\end{gather}
Note that by (\ref{eqn: uniform conv_alpha, beta}) and (\ref{eqn: h_k(t,.) weakly conv to h_*(t,.)}),
\begin{equation}\label{eqn: h_k times characteristic function weak convergence}
    h_k(t, \cdot) \chi_{(\alpha_k(t) -(M_1 \tau_k)^{1/2},  \beta_k(t) +(M_1 \tau_k)^{1/2})} \rightharpoonup h_*(t,\cdot) \chi_{(\alpha(t), \beta(t))} \text{ in } L^2(\R),
\end{equation}
where $\chi_{(\alpha_k(t) -(M_1 \tau_k)^{1/2},  \beta_k(t) +(M_1 \tau_k)^{1/2})}, \chi_{(\alpha(t), \beta(t))} $ are the characteristic functions of the respective intervals in $\R$. Now, by (\ref{eqn: def h_k}), (\ref{eqn: bounds on min max of alpha_ks}) and (\ref{eqn: bounds on min max of beta_ks}), $h_k(t,x)=0$ for $x \in (\alpha_k(t) -(M_1 \tau_k)^{1/2},  \beta_k(t) +(M_1 \tau_k)^{1/2})$. Hence, for any $\varphi \in C_c^\infty(\R)$ such that supp$(\varphi) \subseteq \R \setminus (\alpha (t), \beta(t))$, we can show that 
\begin{equation*}
    \int_\R h_*(t,x) \varphi(t,x) dx =0
\end{equation*}
by (\ref{eqn: uniform conv_alpha, beta}), (\ref{eqn: h_k(t,.) weakly conv to h_*(t,.)}) and (\ref{eqn: h_k times characteristic function weak convergence}). 
Therefore, $h_*(t,x)=0$ for every $x \notin (\alpha (t), \beta(t))$, for every $t \in[0, \infty)$.
\paragraph{} Let $0<\delta<1$ and fix $t \in [0, \infty)$. By (\ref{eqn: uniform conv_alpha, beta}), we know that there exists $N_0 \in \N$ such that 
\begin{equation}\label{eqn: h_k=0 outside alpha-del, beta+del}
    h_k(t,x)=0 \text{ for } x \notin [\alpha(t)- \delta, \beta(t)+ \delta],
\end{equation} for every $k \geq N_0$. Since Lip$(h_k(t, \cdot)) \leq L_0$ for every $k$, we have that $\{h_k(t, \cdot)\}_k$ are uniformly bounded and equicontinuous in $[\alpha(t)-\delta, \beta(t)+\delta]$. Hence, by Ascoli-Arzel\`a theorem, we have 
\begin{equation*}
    h_k(t, \cdot) \rightarrow h_*(t,\cdot) \text{ uniformly in } [\alpha(t) -\delta,  \beta(t) + \delta].
\end{equation*}
Therefore, by (\ref{eqn: h_k=0 outside alpha-del, beta+del}) and the fact that $h_\ast(t,\cdot) =0 $ for $x \notin (\alpha (t), \beta(t))$, we have (\ref{eqn: uniform convergence of h_k in R}). Clearly, (\ref{eqn: uniform convergence of h_k in R}) implies that $\text{Lip}(h_*(t,\cdot)) \leq L_0$ for every $t \in [0, \infty)$.
\paragraph{} Let $h(t,\cdot)$ be the restriction of $h_*(t,\cdot)$ to $[\alpha(t), \beta(t)]$ for every $t \in [0, \infty)$. Then, we have $ h(t, \cdot) \in H^1_0((\alpha(t), \beta(t)))$, and by (\ref{eqn: uniform conv_alpha, beta}) and (\ref{eqn: uniform convergence of h_k in R}) we have 
\begin{eqnarray*}
    A_0 &=& \int_\R h_k(t, \cdot) \chi_{(\alpha_k(t) -(M_1 \tau_k)^{1/2},  \beta_k(t) +(M_1 \tau_k)^{1/2})} dx \\
    &\rightarrow & \int_R h_*(t,\cdot) \chi_{(\alpha(t), \beta(t))} dx =\int_{\alpha(t)}^{\beta(t)} h(t,x) dx \text{ as } k \rightarrow \infty,
\end{eqnarray*}
thus proving (\ref{eqn: integral of h(t,.)=A_0}). 
\paragraph{} By (\ref{eqn: bounds on min max of alpha_ks}) and (\ref{eqn: bounds on min max of beta_ks}), for every $x \in (\alpha_k(t) + (M_1 \tau_k)^{1/2}, \beta_k(t) - (M_1 \tau_k)^{1/2})$,
\begin{equation*}
    h_k(t,x) = h_k^{i-1}(x)+ \frac{t-t_k^{i-1}}{\tau_k}(h_k^{i}(x)-h_k^{i-1}(x)).
\end{equation*}
Fix $a, b$ such that $\alpha(t) < a< b< \beta(t)$. Then, by (\ref{eqn: uniform conv_alpha, beta}), there exists $N_1 \in \N$ such that $\alpha_k(t) + (M_1 \tau_k)^{1/2} <a<b< \beta_k(t) - (M_1 \tau_k)^{1/2}$ for $k\geq N_1$. Therefore, by (\ref{eqn:uni L^2 bound h^i_k''}), $h_k''(t,\cdot) \in L^2((a,b))$ for every $k \geq N_1$ and  
\begin{eqnarray*}
    \int_a^b |h_k''(t,x)|^2 dx &\leq& \frac{t-t^{i-1}_k}{\tau_k} \int_a^b |(h^{i}_k)''(t,x)|^2 dx + \frac{t^i_k-t}{\tau_k}\int_a^b |(h^{i-1}_k)''(t,x)|^2 dx\\
    &\leq& M_3.
\end{eqnarray*}
Since $\{h_k''(t,\cdot)\}_k$ is uniformly bounded in $L^2((a,b))$, using (\ref{eqn: uniform convergence of h_k in R}), there exists a subsequence (not relabeled) such that 
\begin{equation*}
    h_k''(t,\cdot) \rightharpoonup h''(t,\cdot) \text{ in } L^2((a,b)),
\end{equation*}
and by the lower semi-continuity of norms in weak convergence,
\begin{eqnarray*}
    \int_a^b |h''(t,x)|^2 dx \leq \liminf_{k \rightarrow \infty } \int_a^b |h_k''(t,x)|^2 dx \leq M_3.
\end{eqnarray*}
Taking the limit as $a \rightarrow \alpha(t)$ and $b \rightarrow \beta(t)$, we obtain that $h(t,\cdot) \in H^2((\alpha(t), \beta(t)))$ and prove (\ref{eqn: L^2 bound on h''}). In turn, using the fundamental theorem of calculus and H\"older's inequality, we have (\ref{eqn: C^{0,1/2} seminorm bound on h'(t,.) in (alpha(t), beta(t))}).

 \end{proof}

 \begin{prop}\label{prop: 32}
    Let $\{h_k\}_k$ and $h^\ast$ be the subsequence and the function given by Proposition \ref{prop:31}, respectively. Then, there exists a further subsequence (not relabeled) such that for every $T>0$,
    \begin{equation}\label{eqn: time derivative weak conv h_k}
         h_k \rightharpoonup h_* \text{ weakly in } H^1((0,T); L^2(\R)).
     \end{equation}
\end{prop}
\begin{proof}
    We produce a uniform $T$ dependent upper bound in $H^1((0,T);L^2(\R))$  and then use the reflexivity of the space to obtain a convergent subsequence.
    \paragraph{} By (\ref{eqn: bound on time derivative of h_k}) and  (\ref{eqn: bound on L^2 norm of h_k(t,.)}), 
\begin{gather}
    \int_0^T \int_\R |h_k(t,x)|^2 dx dt+ \int_0^T \int_\R |\dot{h}_k(t,x)|^2 dx dt  \leq 2\bigg(\int_0^T ( M_2 t+\|\tilde{h}^0\|_{L^2(\R)}^2)dt\bigg)+ M_2 \nonumber \\
    = M_2 T^2 + 2 \|\tilde{h}^0\|_{L^2(\R)}^2 T + M_2
    =:\overline{M}_2. \label{eqn: h_k H^1 bound}
\end{gather}
Hence, by the Kakutani's theorem, there exists a subsequence (not relabeled) $\{ h_k\}_k$ such that $h_k \rightharpoonup \bar{h}$ for some $\bar{h} \in H^1((0,T);L^2(\R))$. 
\paragraph{} Now we need to show that $\bar{h}= h_\ast $ a.e. in $(0,T) \times \R$. Let $\varphi \in C_c^\infty ((0,T))$ and $\psi \in C_c^\infty (\R)$. Since $ h_k \rightharpoonup \bar{h} \text{ weakly in } H^1((0,T);L^2(\R))$, using Lebesgue's dominated convergence theorem,
\begin{gather*}
    \int_0^T \int_\R \varphi(t) \psi(x) \bar{h}(t,x) dx dt = \lim_{k \rightarrow \infty} \int_0^T \int_\R \varphi(t) \psi(x) h_k(t,x) dx dt\\ 
    =  \lim_{k \rightarrow \infty} \int_0^T \varphi(t) \bigg(\int_\R  \psi(x) h_k(t,x) dx\bigg) dt
    = \int_0^T \varphi(t)  \lim_{k \rightarrow \infty}\bigg(\int_\R  \psi(x) h_k(t,x) dx\bigg) dt,
\end{gather*}
where we used the bound
\begin{gather*}
    \bigg| \int_\R  \psi(x) h_k(t,x) dx\bigg|  \leq  \|\psi\|_{L^2(\R)} \| h_k(t, \cdot)\|_{L^2(\R)} \leq\|\psi\|_{L^2(\R)} (M_2 t^{1/2} + \|\tilde{h}_0\|_{L^2(\R)})
\end{gather*}
which follows from (\ref{eqn: bound on L^2 norm of h_k(t,.)}) and H\"older's inequality. Now by (\ref{eqn: upper bound on |beta_k-alpha_k|}) supp$(h_k(t, \cdot))$ is contained in the compact interval $[\alpha_0 - (tM_1)^{1/2}, \beta_0 + (tM_1)^{1/2}]$. Hence,
\begin{eqnarray*}
    \int_0^T \int_\R \varphi(t) \psi(x) \bar{h}(t,x) dx dt &=&  \int_0^T \varphi(t)  \lim_{k \rightarrow \infty}\bigg(\int_{\alpha_0 - (tM_1)^{1/2}}^{\beta_0 +(tM_1)^{1/2}}  \psi(x) h_k(t,x)dx\bigg) dt\\
    &=& \int_0^T \int_\R \varphi(t) \psi(x) h_\ast(t,x) dx dt
\end{eqnarray*}
since $ h_k (t,\cdot) \rightarrow h_{*}(t,\cdot) \text{ uniformly in } \R$ by (\ref{eqn: uniform convergence of h_k in R}). This implies that $h_\ast \in H^1((0,T);L^2(\R))$ and $h_k \rightharpoonup h_\ast \text{ weakly in } H^1((0,T);L^2(\R)) $.

\end{proof}

\begin{prop}\label{prop 33 & corollary 34}
   Let $\{h_k\}_k$ be the subsequence given by Proposition \ref{prop: 32} and $\hat{h}_k$ be defined as in (\ref{eqn: piecewise constant extension of alpha, beta, h}) for every $k$. Then
   \begin{equation} \label{eqn: h_k hat weak L^2 convergence}
       \hat{h}_k(t,\cdot) \rightharpoonup h_*(t,\cdot) \text{ in } L^2(\R) \text{ for every } t \geq 0. 
   \end{equation}
Moreover, there exists a further subsequence (not relabeled) such that, for fixed $t \geq 0$, 
\begin{gather}
   % \hat{h}_k(t,.) &\rightharpoonup& h_*(t,.) \text{ in } H^2(\R),\label{eqn: h_k hat weak H^2 convergence}\\
     \hat{h}_k(t,\cdot) \rightarrow h(t,\cdot) \text{ uniformly in } \R \label{eqn: uniform conv of hat h_k in R}, \text{ and}\\
     \hat{h}_k'(t,\cdot) \overset{\ast}{\rightharpoonup} h'(t,\cdot) \text{ in } L^\infty(\R).\label{eqn: h_k hat ' weak start L^infinity convergence}
\end{gather}
\end{prop}
\begin{proof}
By (\ref{eqn: def h_k}), (\ref{eqn: piecewise constant extension of alpha, beta, h}) and (\ref{eqn: bound on L^2 norm of h_k(t_2,.)-h_k(t_1,.)}), for every $t \in (t^{i-1}_k, t^i_k]$, 
\begin{eqnarray}
    \|\hat{h}_k(t,\cdot)-h_k(t, \cdot)\|_{L^2(\R)} = \|h_k(t^i_k,.)-h_k(t, \cdot)\|_{L^2(\R)} 
    \leq M_2^{1/2}\tau_k^{1/2} \rightarrow 0 \text{ as } k \rightarrow \infty. \label{eqn: h_k and h_k hat are L^2 close}
\end{eqnarray}
Let $\varphi \in C_c^\infty (\R)$ and let $t \in (0,T)$. Then, for every $k\in \N$, there exists $0\leq i\leq k$ such that $t \in (t^{i-1}_k, t^i_k]$, and
\begin{gather*}
    \bigg| \int_\R (\hat{h}_k(t,x)- h_* (t,x)) \varphi(x) dx \bigg| \leq \bigg| \int_\R (\hat{h}_k(t,x)- h_k(t,x)) \varphi(x) dx \bigg|\\
    +\bigg| \int_\R (h_k(t,x)- h_* (t,x)) \varphi(x) dx \bigg|
    \leq \|\varphi\|_{L^2(\R)}\|\hat{h}_k(t,\cdot)-h_k(t, \cdot)\|_{L^2(\R)}\\
    + \bigg| \int_\R (h_k(t,x)- h_* (t,x)) \varphi(x) dx \bigg| \rightarrow 0 \text{ as } k \rightarrow \infty,
\end{gather*}
by (\ref{eqn: time derivative weak conv h_k}) and (\ref{eqn: h_k and h_k hat are L^2 close}). Hence, (\ref{eqn: h_k hat weak L^2 convergence}) follows. Since $\text{Lip}(\hat{h}_k(t, \cdot) \leq L_0$ for every $t\geq 0$ and $k \in \N$, by the Ascoli-Arzel\`a theorem, Bolzano-Weierstrass theorem and (\ref{eqn: h_k hat weak L^2 convergence}), we can follow the steps in the proof of Proposition \ref{prop:31} to obtain a subsequence such that (\ref{eqn: uniform conv of hat h_k in R}) and (\ref{eqn: h_k hat ' weak start L^infinity convergence}) hold.

%\paragraph{} Further, as $\hat{h}_k'(t,x) \leq L_0$ for every $t, x$, using (\ref{eqn: upper bound on |beta_k-alpha_k|}) and (\ref{eqn: uniform conv_alpha, beta}), we have that $\{\hat{h}_k'\}_k$ is bounded in $L^2(\R)$ for a fixed $t>0$. Hence, $\{\hat{h}_k\}_k$ is bounded in $H^2(\R)$ by (\ref{eqn:uni L^2 bound h^i_k''}) and the results above.

\end{proof}
From here on, we assume that the sequences $\{\alpha_k\}_k, \{\beta_k\}_k,\{\hat{\alpha}_k\}_k, \{\hat{\beta}_k\}_k, \{h_k\}_k, \{\hat{h}_k\}_k$ are the subsequences mentioned in (\ref{eqn: alpha, beta weak conv})-(\ref{eqn: uniform conv_alpha hat, beta hat}), (\ref{eqn: uniform convergence of h_k in R}), (\ref{eqn: weak star convergence of h_k'(t,.) to h_* (t,.)}), (\ref{eqn: time derivative weak conv h_k}) and (\ref{eqn: h_k hat weak L^2 convergence})-(\ref{eqn: h_k hat ' weak start L^infinity convergence}).
\paragraph{}The proof of the following lemma is adapted from the proof of \cite[Lemma 9.7]{dal2025motion}.
\begin{lem}\label{lem: 35}
    Let $\{t_k\}_k$ be a sequence of non-negative numbers converging to some $t_0\geq 0$. Then, there exists a subsequence of $\{h_k\}_{k \in \N}$ (not relabeled) such that
     \begin{equation}\label{eqn:uniform convergence of h_k hat in R t_k t_0}
        h_k(t_k, \cdot) \rightarrow h_*(t_0, \cdot), \,\hat{h}_k(t_k, \cdot) \rightarrow h_*(t_0, \cdot) \text{ uniformly on } \R.
    \end{equation}
Moreover, if $\alpha(t_0)<a<b< \beta(t_0)$, then
 \begin{equation}\label{eqn: uniform convergence of h_k hat'}
     h_k'(t,\cdot) \rightarrow h'(t,\cdot), \,    \hat{h}_k'(t,\cdot) \rightarrow h'(t,\cdot) \text{ uniformly on }[a,b].
   \end{equation}
   Finally, if $\{x_k\}_k$ is a sequence in $\R$ converging to some $x_0 \in \R$ such that $\hat{\alpha}_k(t_k) \leq x_k \leq \hat{\beta}_k(t_k)$, then $\alpha(t_0) \leq x_0 \leq \beta(t_0)$ and 
   \begin{equation}
       \hat{h}_k'(t_k, x_k) \rightarrow h'(t_0, x_0).
   \end{equation}
\end{lem}
\begin{proof}
    We have that Lip $h_k(t, \cdot) \leq L_0$ for every $k \in \N$. By (\ref{eqn: bound on L^2 norm of h_k(t_2,.)-h_k(t_1,.)}), we can apply the Ascoli-Arzel\`a theorem as in the proof of (\ref{eqn: uniform convergence of h_k in R}) in Proposition \ref{prop:31}, to obtain a subsequence (not relabeled) such that 
    \begin{equation*}
        h_k(t, \cdot) \rightarrow g(\cdot) \text{ uniformly in } \R
    \end{equation*}
    for some Lipschitz continuous function $g: \R \rightarrow [0,\infty)$. Now, by (\ref{eqn: bound on L^2 norm of h_k(t_2,.)-h_k(t_1,.)}) and (\ref{eqn: time derivative weak conv h_k}), we have 
    \begin{gather*}
        \|h_k(t_k,\cdot)-h_k(t_0,\cdot)\|_{L^2(R)} \leq M_2^{1/2}|t_k-t_0|^{1/2}, \text{ and},\\
        h_k(t_0, \cdot) \rightharpoonup h_*(t_0, \cdot) \text{ in } L^2(\R).
    \end{gather*}
Therefore, 
\begin{equation*}
     h_k(t_k, \cdot) \rightharpoonup h_*(t_0, \cdot) \text{ in } L^2(\R),
\end{equation*}
and hence $g(\cdot)= h_*(t_0, \cdot)$ for every $x \in \R$. This proves the first equation in (\ref{eqn:uniform convergence of h_k hat in R t_k t_0}). The second equation in (\ref{eqn:uniform convergence of h_k hat in R t_k t_0}) can be proved similarly, using (\ref{eqn: bound on L_2 norm of h_k hat(t_1)-h_k hat(t_2)}), (\ref{eqn: h_k hat weak L^2 convergence}) and (\ref{eqn: uniform conv of hat h_k in R}). The rest of the proof follows as in the proof of Lemma 9.7 in \cite{dal2025motion}.

\end{proof}

The proof of the following theorem is adapted from \cite[Theorem 3.2]{piovano2014evolution}.
 \begin{thm}\label{Thm: Piovano 3.2.2}
       Let $a,b \in \R$, $T>0$ and $I \subseteq [0,T]$ such that $\alpha(t) < a< b< \beta(t)$ for every $t \in I$. Then, $h \in C^{0, s_2}([0,T]; C^{1, s_1}([a,b]))$ and there exists a subsequence of $\{h_k\}_{k \in \N}$ (not relabeled) such that
    \begin{gather}
        h_k \rightarrow h \text{ in } C^{0, s_2}(I; C^{1, s_1}([a,b])) \text{ as } k \rightarrow \infty \label{eqn: strong conv_C^1, alpha h_k},\\
        \hat{h}_k \rightarrow h \text{ in } L^\infty(I; C^{1, s_1}([a,b])) \text{ as } k \rightarrow \infty, \label{eqn: Piovano_h_k hat L^infty C^1, alpha convergence}
    \end{gather}
for $s_1 \in (0, \frac{1}{2})$ and $ s_2 \in (0, \frac{1-2s_1}{8})$. Moreover, for every fixed $k \in \N$, $h_k(t, \cdot) \rightarrow h_0(\cdot)$ in $C^{1, s_1}([a,b])$ as $t \rightarrow 0^+$.
    \end{thm}

\begin{proof}
    Let $t_1, t_2 \in I$. By (\ref{eqn: uniform conv_alpha, beta}), we can choose some $k_* \in \N$ such that $\alpha_k(t) < a<b< \beta_k(t)$ for every $k \geq k_*$. Fix some $k\geq k_*$ and define $g(x):= h_k(t_2,\cdot)-h_k(t_1,\cdot)$ for $x \in (a,b)$.
    \paragraph{} Recall the definition of $\bar{\alpha}, \bar{\beta}$ from (\ref{eqn: def alpha_k bar, beta_k bar}). By Proposition \ref{prop: 30}, 
  \begin{equation}\label{eqn: alpha bar, beta bar convergence}
      \bar{\alpha}_k(t) \rightarrow \alpha(t), \quad \bar{\beta}_k(t) \rightarrow \beta(t)
  \end{equation}
  as $k \rightarrow \infty$ for every $t>0$. Hence, there exists $N_0 \in \N$ such that $\bar{\alpha}_k(t)<a<b<\bar{\beta}_k(t)$ for every $k \geq N_0$. Therefore, by (\ref{eqn: upper bound on H^2 norm of h_k(t,.)}), 
  \begin{equation}\label{eqn: h_k(t,.) H^2 bound}
     \|h_k(t, \cdot)\|_{H^2((a,b))} \leq M_4(1+T^{1/4}+T^{1/2}), 
  \end{equation}
  for every $t \in [0, T]$. This, in turn, implies that $g \in H^2((a,b))$. Since $b-a$ is bounded above by a constant independent of $a, b, I$ by (\ref{eqn: upper bound on |beta_k-alpha_k|}), we can use Theorem \ref{Thm: 7.41, AFCSS} to obtain $C>0$ such that
    \begin{equation}
        \|g'\|_{L^\infty}((a,b)) \leq C(\|g\|_{L^2((a,b))}+ \|g\|_{L^2((a,b))}^{1/4}\|g''\|_{L^2((a,b))}^{3/4}),
    \end{equation}
where $C$ depends only on $T$. Therefore, by (\ref{eqn: bound on L^2 norm of h_k(t_2,.)-h_k(t_1,.)}) and (\ref{eqn:uni L^2 bound h^i_k''}),
\begin{gather}
    \|h_k'(t_2,\cdot)-h_k'(t_1,\cdot)\|_{L^\infty((a,b))} \leq C (\|h_k(t_2,\cdot)-h_k(t_1,\cdot)\|_{L^2((a,b))}\nonumber\\
    +\|h_k(t_2,\cdot)-h_k(t_1,\cdot)\|^{1/4}_{L^2((a,b))}\|h_k''(t_2,\cdot)-h_k''(t_1,\cdot)\|^{3/4}_{L^2((a,b))})\nonumber\\
    \leq  CM_2^{1/2}|t_2-t_1|^{1/2}+ 2C M_3^{3/8}M_2^{1/8}|t_2-t_1|^{1/8} \leq C_0 (|t_2-t_1|^{1/8}+|t_2-t_1|^{1/2}),\label{eqn: Piovano_h_k' difference L^infty bound}
\end{gather}
  for some $C_0>0$. Now, by the mean value theorem, there exists $x_0 \in (a,b)$ such that 
  \begin{equation*}
      g(x_0)= \frac{1}{b-a}\int_a^b g(x) dx.
  \end{equation*}
  Therefore, for every $x \in [a,b]$, 
  \begin{eqnarray*}
      |g(x)|=|g(x)- g(x_0)|+|g(x_0)|\leq (b-a)\|g'\|_{L^\infty((a,b))} + (b-a)^{-1/2}\|g\|_{L^2((a,b))}.
  \end{eqnarray*}
 That is, by (\ref{eqn: bound on L^2 norm of h_k(t_2,.)-h_k(t_1,.)}) and (\ref{eqn: Piovano_h_k' difference L^infty bound}),  
 \begin{equation} \label{eqn: Piovano_h_k difference L^infty bound}
     \|h_k(t_2,\cdot)-h_k(t_1,\cdot)\|_{L^\infty((a,b))} \leq C(|t_2-t_1|^{1/8}+ |t_2-t_1|^{1/2}).
 \end{equation}
  Note that $C$ depends on $a, b$ in the equation above. Now, by (\ref{eqn: h_k(t,.) H^2 bound}) and Morrey's theorem, we have 
  \begin{equation}\label{eqn: g' 1/2 seminorm bound}
      |h_k'(t_2,\cdot)-h_k'(t_1,\cdot)|_{1/2} \leq C(1+T^{1/4}+T^{1/2}),
  \end{equation}
  where $|.|_{s}$ denotes the $s-$H\"older seminorm for $0<s<1$. If $0< s_1< \frac{1}{2}$, 
  \begin{eqnarray*}
      |g'|_{s_1} &=& \sup \bigg\{ \frac{|g'(x)-g'(y)|}{|x-y|^{s_1}}: x,y \in [a,b], x \neq y \bigg\}\\
      &\leq& |g'|_{1/2}^{2s_1}(2 \|g'\|_{L^\infty((a,b))})^{1-2s_1}.
  \end{eqnarray*}
Therefore, by (\ref{eqn: Piovano_h_k' difference L^infty bound}) and (\ref{eqn: g' 1/2 seminorm bound}), 
\begin{equation}\label{eqn: Piovano_h_k difference Holder seminorm bound}
    |h_k'(t_2,x)-h_k'(t_1,x)|_{s_1} \leq C (|t_2-t_1|^{\frac{1-2 s_1}{2}}+|t_2-t_1|^{\frac{1-2 s_1}{8}}),
\end{equation}
where $C$ depends on $T_1$. Hence, by (\ref{eqn: Piovano_h_k' difference L^infty bound}), (\ref{eqn: Piovano_h_k difference L^infty bound}), and (\ref{eqn: Piovano_h_k difference Holder seminorm bound}), we have 
\begin{equation}\label{eqn: h_k difference t difference bound}
    \|h_k(t_2,\cdot)-h_k(t_1,\cdot)\|_{C^{1, \alpha}([a,b])} \leq C |t_2-t_1|^{\frac{1-2s_1}{8}},
\end{equation}
where we assumed without loss of generality that $|t_2-t_1| \leq 1$ and used the fact that $\frac{1-2s_1}{8}<\frac{1-2s_1}{2} $ and $\frac{1-2s_1}{8} < \frac{1}{8} < \frac{1}{4}< \frac{1}{2}$.
Therefore, using the Ascoli-Arzel\`a theorem, there exists a subsequence (not relabeled) such that $h_k \rightarrow h$ in $C^{0, s_2}(I; C^{1,s_1}([a,b]))$ for $0< s_2 < \frac{1-2s_1}{8}$ and $0<s_1< \frac{1}{2}$. Since $h(0,\cdot)=h_0(\cdot)$ in $(\alpha(t), \beta(t))$ by (\ref{eqn: def h_k}) and (\ref{eqn: uniform convergence of h_k in R}), this implies that for every fixed $k \in \N$, $h_k(t, \cdot) \rightarrow h_0(\cdot)$ in $C^{1, s_1}([a,b])$ as $t \rightarrow 0^+$. 
\paragraph{} Now, for every $t \in I$, and $k \in \N$, $t \in (t^{i_k-1}_k, t^{i_k}_k]$ for some $0<i_k \leq k$. Then, since $t^{i_k}_k \rightarrow t$ as $k \rightarrow \infty$,
\begin{gather*}
    \sup_{t \in [0,T]}\|\hat{h}_k(t,\cdot)-h(t,\cdot)\|_{C^{1, \alpha}([a,b])} = \sup_{t \in [0,T]} \|h_k(t^i_k,\cdot)-h(t,\cdot)\|_{C^{1, \alpha}([a,b])} \\
    \leq \sup_{t \in [0,T]}\|h_k(t^i_k,\cdot)-h_k(t,\cdot)\|_{C^{1, \alpha}([a,b])} + \sup_{t \in [0,T]} \|h_k(t,\cdot)-h(t,\cdot)\|_{C^{1, \alpha}([a,b])} \\
    \leq C\tau_k^{\frac{1-2s_1}{8}}+\sup_{t \in [0,T]} \|h_k(t,\cdot)-h(t,\cdot)\|_{C^{1, \alpha}([a,b])} \rightarrow 0
\end{gather*}
as $ k \rightarrow \infty$, by (\ref{eqn: piecewise constant extension of alpha, beta, h}), (\ref{eqn: strong conv_C^1, alpha h_k}) and (\ref{eqn: h_k difference t difference bound}). Hence, (\ref{eqn: Piovano_h_k hat L^infty C^1, alpha convergence}) follows.

\end{proof}
\begin{prop}
    For every fixed $t\in [0, \infty)$ and $\alpha(t)<a<b<\beta(t)$, $h(t,\cdot) \in H^2((a,b))$.
\end{prop}
\begin{proof}
    Fix $t \in [0, \infty)$ and $\alpha(t)<a<b< \beta(t)$. Recall from the proof of Theorem \ref{Thm: Piovano 3.2.2} that there exists $N_0 \in \N$ such that $\bar{\alpha}_k(t)<a<b<\bar{\beta}_k(t)$ for every $k \geq N_0$. Then, by (\ref{eqn: upper bound on H^2 norm of h_k(t,.)}) from Proposition \ref{prop: H^2 regularity of h_k(t,.)}, $\{h_k(t, \cdot)\}_{k \geq N_0}$ is uniformly bounded in $H^2((a,b))$. Therefore, there exists a subsequence, not relabeled, such that 
    \begin{equation*}
        h_k(t,\cdot) \rightharpoonup h(t,\cdot) \text{ in } H^2((a,b)),
    \end{equation*}
    where we used (\ref{eqn: uniform convergence of h_k in R}) and Kakutani's theorem. 
    
\end{proof}
Now, recall the assumptions (\ref{assumption: Lip h< L_0 and h>0}) from Theorem \ref{thm: 4_EL Theorem} and (\ref{eqn: h'(alpha), h'(beta) bounds}), (\ref{eqn: Lower bound on h_EL section-assumption}) from Theorem \ref{thm: 6}. The following result asserts that these assumptions are valid when $t$ is restricted to a finite interval. The proof follows directly from Theorem 9.11, Proposition 9.12, Proposition 9.13, Theorem 9.14 and Proposition 9.15 in \cite{dal2025motion}. The derivatives of $h, h^i_k$ at the endpoints of the intervals are one-sided.
\begin{thm}\label{Res: 39,40,41,42,43}
        Under the assumptions (\ref{assumption: alpha_0, beta_0, h_0 are in A_s})-(\ref{assumption: Lip h_0< L_0}),
        \begin{itemize}
            \item[(i)]there exists $T_0>0$ such that for all $t \in [0, T_0]$, 
        \begin{gather}
            h(t,x)>0 \text{ for all } x \in (\alpha(t), \beta(t)),\label{eqn: h positive in (alpha, beta)}\\
            h'(t, \alpha(t))>0 \text{ and }  h'(t, \beta(t))<0. \label{eqn: h' sign at alpha, beta}
        \end{gather}
        
        \item[(ii)]There exist $k_0 \in \N$ and $0< \eta_0<1$ such that 
    \begin{equation}
        (h^i_k)'(\alpha^i_k) > 2 \eta_0 \text{ and }  (h^i_k)'(\beta^i_k) < -2 \eta_0 
    \end{equation}
    for all $k \geq k_0$ and all $0 \leq i \leq kT_0$.
    \item[(iii)] Let $\delta \in \R$ such that
         \begin{equation}\label{eqn: choose delta prop 41}
        0< \delta < \frac{1}{2} \min_{t \in [0, T_0]} (\beta(t)-\alpha(t)).
    \end{equation}
    Then, there exist $k_1 \geq k_0$ and $0< \eta_1< 1$ such that 
    \begin{gather}
        h^i_k(x) \geq 2 \eta_1 \text{ for all } x \in [ \alpha^i_k +\delta, \beta^i_k - \delta] \text{ and } \\
        h^i_k(x) > 0 \text{ for all } x \in (\alpha^i_k, \beta^i_k)
    \end{gather}
    for all $k \geq k_1$ and $0 \leq i \leq kT_0$.
    \item[(iv)]There exists $0< T_1 \leq T_0$ such that 
        \begin{equation}\label{eqn: Lip constant strict bound on h_*}
            \operatorname{Lip } h_*(t,\cdot) < L_0
        \end{equation}
        for every $t \in [0, T_1]$.
        \item[(v)]There exists $k_2 \geq k_1$ such that 
        \begin{equation}
           \operatorname{Lip } \hat{h}_k(t,\cdot) < L_0
        \end{equation}
for all $k \geq k_2$ and all $t \in [0, T_1]$.
      \end{itemize}  
    \end{thm}
The next result follows from Theorem 9.17 in \cite{dal2025motion}. The proof is dependent on Theorem \ref{thm: 24} and Theorem \ref{Res: 39,40,41,42,43}. 
\begin{thm}\label{Thm: 45}
        Under the assumptions  (\ref{assumption: alpha_0, beta_0, h_0 are in A_s})-(\ref{assumption: Lip h_0< L_0}), let $L_0$, $T_1$ and $k_2$ be as in Theorem \ref{Res: 39,40,41,42,43}. Let $p_1, q_1>1$ be defined as 
        \begin{equation}
            p_1:=\min\{6/5, p_0/(4-2p_0)\} \quad \text{and} \quad q_1:=4p_1/(2+p_1),
         \end{equation}
        where $p_0$ is as in Theorem \ref{thm: 13, elliptic regularity in triangle}. Then, there exists $M_6>0$ such that 
        \begin{eqnarray}\label{eqn: h_k^''' integral bound}
        \int_0^{T_1} \|\hat{h}_k'''(t,\cdot)\|_{L^{q_1}((\hat{\alpha}_k(t), \hat{\beta}_k(t)))}^{p_1} dt \leq M_6 \bigg( \int_0^{T_1} (\hat{\beta}_k(t)-\hat{\alpha}_k(t))^{\frac{2-5p_1}{4-2p_1}} dt\bigg)^{1-\frac{p_1}{2}} + M_6,
    \end{eqnarray}
    \begin{equation}\label{eqn: h_k ^ (iv) integral bound}
        \int_0^{T_1} \int_{\hat{\alpha}_k(t)}^{\hat{\beta}_k(t)}|\hat{h}_k^{(iv)}(t,x)|^{p_1} dx dt \leq M_6,
    \end{equation}
    for all $k \geq k_2$.
    \end{thm}

\begin{cor} \label{Prop: Regularity of h}
    Let $T_1$ be as in Theorem \ref{Res: 39,40,41,42,43}. Then, $h^{(iv)} \in L^{p_1}([0,T_1] \times [\alpha(t), \beta(t)])$. In particular, $h^{(iv)} \in L^{1}([0,T_1] \times [\alpha(t), \beta(t)])$ and there exists $\overline{M}_6>0$ such that 
    \begin{equation}\label{eqn: L^1 bound on h^(iv)}
         \int_0^{T_1} \int_{\alpha(t)}^{\beta(t)}|h^{(iv)}(t,x)|dx dt \leq \overline{M}_6.
    \end{equation}
\end{cor}

\begin{proof}
    By (\ref{eqn: h_k ^ (iv) integral bound}) from Theorem \ref{Thm: 45}, we know that
    \begin{equation*}
         \int_0^{T_1} \int_\R |\hat{h}_k^{(iv)}(t,x)\chi_{[\alpha_k(t), \beta_k(t)]}(x)|^{p_1} dx dt \leq M_6,
    \end{equation*}
    for every $k \geq k_2$, where $\chi_{[\alpha_k(t), \beta_k(t)]}$ is the characteristic function of the interval $[\alpha_k(t), \beta_k(t)]$ and $k_2$ is defined in Theorem \ref{Res: 39,40,41,42,43}. Since $p_1>1$, there exists a further subsequence $\{\hat{h}_{k_j}^{(iv)}\chi_{[\alpha_{k_j}(t), \beta_{k_j}(t)]}\}_{k_j}$ that converges weakly to some $g \in  L^{p_1}([0,T_1] \times \R)$. But, as $\hat{h}_k(t,\cdot) \rightarrow h_*(t,\cdot)$ uniformly in $\R$ for every $t \geq 0$ by (\ref{eqn: uniform conv of hat h_k in R}) and $\chi_{[\alpha_k(t), \beta_k(t)]} \rightarrow \chi_{[\alpha(t), \beta(t)]}$ uniformly in $\R$ by (\ref{eqn: uniform conv_alpha, beta}), we know that $g(t, x)= h^{(iv)}(t,x) \chi_{[\alpha(t), \beta(t)]}(x)$ a.e. in $[0,T_1] \times \R$. Hence, $h^{(iv)}\chi_{[\alpha(t), \beta(t)]} \in L^{p_1}([0,T_1] \times \R)$. Also note that we have the bound
    \begin{gather}
        \int_0^{T_1} \int_\R  |h^{(iv)}(t,x)\chi_{[\alpha(t), \beta(t)]}(x)|^{p_1} dx dt = \int_0^{T_1} \int_{\alpha(t)}^{\beta(t)}|h^{(iv)}(t,x)|^{p_1}dx dt \nonumber \\
        \leq \liminf_{j \rightarrow \infty} \int_0^{T_1} \int_{\hat{\alpha}_k(t)}^{\hat{\beta}_k(t)} |\hat{h}_{k_j}^{(iv)}(t,x)|^{p_1} dx dt \leq M_6. \label{eqn: h^iv L^p_1 bound}
    \end{gather}
   By (\ref{eqn: upper bound on |beta_k-alpha_k|}) and (\ref{eqn: uniform conv_alpha, beta}), we have that 
   \begin{equation}
       \beta(t)-\alpha(t) \leq \beta_0-\alpha_0+2(T_1M_1)^{1/2}
   \end{equation}
   for every $t \in [0, T_1]$. Hence, as $p_1>1$, by H\"older's inequality, we have that 
    \begin{eqnarray*}
          \int_0^{T_1} \int_{\alpha(t)}^{\beta(t)}|h^{(iv)}(t,x)|dx dt &\leq& (T_1(\beta_0-\alpha_0+2(T_1M_1)^{1/2}))^{1/p_1'}\bigg(\int_0^{T_1} \int_{\alpha(t)}^{\beta(t)} |h^{(iv)}(t,x)|^{p_1} dx dt\bigg)^{1/p_1}\\
         &\leq& (T_1(\beta_0-\alpha_0+2(T_1M_1)^{1/2}))^{1/p_1'}M_6^{1/p_1}=: \overline{M}_6,
    \end{eqnarray*}
    where $p_1'>1$ is the conjugate of $p_1$.
    
\end{proof}
Now, we study the convergence of the elastic energy term given in (\ref{eqn: def elastic energy}). Let $T_1$ be as in Theorem \ref{Res: 39,40,41,42,43} (iv). For every $t \in [0, T_1]$, let $u(t,\cdot,\cdot)$ be the unique minimizer of the problem 
\begin{equation}\label{problem: minimization problem on h(t,.)}
        \text{min}\bigg\{\int_{\Omega_{h(t,\cdot)}}W(Ev(x,y))dxdy:v \in \a_e(\alpha(t),\beta(t), h(t,\cdot))\bigg\}.
    \end{equation}
Note that the minimizer is unique due to the strict convexity of the functional. We state the following proposition without proof as it follows directly from \cite[Proposition 9.16]{dal2025motion}.
\begin{prop}\label{prop: 44}
    Let $T_1$ be as in Theorem \ref{Res: 39,40,41,42,43}(iv) and $\{t_k\}_k$ be a sequence in $[0, T_1]$ converging to $t_0 \in [0, T_1]$. Assume that for every $k$ there exists $i_k \in \N \cup \{0\}$ such that $t_k=t^{i_k}_k$. Then, $\{u^{i_k}_k\}_{k \in \N}$ converges to $u(t_0, \cdot,\cdot)$ weakly in $H^1(\tilde{\Omega}; \R^2)$ for every open set $\tilde{\Omega} \subset \Omega_{h(t_0, \cdot)}$ with ${\operatorname*{dist}}(\tilde{\Omega}, {\operatorname*{graph}}(h(t_0, \cdot)))>0$, where $u^i_k$ for $1 \leq i \leq k$, $i,k \in \N$ is defined in section \ref{sec: Discretization}.
\end{prop}
Now, we prove the following proposition.
\begin{prop} \label{prop: subdivide [0,T_1]}
    Let $T_1$ be as in Theorem \ref{Res: 39,40,41,42,43} (iv). Then, there exists $n_0 \in \N$ such that $[0,T_1]$ can be divided into sub intervals $I_1,I_2,..., I_{n_0}$ such that $\cup_{k=1}^{n_0} I_k = [0,T_1], I_k\cap I_{k'} = \phi \text{ when } k \neq k', 1\leq k, k' \leq n_0,$ and for every $1 \leq k \leq n_0$, there exist $a,b \in \R$ such that 
    \begin{equation}\label{eqn: a,b between alpha beta functions}
        \alpha(t) <a<b< \beta(t) \text{ for every } t \in I_k.
    \end{equation}
\end{prop}
\begin{proof}
Let $A_1:=\sqrt{\frac{2A_0}{L_0}}$. Then, by (\ref{eqn: alpha(0), beta(0), beta-alpha bdd below}), $\beta(t)- \alpha(t) \geq A_1$ for every $t \in [0, \infty)$.
    By (\ref{eqn: alpha, beta weak conv}), we know that $\alpha(t), \beta(t) \in H^1((0, T_1))$.  In turn, by Morrey's theorem, $\alpha, \beta : [0,T_1] \rightarrow \R$ are uniformly continuous. Now, choose $\delta>0$ such that 
    \begin{eqnarray*}
        |\alpha(t)-\alpha(s)| \leq  \frac{A_1}{8}, \quad
        |\beta(t)-\beta(s)| \leq \frac{A_1}{8}
    \end{eqnarray*}
    for every $t, s \in [0,T_1]$ such that $|t-s| \leq \delta$. Let $n_0:= \big\lfloor \frac{T_1}{\delta}\big\rfloor +1$. Partition $[0, T_1]$ into $n_0$ intervals of equal length, $I_1,..., I_{n_0}$, where $I_k= \big[\frac{(k-1)T_1}{n_0}, \frac{kT_1}{n_0}\big]$, $1 \leq k \leq n_0$, $k \in \N$. 
    \paragraph{} Now, for every $I_k$, choose $a,b \in \R$ such that 
    \begin{equation*}
        \alpha\bigg(\frac{kT_1}{n_0} \bigg) + \frac{A_1}{4}\leq a <b \leq \beta \bigg( \frac{kT_1}{n_0}\bigg) -\frac{A_1}{4}.
    \end{equation*}
    Then, $\alpha(t)<a<b<\beta(t)$ for every $t \in I_k$.
    
\end{proof}

We use Proposition \ref{prop: 44} and Proposition \ref{prop: subdivide [0,T_1]} to prove uniform convergence of $u$ in subsets of $\Gamma_{h(t, \cdot)}$ that are away from the end points.

\begin{prop}\label{prop: uniform conv u}
    Let $T_1>0$ be as in Theorem \ref{Res: 39,40,41,42,43} (iv) and $I$ be a sub-interval of $[0, T_1]$ as in Proposition \ref{prop: subdivide [0,T_1]}. Let $a, b \in \R$ such that $\alpha(t)< a < b < \beta(t)$ for every $t \in I$. Let $\hat{u}_k (t,\cdot,\cdot)$  minimize the functional $\e(\hat{\alpha}_k(t), \hat{\beta}_k(t), \hat{h}_k(t,\cdot), u)= \int_{\Omega_{\hat{h}_k(t,\cdot)}} W(Eu(t,x,y)) dx dy$. Then, there exists a subsequence of $\{\hat{u}_k(t,\cdot,\cdot)\}_k$ (not relabeled) such that
    \begin{equation}\label{eqn: nabla u_k uniform conv nabla u on graph of h}
        \nabla \hat{u}_k(t,\cdot, h_k(t, \cdot)) \rightarrow \nabla u(t,\cdot, h(t,\cdot)) \text{ uniformly in } [a,b]
    \end{equation}
    for every $t \in I$, where $u(t,\cdot, \cdot)$ is the unique minimizer of (\ref{problem: minimization problem on h(t,.)}).
\end{prop}
\begin{remark}
    Note that 
    \begin{equation}\label{eqn: def u_k hat}
        \hat{u}_k(t,x,y)= u^i_k(x,y)
    \end{equation}for $(x,y) \in \Omega_{\hat{h}_k(t,\cdot)}$, $t \in (t^{i-1}_k, t^i_k]$, for $i,k \in \N$, $1\leq i \leq k$, where $u^i_k$ is defined in section (\ref{sec: Discretization}).
\end{remark}
Before we start the proof, we state a lemma that follows from Theorem 5.19 in \cite{adams2003sobolev} and a reflection argument.
\begin{lem}\label{lemma: Exc 13.3 AFCSS}
    Let $1\leq p \leq \infty$ and $u \in W^{2,p}(\R^2_-)$, where $\R^2_-:= \R \times (0, -\infty)$. Then, there exists $c_1, c_2 \in \R$ such that the function
    \begin{equation}
        v(x):= \begin{cases}
             c_1 u(x,y)+c_2u(x, -2y) \text{ if } y>0,\\
            u(x,y) \hspace{3 cm} \text{if }y < 0
        \end{cases}
    \end{equation}
    is well defined and belongs to $W^{2,p}(\R^2)$. Further, for every $k=0,1,2$,
    \begin{equation}
        \|\nabla^k v \|_{L^p(\R^2)} \leq c \|\nabla^k u\|_{L^p(\R^2_-)} 
    \end{equation}
    for some constant $c=c(p)>0$.
\end{lem}

\begin{proof}[Proof of Proposition \ref{prop: uniform conv u}] By (\ref{eqn: uniform conv_alpha, beta}) and (\ref{eqn: a,b between alpha beta functions}), there exists $k_3 \geq k_2 \in \N$ such that $ \alpha_k(t)< a < b < \beta_k(t)$ for every $k\geq k_3$, for every $t \in I$. Now, fix $a_0, b_0 \in \R$ such that $\alpha(t)<  a_0< a < b <b_0 < \beta(t)$ and $\alpha_k(t)< a_0< a < b < b_0 < \beta_k(t)$ for all $k \geq k_3$, for all $t \in I$. Then, in particular, we have $h(t, x)>0$ and $h_k(t, x) >0$ for all $x \in [a_0,b_0]$ and $t \in I$. Fix $t \in I$ and choose $0<\varepsilon < \min\{ \frac{1}{8}\min\{h(t,x): x \in [a_0,b_0]\}, \frac{1}{2}(a-a_0), \frac{1}{2}(b_0-b)\} $.
\paragraph{} Now, by (\ref{eqn: uniform conv of hat h_k in R}), we choose $k_4 \geq k_3$, $k_4 \in \N$ such that $\|\hat{h}_k(t,\cdot)-h(t,\cdot)\|_{L^{\infty}(\R)} < \varepsilon$ for every $k \geq k_4$. That is, $\hat{h}_k (t,x) > h(t,x)-\varepsilon >7 \varepsilon >0$, for every $x \in [a_0,b_0]$, for every $k \geq k_4$. Suppose $z_k:=(x_k, \hat{h}_k(t, x_k)) \in \Gamma_{\hat{h}_k(t,\cdot)}$ such that $a_0< x_k - 2 \varepsilon < x_k +2\varepsilon< b_0$. Note that $2\varepsilon < \frac{1}{2}\hat{h}_k(t, x_k)$ by the choice of $\varepsilon$. Then, by Proposition 8.9 in \cite{fusco2012equilibrium}, we have that there exists $C_k$, depending only on $\varepsilon$, the Lam{\'e} coefficients $\lambda, \mu$ and $C^{1, \alpha}$ norm of $\hat{h}_k(t,\cdot)$ in $[x_k - 2 \varepsilon, x_k + 2 \varepsilon]$, such that 
\begin{equation}
    \|\nabla \hat{u}_k(t,\cdot,\cdot)\|_{C^{0, \alpha}(\overline{\Omega}_{\hat{h}_k(t,\cdot)} \cap \overline{B}_{\varepsilon}(z_k))} \leq C_k,
\end{equation}
where $B_\varep(z_k)$ is the ball of radius $\varep$ around $z_k$ in $\R^2$. Since $\hat{h}_k(t,\cdot)$ is Lipschitz continuous with Lipschitz constants less than or equal to $L_0$ for every $k$, there exists $N >0$ independent of $k$ such that for every $k \geq k_4$ there exists $z_{k,1}, z_{k,2},..., z_{k,N} \in \Gamma_{\hat{h}_k(t,\cdot)}$ where 
\begin{eqnarray*}
    z_{k,i}:= (x_{k,i}, \hat{h}_k(t, x_{k,i})), \quad
    \bigcup_{i=1}^{N} [x_{k,i}- \varepsilon, x_{k,i}+ \varepsilon]=[a,b].
\end{eqnarray*}
Note that we used the choice of $ \varepsilon$ to obtain the $x_{k,i}$ that satisfy the second condition above. Hence, as $\hat{h}_k(t,\cdot) \rightarrow h(t,\cdot)$ in $C^{1, \alpha}([a_0,b_0])$ by Theorem \ref{Thm: Piovano 3.2.2}, we obtain $C_0>0$ independent of $k$ such that for every $k \geq k_4$,
\begin{equation}\label{eqn:elliptic regularity union of balls}
     \|\nabla \hat{u}_k(t,\cdot,\cdot)\|_{C^{0, \alpha}(\bigcup_{i=1}^N(\overline{\Omega}_{\hat{h}_k(t,\cdot)} \cap \overline{B}_{\varepsilon}(z_{k,i})))} \leq C_0.
\end{equation}
Define $R^k_\varepsilon:= \{(x,y): a \leq x\leq b, \hat{h}_k(t,x)- \varepsilon\leq y \leq \hat{h}_k(t,x)\}$. By choosing a smaller $\varepsilon$ if required, (\ref{eqn:elliptic regularity union of balls}) gives
\begin{equation}
    \|\nabla \hat{u}_k(t,\cdot,\cdot)\|_{C^{0, \alpha}(R^k_{2\varepsilon})} \leq C_0
\end{equation}
for every $k \geq k_4$. 
\paragraph{}Define $G_k(x,y):= (x, y-\hat{h}_k(t, x))$ for $(x,y) \in [a,b] \times \R$. Clearly, for every $k$, $G_k(R^k_{2\varepsilon})= Q_{2\varepsilon}$ where $Q_\varepsilon:=\{(x,y): a\leq x \leq b, -\varepsilon \leq y \leq 0\} \subset \R^2$. Let $v_k(t, x, y):= \hat{u}_k(t, G^{-1}(x,y))$ for $(x,y) \in Q_{2\varepsilon}$. Since $G_k : R^k_{2\varep} \rightarrow Q_{2\varep}$ is a $C^1$ homeomorphism and since the Lipschitz constants of $\hat{h}_k(t,\cdot)$, $h(t,\cdot)$ are uniformly bounded by $L_0$, $v_k(t,\cdot,\cdot) \in C^{1, \alpha}(Q_{2\varep})$ for every $k \geq k_3$ and there exists $C_1>0$ such that for every $k \geq k_4$,
\begin{equation}
    \|\nabla v_k(t,\cdot,\cdot)\|_{C^{0, \alpha}(Q_{2\varep})} \leq C_1.
\end{equation}
Using Lemma \ref{lemma: Exc 13.3 AFCSS}, we have $\overline{v}_k(t,\cdot,\cdot)$ defined on $Q^0_{2\varep}:= \{(x,y): a\leq x\leq b, -2\varep \leq y \leq 2\varep) \subset \R^2$ and $C_2>0$ such that
\begin{equation*}
    \|\nabla \overline{v}_k(t,\cdot,\cdot)\|_{C^{0, \alpha}(Q^0_{2\varep})} \leq C_2
\end{equation*}
for every $k \geq k_4$. Let $S^k_{2\varep} := \{(x,y): a \leq x \leq b, \hat{h}_k(t,x)- 2\varep \leq y \leq \hat{h}_k(t,x)+ 2\varep\} \subset \R^2$. Note that $G_k(Q^0_{2\varep})= S^k_{2\varep}$. Define $\overline{u}_k(t,\cdot,\cdot)$ as 
\begin{equation*}
    \overline{u}_k(t,x,y):= \overline{v}_k(t, G_k(x,y))
\end{equation*}
for $(x, y) \in S^k_{2\varep}$ and $k \geq k_4$. Again, since $G_k$ is a $C^1$ homeomorphism and since the Lipschitz constants of $\hat{h}_k(t,\cdot)$, $h(t,\cdot)$ are uniformly bounded by $L_0$, there exists $C>0$ such that
\begin{equation}
     \|\nabla \overline{u}_k(t,\cdot,\cdot)\|_{C^{0, \alpha}(S^k_{2\varep})} \leq C,
\end{equation}
for every $k \geq k_4$. Now, by the choice of $\varep$ and $k_4$ we have the estimate
\begin{equation}\label{eqn: u_k bar C^1, alpha bound}
     \|\nabla \overline{u}_k(t,\cdot,\cdot)\|_{C^{0, \alpha}(S_{\varep})} \leq C
\end{equation}
for every $k \geq k_3$ where $S_\varep:= \{(x,y): a\leq x \leq b, h(t,x)-\varep \leq h(t,x)+\varep\} $. Hence, by Ascoli-Arzel\`a theorem there exists a subsequence $\{\overline{u}_k(t,\cdot,\cdot)\}_{k\geq k_4}$ (not relabeled) that converges to a limit function $\overline{u}(t,\cdot,\cdot)$ in $C^{1, \alpha}(S_\epsilon)$. By Proposition \ref{prop: 44} it is clear that the restriction of $\overline{u}(t,\cdot,\cdot)$ to $R_\varep:=\{(x,y): a\leq x \leq b, h(t,x) - \varep \leq y \leq h(t,x)\}$ is $u(t,\cdot,\cdot)$, the unique minimizer of (\ref{problem: minimization problem on h(t,.)}) and that (\ref{eqn: nabla u_k uniform conv nabla u on graph of h}) holds.

\end{proof}
Now, we state a theorem that describes the dynamics of $\alpha, \beta$ in time.
\begin{thm}\label{thm: contact points ODE final}
    Under the assumptions (\ref{assumption: alpha_0, beta_0, h_0 are in A_s})-(\ref{assumption: Lip h_0< L_0}), let $T_1$ be as in Theorem \ref{Res: 39,40,41,42,43} (iv). Then, for a.e. $t \in (0, T_1)$, we have 
    \begin{eqnarray}
        \sigma_0 \dot{\alpha}(t)&=& \frac{\gamma}{J(t,\alpha(t))}-\gamma_0+ \nu_0 \frac{h'(t,\alpha(t))}{(J(t,\alpha(t))^2}\bigg(\frac{h''(t,\cdot)}{(J(t,\cdot))^3}\bigg)'(\alpha(t)), \label{eqn: alpha ODE}\\
        \sigma_0 \dot{\beta}(t)&=& -\frac{\gamma}{J(t,\beta(t))}+\gamma_0- \nu_0 \frac{h'(t,\beta(t))}{(J(t,\beta(t))^2}\bigg(\frac{h''(t,\cdot)}{(J(t,\cdot))^3}\bigg)'(\beta(t)), \label{eqn: beta ODE}
    \end{eqnarray}
    where $J(t,x)$ is defined in (\ref{eqn: def J(t,x)}).
\end{thm}
\begin{proof}
Let $t \in [0,T_1]$. By (\ref{eqn: uniform conv_alpha hat, beta hat}), there exists $k_5 \geq k_4$ such that $|\hat{\alpha}_k(t^i_k)-\hat{\alpha}_k(t^{i-1}_k)| \leq \delta_0$ and $|\hat{\beta}_k(t^i_k)-\hat{\beta}_k(t^{i-1}_k)| \leq \delta_0$ for every $k\geq k_5$ and $1\leq i\leq k$, where $k_4$ is as in Proposition \ref{prop: uniform conv u}. Find $1\leq i\leq k$ where $t \in [t^{i-1}_k, t^i_k]$. Define 
\begin{equation}\label{eqn: def J^i_k}
    J^i_k(x):= (1+((\htil^i_k)'(x))^2)^{1/2}
\end{equation}
for $x \in \R$ and $0\leq i \leq k$, $i,k \in \N$.  Now, by (\ref{eqn: alpha ode precursor}) and (\ref{eqn: beta ode precursor}) in Theorem \ref{thm: 16}, 
\begin{eqnarray*}
    \sigma_0 \frac{\alpha^i_k -\alpha^{i-1}_k}{\tau_k}= \frac{\gamma}{J^i_k(\alpha^i_k)}-\gamma_0 +\nu_0 \frac{(h^i_k)'(\alpha^i_k)}{(J^i_k(\alpha^i_k))^2}\bigg(\frac{(h^i_k)''}{(J^i_k)^3}\bigg)'(\alpha^i_k),\\
     \sigma_0 \frac{\beta^i_k -\beta^{i-1}_k}{\tau_k}= -\frac{\gamma}{J^i_k(\beta^i_k)}+\gamma_0 -\nu_0 \frac{(h^i_k)'(\beta^i_k)}{(J^i_k(\beta^i_k))^2}\bigg(\frac{(h^i_k)''}{(J^i_k)^3}\bigg)'(\beta^i_k).
\end{eqnarray*}
Now, we can proceed as in the proof of Theorem 9.18 in \cite{dal2025motion} to obtain the result.

\end{proof}

\begin{remark}\label{rem: intrinsic alpha, beta odes}
    We can express (\ref{eqn: alpha ODE}) and (\ref{eqn: beta ODE}) in terms of the arc length of the curve $h(t,\cdot)$ and $\theta_\alpha(t), \theta_\beta(t)$ which are the oriented angles between the oriented x axis and the tangent to the graph of of $h(t,\cdot)$ at $(\alpha(t), 0), (\beta(t),0)$ respectively. We can write $\theta_\alpha(t), \theta_\beta(t)$ as 
    \begin{eqnarray}
        \theta_\alpha(t) &=& \operatorname*{arc sin} \frac{h'(t, \alpha(t))}{(1+(h'(t,\alpha(t)))^2)^{1/2}}, \label{eqn: def theta_alpha}\\
        \theta_\beta(t) &=& \operatorname*{arc sin} \frac{h'(t, \beta(t))}{(1+(h'(t,\beta(t)))^2)^{1/2}}. \label{eqn: def theta_beta}
    \end{eqnarray}
   Then, 
   \begin{eqnarray*}
       \cos{\theta_\alpha(t)}= \frac{1}{J(t, \alpha(t))} \text{ and }
       \cos{\theta_\beta(t)}= \frac{1}{J(t, \beta(t))}.
   \end{eqnarray*}
Now, let $t \in [0, T_1]$ and $x \in [\alpha(t), \beta(t)]$. Then, define the arc length 
    \begin{equation} \label{eqn: def arc length s}
        s(t,x):= \int_{\alpha(t)}^x \sqrt{1+(h'(t,r))^2} dr.
    \end{equation}
  Fix $t \in [0, T_1]$ and define $x(t,\cdot): [s(t,\alpha(t)), s(t, \beta(t))]\rightarrow [\alpha(t), \beta(t)]$ as the inverse of the increasing function $s(t,\cdot)$. Let $\kappa(t,\cdot): [s(t,\alpha(t)), s(t, \beta(t))] \rightarrow \R$ be the signed curvature of the graph of $h(t,\cdot)$ considered as the function of arc-length, that is, 
  \begin{equation}\label{eqn: def curvature kappa}
      \kappa(t,s):= \frac{h''(t, x(t,s))}{(1+(h'(t, x(t,x)))^2)^{3/2}}.
  \end{equation}
  Then, $\kappa(t, s(t,x))= \frac{h''(t,x)}{J(t,x)^3}$. Differentiating with respect to $x$, and using chain rule, we have 
  \begin{eqnarray}\label{eqn: kappa, s chain rule}
    \pa_x \bigg( \frac{h''(t,x)}{J(t,x)^3}\bigg) = \pa_s \kappa(t,s)  \pa_x s(t,x) = \pa_s \kappa(t,s) J(t,x).
  \end{eqnarray}
 Therefore, (\ref{eqn: alpha ODE}) and (\ref{eqn: beta ODE}) can be written in the form
 \begin{eqnarray*}
      \sigma_0 \dot{\alpha}(t)&=& \gamma \cos{\theta_\alpha(t)}-\gamma_0 + \nu_0 \pa_s \kappa(t,s(t, \alpha(t))) \sin{\theta_\alpha(t)},\\
      \sigma_0 \dot{\beta}(t)&=& -\gamma\cos{\theta_\beta(t)}+\gamma_0-\nu_0 \pa_s \kappa(t,s(t, \beta(t))) \sin{\theta_\beta(t)}.
 \end{eqnarray*}
\end{remark}

\begin{prop}
    For a.e. $t \in [0, \infty)$, there exists $\dot{h}(t,\cdot) \in L^2((\alpha(t), \beta(t)))$, such that 
\begin{equation}\label{eqn: time derivative convergence of h}
    \frac{h(s,\cdot)- h(t,\cdot)}{s-t} \rightarrow \dot{h}(t,\cdot) 
\end{equation}
in $L^2((a,b))$ as $s \rightarrow t$ for every $\alpha(t)< a < b < \beta(t)$.
\end{prop}

\begin{proof}
    Fix $T>0$. By Proposition \ref{prop: 32}, 
\begin{equation}
    h_* \in H^1((0,T);L^2(\R)).
\end{equation}
Since difference quotients of an $H^1$ function converges to it's weak derivative in $L^2$, we have 
\begin{equation}\label{eqn: converging difference quotients}
    \lim_{s \rightarrow t} \frac{h_*(s,\cdot)-h_*(t,\cdot)}{s-t} = \dot{h}_*(t,\cdot) \text{ in } L^2(\R).
\end{equation}
In particular, if $t \in [0,\infty)$ satisfies (\ref{eqn: converging difference quotients}) and $\alpha(t) < a<b< \beta(t)$, since $\alpha, \beta$ are continuous by Proposition \ref{prop: 30}, we have $\alpha(s) < a<b< \beta(s)$ for all $s$ sufficiently close to $t$. Hence, (\ref{eqn: time derivative convergence of h}) follows.

\end{proof}

\begin{thm}\label{thm: h PDE final}
    Under the assumptions (\ref{assumption: alpha_0, beta_0, h_0 are in A_s})-(\ref{assumption: Lip h_0< L_0}), let $T_1$ be as in Theorem \ref{Res: 39,40,41,42,43}. Then, for a.e. $t \in (0, T_1)$, we have 
    \begin{equation}\label{eqn: final thm PDE}
        \dot{h} = J\bigg[  \gamma \bigg(\frac{h'}{J}\bigg)'-\nu_0 \bigg(\frac{h''}{J^5}\bigg)'' - \frac{5 \nu_0}{2}\bigg(\frac{h'(h'')^2}{J^7}\bigg)'-\overline{W}+m\bigg]
    \end{equation}
    in $\d'((\alpha(t), \beta(t)))$, where $J: [0, \infty) \times \R \rightarrow \R$ is defined in (\ref{eqn: def J(t,x)}), $\overline{W}(t,x):=W(Eu(t,x,h(t,x)))$ for $t \in [0,T_1]$ and $x \in (\alpha(t), \beta(t))$ and the Lagrange multiplier $m \in L^\infty ((0, T_1))$ has the form
    \begin{gather}
        m(t)= \nu_0 \int_\alpha^\beta \frac{h''(t,x)}{(J(t,x)^5}\phi_0''(t,x) dx + \gamma \int_\alpha^\beta \bigg(\frac{h'(t,x)}{J(t,x)}\bigg)'\phi_0'(t,x) dx \nonumber\\ -\frac{5}{2} \nu_0 \int_\alpha^\beta \frac{h'(t,x)(h''(t,x))^2}
{(J(t,x))^7} \phi_0'(t,x) dx 
-\int_\alpha^\beta W(Eu(t,x,h(t,x))) \phi_0(t,x) dx \nonumber \\-\int_\alpha^\beta\frac{1}{J(t,x)}\dot{h}(t,x) \phi_0(t,x) dx, \label{eqn: def m(t)}
    \end{gather}
for $t\geq 0$ and $x \in \R$, where $\phi_0: (0, \infty) \times \R \rightarrow \R$ is defined as in (\ref{eqn: def phi_0(t,x)}), in Corollary \ref{cor: a_0, b_0, phi_0 conv}.

\end{thm}
\begin{proof}
By Proposition \ref{prop: subdivide [0,T_1]}, we can subdivide $[0,T_1]$ into $n_0$ intervals $\{I_j\}_{j=1}^{n_0}$ such that there exists $a,b \in \R$ such that $\alpha(t) <a<b< \beta(t)$ for every $t \in I_j$, for every $j$. It is enough to prove convergence in one such $I_{j_0}=:I$. Let $\delta>0$ be such that $\alpha(t)+3 \delta < a< b< \beta(t)-3 \delta$ for all $t \in I$ and (\ref{eqn: choose delta prop 41}) is satisfied. By (\ref{eqn: uniform conv_alpha, beta}), there exists $k_6 \in \N$, $k_6 \geq k_5$ such that $\alpha_k(t)+ 2 \delta < a<b< \beta_k(t) - 2 \delta$ for all $k \geq k_6$. Now, we use Theorem \ref{Res: 39,40,41,42,43} (iii) to obtain $\zeta>0$ and $k_7 \geq k_6$ such that 
\begin{equation*}
    h_k(t, x) \geq \zeta \text{ for all } t \in I, \, x \in [a-\delta, b+\delta], \text{ for all } k \geq k_7.
\end{equation*}
As $ k \rightarrow \infty$, by (\ref{eqn: uniform convergence of h_k in R}), we have 
\begin{equation*}
     h(t, x) \geq \zeta \text{ for all } t \in I, \, x \in [a-\delta, b+\delta].
\end{equation*}

Fix $\varphi \in C_c^ \infty ((a,b) \times (t_1, t_2))$. Define
\begin{gather}
     V^i_k (x) = \frac{1}{\tau_k}\frac{\htil^i_k(x)- \htil^{i-1}_k(x)}{J^{i-1}_k(x)}
\end{gather}
for $x \in \R$, $1\leq i\leq k, k \in \N$. Recall the definition of $\hat{u}_k$ from (\ref{eqn: def u_k hat}). Now, define $\hat{V}_k$ and $\hat{\overline{W}}_k$ as 
\begin{gather}
    V_k(t,x):= V^i_k (x), \quad \hat{\overline{W}}_k(t,x):= W (E\hat{u}_k (t, x, \hat{h}_k(t,x))),
\end{gather}
for $t \in (t^{i-1}_k, t^i_k]$ and $x \in (\alpha(t), \beta(t))$. Now, for $t \in I$, by (\ref{eqn: EL strong equation}) and (\ref{eqn: def m^i_k}), we have 
\begin{eqnarray*}
    -\gamma \int_a^b \frac{\hat{h}_k'\varphi'}{\hat{J}_k} dx+ \nu_0 \int_a^b \frac{\hat{h}_k'' \varphi''}{\hat{J}_k^5}dx+\frac{5 \nu_0}{2} \int_a^b \frac{\hat{h}_k'(\hat{h}_k'')^2\varphi'}{\hat{J}_k^7} dx + \int_a^b \hat{\overline{W}}_k \varphi dx + \int_a^b V_k\varphi  dx = \hat{m}_k \int_a^b \varphi dx
\end{eqnarray*}
where $\hat{m}_k(t)$ is defined in (\ref{eqn: def a_0k hat, b_0k hat, phi_0,k hat, m_k hat}). Integrating with respect to time on both sides, we obtain
\begin{eqnarray*}
    -\gamma \int_I\int_a^b \frac{\hat{h}_k'\varphi'}{\hat{J}_k} dx dt+ \nu_0 \int_I \int_a^b \frac{\hat{h}_k'' \varphi''}{\hat{J}_k^5}dx dt+\frac{5 \nu_0}{2} \int_I\int_a^b \frac{\hat{h}_k'(\hat{h}_k'')^2\varphi'}{\hat{J}_k^7} dx dt \\
    +\int_I \int_a^b \hat{\overline{W}}_k \varphi dx dt  + \frac{1}{\tau_k}\int_I \int_a^b V_k\varphi  dx dt = \int_I \hat{m}_k \int_a^b \varphi dx dt.
\end{eqnarray*}
\begin{itemize}
    \item \textbf{Weak convergence of $V_k$:}\\
    Note that 
    \begin{equation}\label{eqn: def V_k hat}
        V_k(t, x) = \frac{1}{\hat{J}_{k}(t-\tau_k, x)} \dot{h}_k(t,x)
    \end{equation}
    for $(t^{i-1}_k, t^i_k)$, where $\hat{J}_k(t,x):= J^i_k(x)$ for $t \in (t^{i-1}_k, t^i_k]$ and $x \in \R$. Note that $J^i_k$ is defined in (\ref{eqn: def J^i_k}) for $t \in [0, \infty)$ and $x \in \R$. By (\ref{eqn: time derivative weak conv h_k}), we have 
    \begin{equation}\label{eqn: dot h_k weak conv}
        \dot{h}_k \rightharpoonup \dot{h} \text{ in }L^2 (I; L^2((a,b))).
    \end{equation}
We claim that $\hat{J}_k \rightarrow J$ in $L^\infty (I; C([a,b]))$, where $J(t, x)$ is defined in (\ref{eqn: def J(t,x)}). Note that the function $s \mapsto \frac{1}{(1+s^2)^{1/2}}$ is a 1-Lipschitz function. Then, the claim follows from (\ref{eqn: Piovano_h_k hat L^infty C^1, alpha convergence}) in Theorem \ref{Thm: Piovano 3.2.2}. In turn, by (\ref{eqn: def V_k hat}) and (\ref{eqn: dot h_k weak conv}), we have
\begin{equation}\label{eqn: V_k weak L^2 convergence}
    V_k \rightharpoonup V \text{ in } L^2 (I; L^2 ([a,b])),
\end{equation}
where $V(t,x):= \frac{1}{J(t,x)}\dot{h}(t,x)$ for $t \in I$ and $ x \in [a,b]$.

\item \textbf{Convergence of $\frac{\hat{h}_k'\varphi'}{\hat{J}_k}$:}\\
For $t \in I$ fixed, we have that $\hat{h}_k'(t,\cdot) \rightarrow h'(t,\cdot)$ uniformly on $[a,b]$ by (\ref{eqn: uniform convergence of h_k hat'}). As shown above, $\hat{J}_k \rightarrow J$ in $L^\infty (I; C([a,b]))$. Hence, $  \frac{\hat{h}_k'(t,\cdot)\varphi'(t,\cdot)}{\hat{J}_k(t,\cdot)}  \rightarrow  \frac{h'(t,\cdot)\varphi'(t,\cdot)}{J(t,\cdot)}  $ in $L^1 ((a,b))$ for every $t \in I$. Now, by (\ref{eqn: upper bound on |beta_k-alpha_k|}), there exists $C=C(T_1)>0$ such that
\begin{eqnarray*}
    \bigg| \int_ a^b \frac{\hat{h}_k'(t,x)\varphi'(t,x)}{\hat{J}_k(t,x)} dx\bigg| \leq L_0\|\varphi'\|_{L^\infty (I \times [a,b])} |b-a| \leq C \|\varphi'\|_{L^\infty (I \times [a,b])}.
\end{eqnarray*}
Hence, by Lebesgue's dominated convergence theorem, we have that 
\begin{equation} \label{eqn: h'/J convergence}
    \frac{\hat{h}_k'\varphi'}{\hat{J}_k} \rightarrow \frac{h'\varphi'}{J} \text{ in } L^1 (I; L^1((a,b))).
\end{equation}

\item \textbf{Convergence of $ \frac{\hat{h}_k'' \varphi''}{\hat{J}_k^5}$:}\\
From (\ref{eqn: h_k^''' integral bound}) and (\ref{eqn: h_k ^ (iv) integral bound}), we have that $\{\hat{h}_k''(t,\cdot)\}_{k \geq k_7}$ is uniformly bounded in $W^{2,p_1}((a,b))$ for every $t \in I$. Therefore, using Morrey's theorem, $\{\hat{h}_k''(t,\cdot)\}_{k \geq k_7}$ is uniformly bounded in $C^{1, \alpha }([a,b])$ for some $\alpha \in (0,1)$. Hence, using Ascoli-Arzel\`a theorem and uniqueness of limits, we have that $\hat{h}_k (t,\cdot) \rightarrow h(t,\cdot)$ in $C^2([a,b])$ for every $t \in I$.
\paragraph{} Using the fact that the map $s\mapsto \frac{1}{(1+s^2)^{5/2}}$ is a Lipschitz function, as shown above, we have $\hat{J}_k \rightarrow J$ in $L^\infty (I; C([a,b]))$. Therefore, 
\begin{equation*}
    \int_a^b \frac{\hat{h}_k'' \varphi''}{\hat{J}_k^5} dx \rightarrow \int_a^b \frac{h'' \varphi''}{J^5} dx
\end{equation*}
for every $t \in I$. Further, by (\ref{eqn:uni L^2 bound h^i_k''}) and H{\"o}lder's inequality,
\begin{eqnarray*}
    \bigg|  \int_a^b \frac{\hat{h}_k'' \varphi''}{\hat{J}_k^5} dx\bigg| &\leq&  \|\varphi''\|_{L^\infty (I \times [a,b])} |b-a|^{1/2}M_3^{1/2}\\
    &\leq& C \|\varphi''\|_{L^\infty (I \times [a,b])} 
\end{eqnarray*}
for some $C>0$ independent of $t$, where have we used (\ref{eqn: upper bound on |beta_k-alpha_k|}) again. Hence, by Lebesgue's dominated convergence theorem, we have 
\begin{equation}\label{eqn: h''/J^5 convergence}
    \frac{\hat{h}_k'' \varphi''}{\hat{J}_k^5} \rightarrow \frac{h'' \varphi''}{J^5} \text{ in } L^1(I; L^1((a,b))).
\end{equation}

    \item  \textbf{Convergence of $ \frac{\hat{h}_k'(\hat{h}_k'')^2\varphi'}{\hat{J}_k^7} $:} Following the same arguments as above, we have that 
    \begin{equation}\label{eqn: h'(h'')^2/J^7 convergence}
         \frac{\hat{h}_k'(\hat{h}_k'')^2\varphi'}{\hat{J}_k^7} \rightarrow  \frac{h'(h'')^2\varphi'}{J^7}  \text{ in } L^1(I; L^1((a,b))).
    \end{equation}

    \item \textbf{Convergence of $\hat{\overline{W}}_k$:} We have proved that $\nabla \hat{u}_k(t,\cdot,\hat{h}_k(t,\cdot))$ converges to $\nabla u(t,\cdot,h(t,\cdot))$ uniformly in $[a,b]$ in Proposition \ref{prop: uniform conv u}. Hence, by definition of $W$ in (\ref{eqn: def W}) we have that 
    \begin{equation}\label{eqn: W on graph of h convergence}
        W(E\hat{u}_k(t,\cdot,\hat{h}_k(t,\cdot))) \rightarrow W(Eu(t,\cdot,h(t,\cdot))) \text{ uniformly in } [a,b].
    \end{equation}
Therefore, $\int_a^b \hat{\overline{W}}_k \rightarrow \int_a^b W(Eu(t,x,h(t,x)))dx$ as $k \rightarrow \infty$. In turn, by Lebesgue's dominated convergence theorem, we have 
\begin{equation}\label{eqn: W_k bar hat convergence}
     \hat{\overline{W}}_k \rightarrow \overline{W}  \text{ in } L^1(I; L^1((a,b))).
\end{equation}
where $\overline{W}(t,x):=  W(Eu(t,x,h(t,x)))$ for $x \in (\alpha(t), \beta(t))$ and $t \in [0, T]$.
\item \textbf{Convergence of $\hat{m}_k$:} Recall the definition of $\phi_0: (0, \infty) \times \R \rightarrow \R$ from (\ref{eqn: def phi_0(t,x)}). By (\ref{eqn: phi_0,k hat with two derivatives uniform convergence}), (\ref{eqn: V_k weak L^2 convergence})-(\ref{eqn: h'(h'')^2/J^7 convergence}) and (\ref{eqn: W_k bar hat convergence}), we have 
\begin{equation}\label{eqn: m_k hat convergence}
    \hat{m}_k \rightharpoonup m \text{ in } L^1(I),
\end{equation} where $m$ is defined in (\ref{eqn: def m(t)}). By (\ref{eqn: def psi_0}), (\ref{eqn: phi_0'(t,.)}), (\ref{eqn: phi_0''(t,.)}), (\ref{eqn: L^2 bound on h''}), (\ref{eqn: V_k weak L^2 convergence}) and (\ref{eqn: W_k bar hat convergence}), we have that $m \in L^\infty((0, T_1))$.
\end{itemize}
Combining the results above, we obtain 
\begin{eqnarray*}
    \int_I \int_a^b \bigg(\nu_0 \frac{h'' \varphi''}{J^5}- \gamma \frac{h'\varphi'}{J}+ \frac{5 \nu_0}{2} \frac{h'(h'')^2\varphi'}{J^7}+ \overline{W} \varphi + V \varphi-m \varphi \bigg)  dx dt =0,
\end{eqnarray*}
for every $\varphi \in C_c^\infty (I  \times (a,b))$. Now, integrate by parts the equation above to obtain
\begin{equation} \label{eqn: final weak form of PDE}
    \int_I \int_a^b g \varphi dx =0,
\end{equation}
where $g(t,x):= \nu_0 \bigg(\frac{h''(t,x) }{(J(t,x))^5}\bigg)'' + \gamma \bigg(\frac{h'(t,x)}{J(t,x)}\bigg)'- \frac{5 \nu_0}{2} \bigg(\frac{h'(t,x)(h''(t,x))^2}{(J(t,x))^7}\bigg)'+  \overline{W}(t,x) + V(t,x) -m(t) $ for $(t,x) \in I \times [a,b]$. We claim that $g \in L^1(I \times (a,b))$.
\paragraph{} Notice that
\begin{eqnarray*}
   \bigg(\frac{h''(t,x) }{(J(t,x))^5}\bigg)''&=& \frac{h^{(iv)}(t,x)}{(J(t,x))^5}-15 \frac{h'(t,x)h''(t,x)h'''(t,x)}{(J(t,x))^7}-5 \frac{(h''(t,x))^3}{(J(t,x))^7}\\
   & & +35\frac{(h'(t,x))^2(h''(t,x))^3}{(J(t,x))^9},\\
  \bigg(\frac{h'(t,x)}{J(t,x)}\bigg)' &=& \frac{h''(t,x)}{J(t,x)}- \frac{(h'(t,x))^2h''(t,x)}{(J(t,x))^3}, \text{ and }\\
    \bigg(\frac{h'(t,x)(h''(t,x))^2}{(J(t,x))^7}\bigg)'&=& \frac{(h''(t,x))^3}{(J(t,x))^7}+ 2\frac{h'(t,x)h''(t,x)h'''(t,x)}{(J(t,x))^7}-7\frac{(h'(t,x))^2(h''(t,x))^3}{(J(t,x))^9}
\end{eqnarray*}
for a.e. $t \in I$ and $x \in (a,b)$. Since $J(t,x) \geq 1$ and $h'(t,x) \leq L_0$ for every $t,x$ by (\ref{eqn: Lip h* and h* support}), we obtain $C= C(\nu_0, \gamma, L_0)>0$ such that 
\begin{eqnarray*}
   \bigg| \nu_0\bigg(\frac{h''(t,x) }{(J(t,x))^5}\bigg)''+  \gamma\bigg(\frac{h'(t,x)}{J(t,x)}\bigg)'+ \frac{5\nu_0}{2}\bigg(\frac{h'(t,x)(h''(t,x))^2}{(J(t,x))^7}\bigg)'\bigg| &\leq& C(|h^{(iv)}(t,x)|+  |h''(t,x) h'''(t,x)|\\
   & & + |h''(t,x)|^3+ |h''(t,x)|)
\end{eqnarray*}
for a.e. $t \in I$ and $x \in (a,b)$.

\paragraph{} Now, note that $\int_I \int_a^b |h^{(iv)}(t,x)|dx dt \leq \overline{M}_6$ where $\overline{M}_6>0$ is independent of $I, a, b$ by Corollary \ref{Prop: Regularity of h}. Further, by (\ref{eqn: upper bound on |beta_k-alpha_k|}), (\ref{eqn: uniform conv_alpha, beta}), (\ref{eqn: L^2 bound on h''}) and H\"older's inequality, 
\begin{gather*}
   \int_I \int_a^b |h''(t,x)| dx dt \leq  \int_I |b-a|^{1/2} \bigg(\int_a^b |h''(t,x)|^2 dx \bigg)^{1/2} dt \\
   \leq \int_I (\beta_0-\alpha_0+2(T_1M_1)^{1/2}) M_3^{1/2} dt \leq  (\beta_0-\alpha_0+2(T_1M_1)^{1/2})T_1 M_3^{1/2}< \infty.
\end{gather*}
Further, by Theorem \ref{Thm: 7.40, AFCSS}, for every $t \in I$,
\begin{eqnarray*}
    \|h''(t,\cdot)\|_{L^{\infty}([a,b])} &\leq& C (\|h'(t,\cdot)\|_{L^{\infty}([a,b])}+ \|h^{(iv)}(t,\cdot)\|_{L^1([a,b])})\\
    &\leq& C(L_0+ s(t))
\end{eqnarray*}
where $C>0$ is independent of $t$ and $s(t):= \int_a^b |h^{(iv)}(t,x)| dx \in L^1(I)$, by (\ref{eqn: L^1 bound on h^(iv)}) from Proposition \ref{Prop: Regularity of h}. Note that $C$ can be chosen to be independent of $a, b$ by (\ref{eqn: upper bound on |beta_k-alpha_k|}) and (\ref{eqn: uniform conv_alpha, beta}). Hence, by using (\ref{eqn: L^2 bound on h''}), (\ref{eqn: L^1 bound on h^(iv)}) and H\"older's inequality, 
\begin{gather*}
    \int_I \int_a^b |h''(t,x)|^3 dx \leq  \int_I C (L_0+ s(t)) \bigg(\int_a^b |h''(t,x)|^2 dx\bigg) dt \\
    \leq  \int_I C M_3 (L_0+ s(t)) dt 
    \leq  CM_3 (L_0 T_1 + \overline{M}_6)< \infty.
\end{gather*}
Using Theorem \ref{Thm: 7.40, AFCSS} again, we have $ \|h'''(t,\cdot)\|_{L^{\infty}([a,b])} \leq C(L_0+ s(t))$  for every $t \in I$ where $C>0$ is independent of $t, a, b$ and $s$ is defined as above. Then,
\begin{gather*}
    \int_I \int_a^b  |h''(t,x) h'''(t,x)| dx dt \leq  \int_I  \|h'''(t,\cdot)\|_{L^{\infty}([a,b])} \bigg(\int_a^b |h''(t,x)|dx\bigg) dt\\
    \leq \int_I C M_3^{1/2} (L_0+ s(t)) dt
    \leq CM_3^{1/2} (L_0 T_1 + \overline{M}_6)< \infty.
\end{gather*}

By (\ref{eqn: W_k bar hat convergence}), it is clear that $\overline{W}(t,x) \in L^1 (I \times (a,b))$. By (\ref{eqn: m_k hat convergence}), $m \in L^1(I)$. Further, by (\ref{eqn: upper bound on |beta_k-alpha_k|}), (\ref{eqn: uniform conv_alpha, beta}), (\ref{eqn: V_k weak L^2 convergence}) and H\"older's inequality,
\begin{gather*}
    \int_I \int_a^b |V(t,x)| dx dt  \leq \int_I |b-a|^{1/2} \bigg(\int_a^b |V(t,x)|^2 dx \bigg)^{1/2} dt \\
    \leq C T_1^{1/2}  \bigg(\int_I \int_a^b |V(t,x)|^2 dx dt\bigg)^{1/2} < \infty,
\end{gather*}
where $C$ is independent of $I, a, b$. Hence, we have proved the claim. \paragraph{}Thus, by (\ref{eqn: final weak form of PDE}), we have $g =0 $ a.e. in $I \times (a,b)$. That is,
\begin{equation*}
    V= \gamma \bigg(\frac{h'}{J}\bigg)'-\nu_0 \bigg(\frac{h''}{J^5}\bigg)'' - \frac{5 \nu_0}{2}\bigg(\frac{h'(h'')^2}{J^7}\bigg)'-\overline{W}+m
\end{equation*}
a.e. in $I \times (a,b)$. Since $V= \frac{\dot{h}}{J}$, we have (\ref{eqn: final thm PDE}).

\end{proof}
\begin{remark}
    Recall the definitions of $s, \kappa, x$ from Remark \ref{rem: intrinsic alpha, beta odes}. We consider the normal velocity of the graph of $h(t,\cdot)$ as a function of arc length $s$. For $t \in [0, T_1]$, let $\tilde{V}(t,s)$ denote the normal velocity of the time-dependent curve $\Gamma_{h(t,\cdot)}$ at $x(t,s) \in [\alpha(t), \beta(t)]$. Then, 
    \begin{equation*}
        \tilde{V}(t,s):= \frac{\dot{h}(t, x(t,s))}{(1+(h'(t, x(t,s)))^2)^{1/2}}.
    \end{equation*}
    We know that $\kappa(t,s)= \frac{h''(t, x(t,s))}{(J(t, x(t,s)))^3}= \bigg(\frac{h'(t, x(t,s))}{J(t, x(t,s))}\bigg)'$. Further, using (\ref{eqn: kappa, s chain rule}) we have 
    \begin{equation*}
        \bigg(\frac{h''(t, x(t,s))}{(J(t, x(t,s)))^5}\bigg)'' + \frac{5 }{2}\bigg(\frac{h'(t, x(t,s))(h''(t, x(t,s)))^2}{(J(t, x(t,s)))^7}\bigg)' = \pa_{ss}\kappa(t, x(t,s))+ \frac{1}{2}(\kappa (t, x(t,s)))^3.
    \end{equation*}
    Let $\tilde{W}(t,s):= W(Eu(t, x(t,x), h(t, x(t,s))))$ for every $t \in [0, T_1]$. Hence, by (\ref{eqn: final thm PDE}), we have 
    \begin{equation}\label{eqn: normal velocity intrinsic de rem}
        \tilde{V}(t,\cdot)= \gamma \kappa -\nu_0 \bigg( \pa_{ss}\kappa + \frac{1}{2} \kappa^3 \bigg)-\tilde{W}+m  
    \end{equation}
    in $\d((s(t, \alpha(t)), s(t, \beta(t))))$.
\end{remark}
\vspace{0.2 cm}
Now, we shall prove Theorem \ref{thm: main thm}.
\begin{proof}[Proof of Theorem \ref{thm: main thm}]
 Let $T=T_1$ from Theorem \ref{Res: 39,40,41,42,43} (iv). 
    \begin{enumerate}
        \item[(i)] The claims follow directly from Proposition \ref{prop: 30}.
        \item[(ii)] Property (\ref{eqn: h(0)=h_0 mt}) follows from (\ref{eqn: def h_k}) and (\ref{eqn: uniform convergence of h_k in R}).  Property (\ref{eqn: h=0 at alpha, beta mt}) follows from (\ref{eqn: Lip h* and h* support}), and (\ref{eqn: h positive in (alpha, beta)}) and (\ref{eqn: h' sign at alpha, beta}) implies (\ref{eqn: h>0 in (alpha, beta) mt}) and (\ref{eqn: h' sign at alpha, beta mt}) respectively. Property (\ref{eqn: Lip h strict inequality mt}) follows from (\ref{eqn: Lip constant strict bound on h_*}) and property (\ref{eqn: area under h is constant mt}) follows from (\ref{eqn: integral of h(t,.)=A_0}). (\ref{eqn: h is in H^2 mt}) and (\ref{eqn: h_* regularity mt}) are proved in Proposition \ref{prop:31} and (\ref{eqn: h is in W^4, p_1 mt}) follows from (\ref{eqn: h^iv L^p_1 bound}) in Corollary \ref{Prop: Regularity of h}.
        \item[(iii)] Let $u$ be the function that minimizes (\ref{problem: minimization problem on h(t,.)}). Then, by (\ref{eqn: h is in W^4, p_1 mt}) and Theorem 9.3 in \cite{agmon1964estimates}, (\ref{eqn: u regularity mt}) holds. Moreover, (\ref{eqn: u(t,.,.) bvp mt}) can be obtained as the Euler Lagrange equation of the problem (\ref{problem: minimization problem on h(t,.)}). 
        \item[(iv)] The results (\ref{eqn: alpha ODE final})-(\ref{eqn: V de final}) follow directly from Theorem (\ref{thm: contact points ODE final}) and Theorem (\ref{thm: h PDE final}).
    \end{enumerate}
\end{proof}

\vspace{0.2 cm}
\textbf{Acknowledgments: }This work is part of a thesis to be submitted by the author in partial fulfillment of the requirements for the degree of Doctor of Philosophy in Mathematics under the supervision of Prof. Giovanni Leoni at Carnegie Mellon University. The author is profoundly grateful to her advisor for suggesting the topic and for his guidance throughout the project. The author would like to thank Prof. Irene Fonseca on her valuable remarks on the paper.

\bibliographystyle{siam}
\bibliography{Reference}

\end{document}